\documentclass[5p,twocolumn,authoryear,sort&compress]{elsarticle}
\journal{Annual Reviews in Control}
\pdfoutput=1

\usepackage[latin1]{inputenc}
\usepackage[english]{babel}

\usepackage{setspace}

\usepackage[titletoc,title,page]{appendix}

\usepackage{footnote}

\usepackage{amsmath}
\usepackage{amsfonts}
\usepackage{amssymb}
\usepackage{amsthm}
\usepackage{bbm}		% bold math for numbers (characteristic functions as \mathbbm{1})
\usepackage{mathtools}	% adjust limits for inf sup
\usepackage{bbm}
\usepackage{xcolor}
\usepackage{tikz}
\usepackage[export]{adjustbox}

\usepackage{algorithm}
\usepackage{algpseudocode}
\usepackage{graphicx}
\usepackage{subfigure}
\usepackage[suffix=]{epstopdf}
\usepackage{booktabs}
\usepackage[font=small,labelfont=bf]{caption} 
\usepackage{multirow}
\usepackage{multicol}
\usepackage{adjustbox}
\usepackage{tabularx}

\usepackage{changepage}
\newlength{\offsetpage}
{\end{adjustwidth}}

\usepackage{balance}
\usepackage{makecell}
\usepackage{pifont}% http://ctan.org/pkg/pifont

\newcommand{\PreserveBackslash}[1]{\let\temp=\\#1\let\\=\temp}
\newcolumntype{C}[1]{>{\PreserveBackslash\centering}p{#1}}
\newcolumntype{R}[1]{>{\PreserveBackslash\raggedleft}p{#1}}
\newcolumntype{L}[1]{>{\PreserveBackslash\raggedright}p{#1}}

\usepackage{tikz}
\usetikzlibrary{shapes,matrix,decorations.shapes}
\usetikzlibrary{arrows,decorations.pathmorphing,backgrounds,positioning,fit,matrix}
\usetikzlibrary{decorations.pathreplacing}
\usetikzlibrary{shapes,arrows,chains}
\usetikzlibrary[calc]
\usetikzlibrary[plotmarks]
\usetikzlibrary[patterns,decorations]
\definecolor{thelinkcolor}{RGB}{0,0,150}
\usepackage[
  colorlinks=true,
  hyperindex=true,
  colorlinks=true,
  pagebackref=false,
  plainpages=false,
  pdfpagelabels,
  ]{hyperref}
\makeatletter
    \def\cl@chapter{\@elt {theorem}}
\makeatother

\usepackage[nameinlink]{cleveref}
\crefname{equation}{}{}
\crefname{theorem}{Theorem}{Theorems}
\crefname{corollary}{Corollary}{Corollaries}
\crefname{example}{Example}{Examples}
\crefname{remark}{Remark}{Remarks}
\crefname{lemma}{Lemma}{Lemmas}
\crefname{proposition}{Proposition}{Propositions}
\crefname{figure}{Figure}{Figures}
\crefname{table}{Table}{Tables}
\crefname{section}{Section}{Sections}
\crefname{appendix}{}{}
\Crefname{equation}{}{}
\Crefname{theorem}{Theorem}{Theorems}
\Crefname{corollary}{Corollary}{Corollaries}
\Crefname{example}{Example}{Examples}
\Crefname{lemma}{Lemma}{Lemma}
\Crefname{proposition}{Proposition}{Proposition}
\Crefname{figure}{Figure}{Figures}
\Crefname{table}{Table}{Tables}
\Crefname{section}{Section}{Sections}
\Crefname{appendix}{}{}

\hypersetup{colorlinks=true,
	urlcolor=blue,
	linkcolor=blue,
	citecolor=blue,
	bookmarksdepth=paragraph,
	pdftitle={Chordal and FW decomposition in semidefinite and polynomial optimization},
	pdfauthor={Y. Zheng, G. Fantuzzi, and A. Papachristodoulou},}
	
\usepackage[shortlabels]{enumitem}
\setlist[itemize]{itemindent=0ex,itemsep=-0.5ex,leftmargin=2ex,topsep=5pt}
\setlist[enumerate]{label={\arabic*)},leftmargin=*,nolistsep,noitemsep,topsep=1ex}

\newcommand{\abs}[1]{\left\vert #1 \right\vert}
\newcommand{\tr}{{{\mathsf T}}}

\newcommand{\SDSOS}{\mbox{\it SDSOS}}

\newcommand{\RELU}{\mbox{ReLU}}
\newcommand{\software}[1]{\texttt{#1}}

\newcommand\xqed[1]{\leavevmode\unskip\penalty9999 \hbox{}\nobreak\hfill\quad\hbox{#1}}
\newcommand\markendexample{\xqed{{\small$\blacksquare$}}}

\newcommand{\xmark}{\ding{55}}%

\definecolor{fwFillBlue}{rgb}{0.803921580314636   0.878431379795074   0.968627452850342}
\definecolor{fwFillRed}{rgb}{0.992156863212585   0.917647063732147   0.796078443527222}
\definecolor{fwFillGreen}{rgb}{0.839215695858002   0.909803926944733   0.850980401039124}
\definecolor{fwFillMagenta}{rgb}{0.937254905700684   0.866666674613953   0.866666674613953}
\definecolor{fwLineBlue}{rgb}{0.000000000000000   0.447058826684952   0.741176486015320}
\definecolor{fwLineRed}{rgb}{0.850980401039124   0.325490206480026   0.098039217293262}
\definecolor{fwLineGreen}{rgb}{0.000000000000000   0.498039215803146   0.000000000000000}
\definecolor{fwLineMagenta}{rgb}{0.494117647409439   0.184313729405403   0.556862771511078}
\definecolor{mygreen}{RGB}{77,175,74}
\definecolor{myred}{RGB}{228,26,28}
\definecolor{myblue}{RGB}{55,126,184}
\definecolor{matlabgray}{RGB}{127,127,127}
\definecolor{matlabblue}{RGB}{0,113,188}
\definecolor{matlabred}{RGB}{216,82,24}
\definecolor{matlabgreen}{rgb}{0.4660,0.6740,0.1880}
\definecolor{matlabcyan}{rgb}{0.3010,0.7450,0.9330}   
\definecolor{matlabyellow}{rgb}{0.9290,0.6940,0.1250}
\definecolor{matlaborange}{RGB}{255,153,102}
\definecolor{matlabpurple}{rgb}{0.4940,0.1840,0.5560}
\definecolor{matlabsafered}{RGB}{215,25,28}
\definecolor{matlabsafegreen}{RGB}{171,221,164}
\definecolor{lightbrown}{RGB}{149,99,99}
\definecolor{darkbrown}{RGB}{82,48,48}
\newcommand\solidrule[1][10pt]{\rule[0.5ex]{#1}{1.5pt}}

\newcommand{\mycross}[1]{%
	\protect\begin{tikzpicture}%
	\protect\draw[thick,color=#1] (0,0) -- (1ex,1ex);
	\protect\draw[thick,color=#1] (0,1ex) -- (1ex,0);
	\protect\end{tikzpicture}%
}
\newcommand{\mycrossLR}[2]{%
	\protect\begin{tikzpicture}%
	\protect\draw[thick,color=#1] (0,0) -- (0.5ex,0.5ex);
	\protect\draw[thick,color=#1] (0.5ex,0.5ex) -- (0,1ex);
	\protect\draw[thick,color=#2] (1ex,0) -- (0.5ex,0.5ex);
	\protect\draw[thick,color=#2] (0.5ex,0.5ex) -- (1ex,1ex);
	\protect\end{tikzpicture}%
}

\DeclareMathOperator{\supp}{supp}							% support
\DeclareMathOperator*{\nnz}{nnz}                            % nonzeros
\DeclareMathOperator*{\argmin}{argmin}                      % argmin
\DeclareMathOperator*{\New}{New}                            % Newton polytope
\DeclareMathOperator*{\csp}{csp}                            % correlative sparsity graph
\DeclareMathOperator*{\var}{var}                            % correlative sparsity graph

\DeclareMathOperator{\trace}{trace}

\newtheorem{theorem}{Theorem}

\newtheorem{proposition}{Proposition}
\newtheorem{definition}{Definition}
\newtheorem{corollary}{Corollary}
\theoremstyle{definition}
\newtheorem{example}{Example}
\newtheorem{remark}{Remark}
\newtheorem{assumption}{Assumption}

\numberwithin{equation}{section}
\numberwithin{figure}{section}
\numberwithin{theorem}{section}
\numberwithin{corollary}{section}
\numberwithin{lemma}{section}
\numberwithin{proposition}{section}
\numberwithin{definition}{section}
\numberwithin{remark}{section}
\numberwithin{example}{section}

\makeatletter
\newcommand{\subalign}[1]{%
	\vcenter{%
		\Let@ \restore@math@cr \default@tag
		\baselineskip\fontdimen10 \scriptfont\tw@
		\advance\baselineskip\fontdimen12 \scriptfont\tw@
		\lineskip\thr@@\fontdimen8 \scriptfont\thr@@
		\lineskiplimit\lineskip
		\ialign{\hfil$\m@th\scriptstyle##$&$\m@th\scriptstyle{}##$\crcr
			#1\crcr
		}%
	}
}
\newcommand{\pushright}[1]{\ifmeasuring@#1\else\omit\hfill$\displaystyle#1$\fi\ignorespaces}
\newcommand*\fsize{\dimexpr\f@size pt\relax}
\renewcommand\p@subfigure{\thefigure}
\makeatother            % definitions, newcommands etc

\biboptions{authoryear}

\begin{document}

\begin{frontmatter}
  \title{Chordal and factor-width decompositions for scalable semidefinite \\ and polynomial optimization}
  \author[ucsd]{Yang Zheng} \ead{zhengy@eng.ucsd.edu}
  \author[imperial]{Giovanni Fantuzzi} \ead{giovanni.fantuzzi10@imperial.ac.uk}
  \author[oxford]{Antonis Papachristodoulou} \ead{antonis@eng.ox.ac.uk}
  \address[ucsd]{Department of Electrical and Computer Engineering, University of California San Diego, CA 92093.}
  \address[imperial]{Department of Aeronautics, Imperial College London, London, SW7 2AZ, UK.}
  \address[oxford]{Department of Engineering Science, University of Oxford, Parks Road, Oxford OX1 3PJ, U.K.}
  
  \begin{abstract}
    Chordal and factor-width decomposition methods for semidefinite programming and polynomial optimization have recently enabled the analysis and control of large-scale linear systems and medium-scale nonlinear systems.     Chordal decomposition exploits the sparsity of semidefinite matrices in a semidefinite program (SDP), in order to formulate an equivalent SDP with smaller semidefinite constraints that can be solved more efficiently. Factor-width decompositions, instead, relax or strengthen SDPs with dense semidefinite matrices into more tractable problems, trading feasibility or optimality for lower computational complexity. This article reviews recent advances in large-scale semidefinite and polynomial optimization enabled by these two types of decomposition, highlighting connections and differences between them. We also demonstrate that chordal and factor-width decompositions allow for significant computational savings on a range of classical problems from control theory, and on more recent problems from machine learning. Finally, we outline possible directions for future research that have the potential to facilitate the efficient optimization-based study of increasingly complex large-scale dynamical systems.  
  \end{abstract}


\end{frontmatter}

\tableofcontents

%%%%%%%%%%%%%%%%%%%%%%%%%%%%%%%%%%%%%%%%%%%%%%%%%%%%
%  Main text
%%%%%%%%%%%%%%%%%%%%%%%%%%%%%%%%%%%%%%%%%%%%%%%%%%%%
%%%%%%%%%%%%%%%%%%%%%%%%%%%%%%%%%%%%%%%%%%%%%%%%%%%%
\section{Introduction}
\label{section:intro}

\begin{figure*}
\setlength{\abovecaptionskip}{0pt}
\setlength{\belowcaptionskip}{0pt}
\centering
     \subfigure[]{
     \raisebox{0.5 em}{\begingroup
\tikzstyle{vertex}=[circle,draw=black,thick, inner sep=0pt,minimum size=5mm,node distance=2cm,font=\footnotesize,fill=white]
\begin{tikzpicture}
    \node (1) at (-1,0) [vertex] {1};
	\node (2) at (-1,1.2) [vertex] {2};
	\node (3) at (0,0) [vertex] {3};
	\node (4) at (1,0) [vertex] {4};
	\node (5) at (1,1.2) [vertex] {5};
	\draw[thick] (1) -- (2);
    \draw[thick] (2) -- (3);
    \draw[thick] (1) -- (3);
    \draw[thick] (3) -- (4);
    \draw[thick] (4) -- (5);
    \draw[thick] (3) -- (5);
    %\node at (0,-0.85) {(a)};
	\begin{pgfonlayer}{background}
    \filldraw[fill=mygreen,draw=none] (1.center) -- (2.center) -- (3.center)  -- cycle;
    \filldraw[fill=myblue,draw=none] (3.center) -- (4.center) -- (5.center)  -- cycle;
    \end{pgfonlayer}
	\end{tikzpicture}
	\endgroup}
   }
   \subfigure[]{
   \begingroup % keep the change local
	\setlength\arraycolsep{2pt}
	\def\arraystretch{0.75}
	\begin{tikzpicture}
	\node at (0,0) {\small$\begin{bmatrix}
		\mycross{mygreen} & \mycross{mygreen} & \mycross{mygreen} \\
		\mycross{mygreen} & \mycross{mygreen} & \mycross{mygreen} \\
		\mycross{mygreen} & \mycross{mygreen} & \mycrossLR{mygreen}{myblue} &\mycross{myblue} & \mycross{myblue} \\
		&&\mycross{myblue} & \mycross{myblue} & \mycross{myblue} \\
		&&\mycross{myblue} & \mycross{myblue} & \mycross{myblue}
		\end{bmatrix} = \begin{bmatrix}
		\mycross{mygreen} & \mycross{mygreen} & \mycross{mygreen} \\
		\mycross{mygreen} & \mycross{mygreen} & \mycross{mygreen} \\
		\mycross{mygreen} & \mycross{mygreen} & \mycross{mygreen} &  &  \\
		& & & &  \\
		& & & &
		\end{bmatrix} + \begin{bmatrix}
		 &  &  \\
		 &  &  \\
		 &  & \mycross{myblue} &\mycross{myblue} & \mycross{myblue} \\
		&&\mycross{myblue} & \mycross{myblue} & \mycross{myblue} \\
		&&\mycross{myblue} & \mycross{myblue} & \mycross{myblue}
		\end{bmatrix}$};
	\node at (-2,-1) {\footnotesize $X\succeq 0$};
	\node at (0.2,-1) {\footnotesize $Y_1\succeq 0$};
	\node at (2.2,-1) {\footnotesize $Y_2\succeq 0$};
	\end{tikzpicture}
	\endgroup
   }
   \subfigure[]{
   \begingroup % keep the change local
	\setlength\arraycolsep{2pt}
	\def\arraystretch{0.75}
	\begin{tikzpicture}
	\node at (0,0) {\small$\begin{bmatrix}
		\mycross{mygreen} & \mycross{mygreen} & \mycross{mygreen} \\
		\mycross{mygreen} & \mycross{mygreen} & \mycross{mygreen} \\
		\mycross{mygreen} & \mycross{mygreen} & \mycrossLR{mygreen}{myblue} &\mycross{myblue} & \mycross{myblue} \\
		&&\mycross{myblue} & \mycross{myblue} & \mycross{myblue} \\
		&&\mycross{myblue} & \mycross{myblue} & \mycross{myblue}
		\end{bmatrix}\in \mathbb{S}^n_+(\mathcal{E},?)
		\Leftrightarrow \begin{bmatrix}
		\mycross{mygreen} & \mycross{mygreen} & \mycross{mygreen} \\
		\mycross{mygreen} & \mycross{mygreen} & \mycross{mygreen} \\
		\mycross{mygreen} & \mycross{mygreen} & \mycross{mygreen} &  &  \\
		& & & &  \\
		& & & &
		\end{bmatrix},
		\begin{bmatrix}
		 &  &  \\
		 &  &  \\
		 &  & \mycross{myblue} &\mycross{myblue} & \mycross{myblue} \\
		&&\mycross{myblue} & \mycross{myblue} & \mycross{myblue} \\
		&&\mycross{myblue} & \mycross{myblue} & \mycross{myblue}
		\end{bmatrix}$};
			\node at (-2.6,-1) {\footnotesize $X$};
	\node at (1.1,-1) {\footnotesize $X_1\succeq 0$};
	\node at (2.85,-1) {\footnotesize $X_2\succeq 0$};
	\end{tikzpicture}
	\endgroup
   }
    \caption{Illustration of chordal decomposition, where $\succeq$ denotes positive semidefiniteness and $X \in \mathbb{S}^n_+(\mathcal{E},?)$ is a positive semidefinite completion constraint (see \Cref{subsection:SparseMatrix} for a precise definition). \textit{(a)} A chordal graph with five vertices, six edges, and two maximal cliques (complete connected subgraphs), $\mathcal{C}_1=\{1,2,3\}$ and $\mathcal{C}_2=\{3,4,5\}$; \textit{(b)} Chordal decomposition of a semidefinite constraint on a sparse matrix $X$ into smaller positive semidefinite constraints on matrices $Y_1, Y_2$ with nonzero entries indexed by the cliques $\mathcal{C}_1$ and $\mathcal{C}_2$; \textit{(c)} Chordal decomposition of a positive semidefinite completion constraint on a sparse matrix $X$ with smaller positive semidefinite constraints on its principal submatrices $X_1$ and $X_2$ indexed by the cliques $\mathcal{C}_1$ and $\mathcal{C}_2$.}
    %Chordal decomposition enables us to more efficiently handle optimization problems that involve sparse positive semidefinite constraints.
    \label{fig:chordal_decomposition_overall}
\end{figure*}
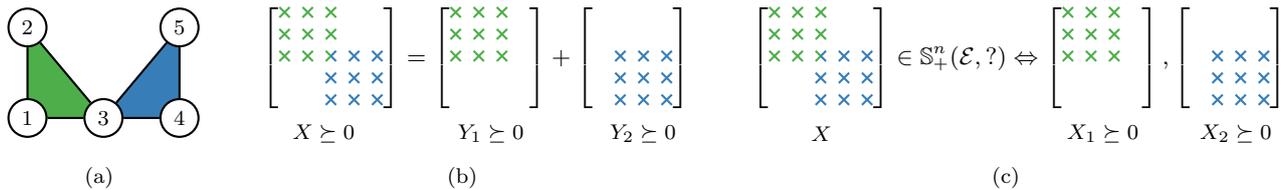

\noindent
The design of innovative technology capable to address the challenges of the 21st century relies on the ability to analyze, predict, and control large-scale complex systems, which are typically nonlinear and may interact over networks~\citep{murray2003future,ASTROM2014control}. Convex optimization is one of the key tools for achieving these goals, because many questions related to the stability and operational safety of dynamical systems, the synthesis of optimal control policies, and the certification of robust performance can be posed as (or relaxed into) convex optimization problems. Very often, these take the form of semidefinite programs (SDPs)---linear optimization problems with positive semidefinite matrix variables.

For linear systems, well-known methods based on linear matrix inequalities (LMIs) enable one to tackle a wide range of problems, including the study of stability, reachability, input-to-state and input-to-output properties, and the design of optimal and robust control strategies~\citep{boyd1994linear,zhou1996robust,kailath1980linear}. Methods based on LMIs have been successfully applied across a broad spectrum of applications, including automotive applications~\citep{rajamani2011vehicle}, flight control~\citep{giulietti2000autonomous}, power grids~\citep{riverso2014plug,sadabadi2016plug}, and traffic systems~\citep{ploeg2013controller,li2017dynamical,zheng2020smoothing}. For nonlinear systems with polynomial dynamics, SDP relaxations based on sum-of-squares polynomials (or, equivalently, moment sequences) enable
stability analysis~\citep{parrilo2000structured,AndP15,PeeP12,HenG05}, the estimation of regions of attractions~\citep{valmorbida2017region,Korda2013roa,topcu2009robust,henrion2014convex,chesi2011domain} and reachable sets~\citep{Jones2019reachable_sets,Magron2019reachable_sets}, safety verification~\citep{PraJP07,Miller2021peak}, analysis of extreme or average behaviour \citep{fantuzzi2016bounds,Fantuzzi2020siads,Korda2021invariant_measures,Kuntz2016,Goluskin2020attractors}, and optimal control~\citep{prajna2004nonlinear,henrion2006convergent,Lasserre2008,Han2018control_om,Majumdar2014occupation_measures,Lasagna2016sos}.

A widespread view since the 1990s is that, once a control problem is reformulated as an SDP or relaxed into one, then the problem is effectively solved~\citep{ParriloLall2003semidefinite,boyd1994linear}. In today's world of large-scale, complex systems, however, this is no longer true, and the formulation of SDPs that can be solved in practice requires further thought. This is because, even though SDPs can theoretically be solved using algorithms with polynomial-time complexity \citep{vandenberghe1996semidefinite,ye2011interior,nesterov1994interior,nesterov2003introductory,nemirovski2006advances}, the very-large-scale SDPs encountered in real-life applications require prohibitively large computational resources in practice.
One particular bottleneck is the complexity of handling large semidefinite constraints; for instance, each iteration of classical interior-point algorithms requires  $\mathcal{O}(n^3m + n^2m^2 + m^3)$ time and $\mathcal{O}(n^2 + m^2)$ memory~\citep[Section 4.3.3]{nesterov2003introductory}, where $n$ is the size of semidefinite constraint and $m$ is the number of equality constraints. The majority of established general-purpose SDP solvers currently available, therefore, cannot handle large problems (e.g., with $n$ larger than a few hundreds and $m$ larger than a few thousands) on a regular computer. Consequently, the application of SDP-based frameworks for analysis and control is currently limited to medium-scale linear systems and small-scale nonlinear ones. 

Overcoming these scalability issues is a problem that has received much attention in recent years~\citep{majumdar2020recent,de2010exploiting,vandenberghe2015chordal,ahmadi2017improving}, and significant progress has been made through a number of different approaches. Most of them are related by a simple, yet powerful, underlying idea: \textit{decompose a large positive semidefinite matrix $X$} %(denoted by $X\in \mathbb{S}^n_+$)
\textit{as a sum of structured ones, for which it is easier to impose positivity}.

One type of structured decomposition considers sums of low-rank matrices \citep{Burer2002,burer2003nonlinear,burer2005local,Burer2006}. Specifically, one writes $X=\sum_{i=1}^t v_i v_i^\tr$ for some vectors $v_1,\ldots,v_t\in \mathbb{R}^n$, where $t \leq n$ is a parameter to be chosen, and optimizes over the choice of such vectors. Such a decomposition is guaranteed to exist for a properly chosen $t$, and there are explicit lower bounds on this parameter ensuring that the global minimum of the decomposed problem coincides with that of the original SDP \citep{Pataki1998}. However, while low-rank decomposition can bring considerable performance gains on large SDPs, it transforms a convex problem into a nonconvex one. Solution algorithms for the latter cannot be guaranteed to converge to the global minimum unless the original SDP is sufficiently ``smooth'' and $t$ is large enough \citep{Waldspurger2020,Boumal2020}.

A second type of structured decomposition, which we focus on in this paper, considers sums of sparse matrices. In this case, one writes $X = \sum_{i=1}^t Y_i$ for positive semidefinite matrices $Y_1,\ldots,Y_t$ that are nonzero only on a certain (and, ideally, small) principal submatrix. The choice of these principal submatrices is crucial in determining the particular type of matrix decomposition, as well as its properties. Two common selection strategies distinguish whether the original matrix $X$ is dense or sparse.

If $X$ is sparse, the principal submatrices are usually indexed by the maximal cliques of the sparsity graph of $X$; these notions will be defined precisely in \Cref{section:chordal-graphs-matrix-decomposition}, but are illustrated in \Cref{fig:chordal_decomposition_overall}. When the sparsity graph is \textit{chordal}, meaning that all cycles of length larger than three have an edge between nonconsecutive vertices, the existence of a clique-based decomposition is guaranteed \citep{agler1988positive,kakimura2010direct,griewank1984existence}. One can therefore replace the optimization of the large matrix $X$ with the optimization of the matrices $Y_1,\ldots,Y_t$ without any loss of generality. Together with a dual result on the existence of positive semidefinite matrix completions~\citep{grone1984positive}, this \textit{chordal decomposition} strategy enables one to significantly reduce the computational complexity of SDPs involving sparse positive semidefinite matrices \citep{fukuda2001exploiting,nakata2003exploiting,kim2011exploiting,vandenberghe2015chordal}. 

When $X$ is dense, instead, each matrix $Y_i$ in the decomposition $X = \sum_{i=1}^t Y_i$ is chosen to be nonzero only on one of the $t=\binom{n}{k}$ possible $k\times k$ principal submatrices of $X$, where the parameter $k \geq 2$ is specified \textit{a priori}. This type of decomposition leads to \textit{factor-width-$k$} inner approximations of the positive semidefinite cone \citep{boman2005factor}, which are conservative but improve as $k$ is increased. When $k \ll n$, optimizing over the matrices $Y_1,\ldots,Y_t$, rather than over the original dense matrix $X$, leads to SDPs with small positive semidefinite cones, which can often be handled efficiently. In the extreme case $k = 2$, one obtains a second-order cone program, for which scalable algorithms exist~\citep{alizadeh2003second}.

This paper offers a comprehensive review of chordal and factor-width-$k$ decomposition methods, as well as of their application to large-scale semidefinite programming and polynomial optimization. Our goal is to introduce practitioners in control theory to the latest advances in these fields, which over the last decade or so have increased the scale of systems for which optimization-based frameworks for analysis and control can be implemented at a reasonable cost. Examples of problems that can now be handled efficiently include the analysis and synthesis of large-scale linear networked systems \citep{mason2014chordal,ZKSP2018scalable,ZMP2018Scalable,andersen2014robust}, the stability analysis and the approximation of regions of attraction for sparse nonlinear systems \citep{Schlosser2020,Tacchi2019,zheng2018sparse,ahmadi2019dsos}, optimal power flow in power grids \citep{andersen2014reduced,jabr2011exploiting,molzahn2013implementation}, and numerous problems in machine learning \citep{BPLZ2021neural,NewP21,kim2009exploiting,dahl2008covariance,zhang2018large,latorre2020lipschitz,chen2020semialgebraic}.
We hope that knowledge of the advanced optimization techniques discussed here can assist control theorists in developing efficient modelling frameworks that can be applied much more widely and, crucially, to increasingly complex large-scale systems.

\subsection{Outline}
\noindent
After introducing relevant graph-theoretic notions in \cref{section:chordal-graphs-matrix-decomposition}, we discuss chordal decomposition methods for general SDPs in \cref{section:sparse-SDPs}. \Cref{section:polynomial_optimization} looks at decomposition methods for sparse polynomial optimization problems, which arise when relaxing analysis and control problems for nonlinear systems. Factor-width-$k$ decompositions for dense matrices are discussed in \cref{section:factor-width-two}. \Cref{section:applications} presents examples of how matrix decomposition can be applied to some classical control problems and to some recent problems in machine learning. \Cref{section:conclusion} draws conclusions and outlines possible directions for future research.

\subsection{Basic notation}
\noindent
Mathematical symbols are defined as necessary in each of the following sections, but we summarize common notation here. The $m$-dimensional Euclidean space, the vector space of $n \times n$ real symmetric matrices, and the cone of $n\times n$ positive semidefinite symmetric matrices are denoted, respectively, by $\mathbb{R}^m$, $\mathbb{S}^n$, and $\mathbb{S}^n_+$. Angled brackets are used to denote the inner product in any of these spaces; in particular, $\langle x, y\rangle = x^\tr y$ when $x, y \in \mathbb{R}^m$ and $\langle X, Y \rangle = \trace(XY)$ when $X, Y \in \mathbb{S}^n$. We often write $X \succeq 0$ instead of $X\in \mathbb{S}^n_+$ when the matrix size is clear from the context or is unimportant, and write $X \succ 0$ if $X$ is strictly positive definite.
%%%%%%%%%%%%%%%%%%%%%%%%%%%%%%%%%%%%%%%%%%%%%%%%%%%%
\section{Chordal graphs and matrix decomposition}
\label{section:chordal-graphs-matrix-decomposition}
\noindent
This section reviews chordal graphs and their applications to sparse matrix decomposition. Matrix decomposition is central to many sparsity-exploiting techniques for semidefinite and polynomial optimization. Detailed introductions to chordal graphs can be found in the surveys by \cite{blair1993introduction} and \cite{rose1970triangulated}, and in the monographs by \cite{vandenberghe2015chordal} and \cite{golumbic2004algorithmic}.
We first introduce some graph-theoretic notions in~\Cref{subsection:chordal_graphs}, and then given an overview of classical matrix decomposition and completion results in~\Cref{subsection:SparseMatrix}. Extensions to sparse block-partitioned matrices are discussed in~\Cref{subsection:block-partition-matrices}.

\subsection{Chordal graphs} \label{subsection:chordal_graphs}
\noindent
A graph $\mathcal{G}(\mathcal{V},\mathcal{E})$ is defined by a set of vertices $\mathcal{V}=\{1,2,\ldots,n\}$ and a set of edges $\mathcal{E} \subseteq \mathcal{V} \times \mathcal{V}$. %We only consider graphs with no self-loops, i.e. $(i,i) \notin \mathcal{E}$.
%We allow
A graph $\mathcal{G}$ is undirected if $(v_i,v_j) \in \mathcal{E}$ implies that $(v_j,v_i) \in \mathcal{E}$. A path in $\mathcal{G}(\mathcal{V},\mathcal{E})$ is a sequence of edges that connect a sequence of distinct vertices. A graph is \emph{connected} if there is a path between
any two vertices, and {\it complete} if any two vertices are connected by an edge, i.e.,
%$ (i,j) \in \mathcal{E}, \forall\, i,j \in \mathcal{V}, i \neq j$.
$\mathcal{E} = \mathcal{V} \times \mathcal{V}$.
The subgraph induced by a subset of vertices $\mathcal{W} \subset \mathcal{V}$ is the undirected graph with vertices $\mathcal{W}$ and edges $\mathcal{E} \cap (\mathcal{W} \times \mathcal{W})$. A subset of vertices $\mathcal{C}\subseteq \mathcal{V} $ is called a \emph{clique} if the subgraph induced by $\mathcal{C}$ is complete.  If $\mathcal{C}$ is not contained in any other clique, it is a \emph{maximal clique}. The number of vertices in $\mathcal{C}$ is denoted by $\vert \mathcal{C} \vert$. %Observe that any connected graph has at most $n-1$ maximal cliques.

\begin{figure}[t]%[!b]
    \centering
    \setlength{\abovecaptionskip}{0em}
    \setlength{\belowcaptionskip}{0em}
    \newcommand{\fighspace}{\hspace{0.2cm}}
    \subfigure[][]
    { \label{Fig:ExampleChordalc}
      \includegraphics[scale=1]{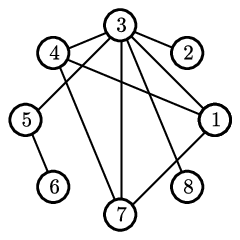}
    } \fighspace
    \subfigure[][]
    { \label{Fig:ExampleChordala}
      \includegraphics[scale=1]{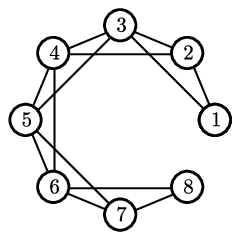}
    } \fighspace
    \subfigure[][]
    { \label{Fig:ExampleChordalb}
      \includegraphics[scale=1]{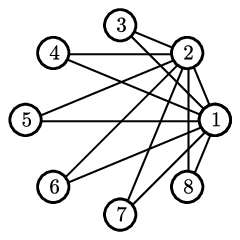}
    }
    \caption{Examples of chordal graphs: {\textit{(a)} A generic chordal graph. \textit{(b)} A ``banded'' chordal graph; \textit{(c)} A ``block-arrow'' graph. The names ``banded'' and ``block-arrow'' are motivated by the fact that, as explained in \cref{subsection:sparse_matrices}, these graphs describe the sparsity patterns of the matrices in~\Cref{Fig:ExampleSparsityPattern}.}}
    \label{Fig:ExampleChordal}
\end{figure}

A \emph{cycle} of length $ k \geq 3 $ in a graph $\mathcal{G}$ is a set of pairwise distinct vertices $ \{v_1,v_2,\ldots,v_k\}\subset \mathcal{V} $ such that $ (v_k,v_1) \in \mathcal{E} $ and $ (v_i,v_{i+1}) \in \mathcal{E} $ for $ i=1,\ldots,k-1 $. A \emph{chord} in a cycle is an edge connecting two nonconsecutive vertices.
%
%Chordal graphs are a special class of undirected graphs, defined as follows.
\begin{definition}
    An undirected graph $\mathcal{G}(\mathcal{V},\mathcal{E})$ is \emph{chordal} if every cycle of length $k\geq 4$ has at least one chord.
\end{definition}
Examples of chordal graphs are given in \cref{Fig:ExampleChordal}. Observe also that many common types of graphs are chordal, including chains, acyclic undirected graphs (i.e., graphs with no cycles, such as trees), undirected graphs with cycles of length no greater than three, and complete graphs.

Chordal graphs have a number of properties that make them easy to handle computationally. For example, a connected chordal graph has at most $n-1$ maximal cliques, and they can be identified in linear time with respect to the number of vertices and edges~\citep{vandenberghe2015chordal} using, for instance, \Cref{algorithm:maximal_clique} in \Cref{app:maximal_cliques}. In addition, any induced subgraph of a chordal graph is chordal because cycles in the subgraph are also cycles in the original graph. This is a useful fact in several induction proofs using chordality in~\citet[Section~2]{blair1993introduction}.
Finally, chordal graphs admit a so-called \textit{perfect elimination ordering} of the vertices, which is central to the zero fill-in property of Cholesky factorizations for sparse matrices. These two properties are reviewed in~\Cref{appendix:zero_fillin}. 

Given the rich structure implied by chordality, it is very often convenient to extend a nonchordal graph $\mathcal{G}(\mathcal{V},\mathcal{E})$ into a chordal graph $\hat{\mathcal{G}}(\mathcal{V}, \hat{\mathcal{E}})$ with larger edge set $\hat{\mathcal{E}} \supset \mathcal{E}$, which is called a \emph{chordal extension} of $\mathcal{G}$. Usually, a graph admits many different chordal extensions, including the trivial one with edge set $\hat{\mathcal{E}}=\mathcal{V} \times \mathcal{V}$ obtained by completion, and the one obtained by completing only the graph's connected components. Finding a \textit{minimal} chordal extension, meaning that the smallest possible number of additional edges has been added, is an NP-complete problem~\citep{yannakakis1981computing}. However, approximately minimal chordal extensions can often be constructed in practice using heuristic strategies such as the maximum cardinality search~\citep{berry2004maximum} and the symbolic Cholesky factorization with approximately minimum degree ordering~\citep{fukuda2001exploiting, vandenberghe2015chordal}.

\begin{figure}[t]
    \centering
    \subfiglabelskip=0pt
    \setlength{\abovecaptionskip}{0em}
    \footnotesize
    \subfigure[][]{%
    \begin{tikzpicture}%[every node/.style={text height=1ex, text width=1em, text centered, align=center}]
	  \matrix (m) [ampersand replacement=\&, matrix of nodes,
	  		       row sep = 1.5 em,	
	  		       column sep = 1.5 em,
  			       nodes={circle, draw=black, anchor=center}] at (0,0)
  		{ 1 \& 2  \\ 4 \& 3 \\};
		\draw (m-1-1) -- (m-1-2);
		\draw (m-2-1) -- (m-2-2);
		\draw (m-1-1) -- (m-2-1);
		\draw (m-1-2) -- (m-2-2);
	\end{tikzpicture}
	}%
	\hfill
	\subfigure[][]{%
    \begin{tikzpicture}%[every node/.style={text height=1ex, text width=1em, text centered, align=center}]
		\matrix (m2) [ampersand replacement=\&, matrix of nodes,
	  		       row sep = 1.5 em,	
	  		       column sep = 1.5 em,
  			       nodes={circle, draw=black, anchor=center}] at (0,0)
        {1 \& 2 \\ 4 \& 3 \\};
		\draw (m2-1-1) -- (m2-1-2);
		\draw (m2-2-1) -- (m2-2-2);
		\draw (m2-1-1) -- (m2-2-1);
		\draw (m2-1-2) -- (m2-2-2);
        \draw (m2-1-2) -- (m2-2-1);
	\end{tikzpicture}
	}%
	\hfill
	\subfigure[][]{%
    \begin{tikzpicture}%[every node/.style={text height=1ex, text width=1em, text centered, align=center}]
        \matrix (m3) [ampersand replacement=\&, matrix of nodes,
	  		       row sep = 1.5 em,	
	  		       column sep = 1.5 em,	
  			       nodes={circle, draw=black, anchor=center}] at (0,0)
        {1 \& 2 \\ 4 \& 3 \\};
		\draw (m3-1-1) -- (m3-1-2);
		\draw (m3-2-1) -- (m3-2-2);
		\draw (m3-1-1) -- (m3-2-1);
		\draw (m3-1-2) -- (m3-2-2);
        \draw (m3-1-1) -- (m3-2-2);
	\end{tikzpicture}
	}
	\hfill
	\subfigure[][]{%
    \begin{tikzpicture}%[every node/.style={text height=1ex, text width=1em, text centered, align=center}]
        \matrix (m4) [ampersand replacement=\&, matrix of nodes,
	  		       row sep = 1.5 em,	
	  		       column sep = 1.5 em,	
  			       nodes={circle, draw=black, anchor=center}] at (0,0)
        {1 \& 2 \\ 4 \& 3 \\};
		\draw (m4-1-1) -- (m4-1-2);
		\draw (m4-2-1) -- (m4-2-2);
		\draw (m4-1-1) -- (m4-2-1);
		\draw (m4-1-2) -- (m4-2-2);
        \draw (m4-1-2) -- (m4-2-1);
        \draw (m4-1-1) -- (m4-2-2);
	\end{tikzpicture}
	}
    \caption{\textit{(a)} A nonchordal graph: the cycle (1-2-3-4) is of length four but has no chords.
    \textit{(b)} Minimal chordal extension obtained by adding edge $(2,4)$. The maximal cliques are $\mathcal{C}_1 = \{1,2,4\}$ and $\mathcal{C}_2 = \{2,3,4\}$.
    \textit{(c)} Minimal chordal extension obtained by adding edge $(1,3)$. The maximal cliques are $\mathcal{C}_1 = \{1,2,3\}$ and $\mathcal{C}_2 = \{1,3,4\}$.
    \textit{(d)} Trivial chordal extension by completion.}
    \label{Fig:ChordalGraphs}
\end{figure}
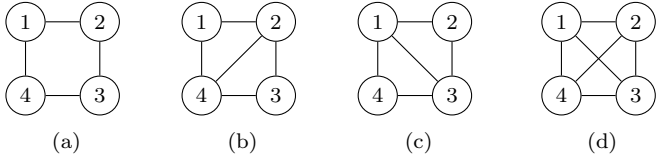

\cref{Fig:ChordalGraphs} illustrates these concepts. The graph in \cref{Fig:ChordalGraphs}(a) is not chordal, but can be extended to the chordal graph in~\cref{Fig:ChordalGraphs}(b) by adding edge $(2,4)$, edge $(1,3)$, or both. The first two extensions are minimal, while the latter is the trivial extension by completion. The minimal chordal extension obtained by adding edge $(2,4)$ has two maximal cliques, $ \mathcal{C}_1 = \{1,2,4\}$ and $\mathcal{C}_2 = \{2,3,4\}$.

\subsection{Sparse matrix decomposition} \label{subsection:SparseMatrix}
\noindent
This subsection reviews two fundamental results on the decomposition of sparse positive semidefinite matrices whose sparsity can be described using chordal graphs. 

\subsubsection{Sparse symmetric matrices} \label{subsection:sparse_matrices}
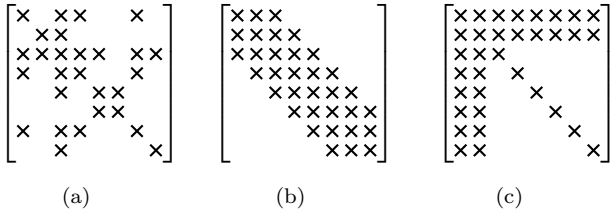
\begin{figure}[t]
    \centering
    \setlength{\abovecaptionskip}{0em}
    \setlength{\belowcaptionskip}{0em}
    \subfigure[]{
    \begingroup % keep the change local
    \setlength\arraycolsep{1pt}
    \def\arraystretch{0.6}
    \begin{tikzpicture}
    \node at (0,0) {$\left[\begin{matrix}\\[-1.2ex]
    	\mycross{black} &  & \mycross{black} &\mycross{black} & & & \mycross{black} & \\
    	  & \mycross{black} & \mycross{black} &  & & & &   \\
    	 \mycross{black} & \mycross{black} &\mycross{black} & \mycross{black} & \mycross{black} & & \mycross{black}& \mycross{black}   \\
    	\mycross{black} &  & \mycross{black} &\mycross{black} & & & \mycross{black} & \\
    	 & &\mycross{black}  && \mycross{black} & \mycross{black} &  &   \\
    	 & &&&\mycross{black} &\mycross{black} &  &   \\
    	 \mycross{black} & &\mycross{black} &\mycross{black} && &\mycross{black} &   \\
    	 & &\mycross{black} &&&&&\mycross{black}    \\%[0.5ex]
    	\end{matrix}\right]$};
    \end{tikzpicture}
    \endgroup
    }
    \hspace{-5mm}
    \subfigure[]{
    \begingroup % keep the change local
    \setlength\arraycolsep{1pt}
    \def\arraystretch{0.6}
    \begin{tikzpicture}
    \node at (0,0) {$\left[\begin{matrix}\\[-1.2ex]
    	\mycross{black} & \mycross{black} & \mycross{black} & & & &  \\
    	 \mycross{black} & \mycross{black} & \mycross{black} & \mycross{black} & & & &   \\
    	 \mycross{black} & \mycross{black} &\mycross{black} & \mycross{black} & \mycross{black} & & &    \\
    	 & \mycross{black} &\mycross{black} &\mycross{black} & \mycross{black} & \mycross{black} & &     \\
    	 & &\mycross{black} &\mycross{black} &\mycross{black} & \mycross{black} & \mycross{black} &   \\
    	 & &&\mycross{black} &\mycross{black} &\mycross{black} & \mycross{black} & \mycross{black}   \\
    	 & &&&\mycross{black} &\mycross{black} &\mycross{black} & \mycross{black}  \\
    	 & &&&&\mycross{black} &\mycross{black} &\mycross{black}    \\%[0.5ex]
    	\end{matrix}\right]$};
    \end{tikzpicture}
    \endgroup
    }
    \hspace{-5mm}
    \subfigure[]{
    \label{fig:sparsematrix_b}
    \begingroup % keep the change local
    \setlength\arraycolsep{1pt}
    \def\arraystretch{0.6}
    \begin{tikzpicture}
    \node at (0,0) {$\left[\begin{matrix}\\[-1.2ex]
    	\mycross{black} & \mycross{black} & \mycross{black} &\mycross{black} & \mycross{black}& \mycross{black}& \mycross{black} & \mycross{black} \\
    	 \mycross{black} & \mycross{black} & \mycross{black} & \mycross{black} & \mycross{black} &\mycross{black} &\mycross{black} &\mycross{black}   \\
    	 \mycross{black} & \mycross{black} &\mycross{black} & &  & & &    \\
    	 \mycross{black}& \mycross{black} &&\mycross{black} & & & &     \\
    	 \mycross{black}& \mycross{black}& &&\mycross{black} & & &   \\
    	 \mycross{black}& \mycross{black}&& && \mycross{black}&  &   \\
    	 \mycross{black}& \mycross{black}&&&&&\mycross{black} &  \\
    	 \mycross{black}&\mycross{black} &&&&&&\mycross{black}    \\%[0.5ex]
    	\end{matrix}\right]$};
    \end{tikzpicture}
    \endgroup
    }
    \caption{Sparsity patterns of $8\times 8$ matrices corresponding to the chordal graphs in~\cref{Fig:ExampleChordala,Fig:ExampleChordalb,Fig:ExampleChordalc}, respectively (throughout this paper, $\mycross{black}$ denotes a real number).}%: (a) a banded sparsity pattern; (b) a block-arrow sparsity pattern; (c) a generic sparsity pattern.}
    \label{Fig:ExampleSparsityPattern}
\end{figure}

\noindent
Fix any positive integer $n$ and set $\mathcal{V} = \{1,\ldots,n\}$. Given an undirected graph $\mathcal{G}(\mathcal{V},\mathcal{E})$, we say that a symmetric matrix $X \in \mathbb{S}^n$ has a sparsity graph $\mathcal{G}$ (alternatively, sparsity pattern $\mathcal{E}$) if $X_{ij} = X_{ji} = 0$ when %$i \neq j$ and
$(i,j) \notin \mathcal{E}$. We denote the space of sparse symmetric matrices by
\begin{equation*}
    \mathbb{S}^n(\mathcal{E},0) \!:=\!\{X \in \mathbb{S}^n \mid X_{ij} = X_{ji} = 0, \text{if}\; %i \neq j \, \text{and} \,
    (i,j) \notin \mathcal{E}\}.
\end{equation*}
For example, the graph\footnote{Throughout, we assume that each vertex has a self-loop, unless otherwise noted. We omit the self-loops when plotting a graph.} in~\cref{Fig:ChordalGraphs}(b) describes the sparsity pattern of the matrix
\begin{equation} \label{Eq:Ch2ExampleMatrix}
    X = \begin{bmatrix} X_{11} & X_{12} & 0 & X_{14} \\
                        X_{21} & X_{22} & X_{23} & X_{24} \\
                        0 & X_{32} & X_{33} & X_{34} \\
                        X_{41} & X_{42} & X_{43} & X_{44}
                        \end{bmatrix} \in \mathbb{S}^4,
\end{equation}
where each entry $X_{ij}$ may be nonzero or zero. Similarly, the symbolic matrices in~\cref{Fig:ExampleSparsityPattern} have sparsity patterns described by the graphs in~\cref{Fig:ExampleChordal}.

Given $X \in \mathbb{S}^n(\mathcal{E},0)$, the diagonal elements $X_{ii}$ and the off-diagonal elements $X_{ij}$ with $(i,j) \in \mathcal{E}$ may be nonzero or zero. Thus, if $X \in \mathbb{S}^n(\mathcal{E},0)$ and $\hat{\mathcal{E}} \supset \mathcal{E}$ is an extension of the edge set, then we also have $X \in \mathbb{S}^n(\hat{\mathcal{E}},0)$. In this paper, we are especially interested in chordal extensions of sparsity pattern. For simplicity, we will say that a matrix $X$ has a chordal sparsity pattern if its corresponding sparsity graph $\mathcal{G}(\mathcal{V},\mathcal{E})$ is chordal. Of course, this can always be achieved via chordal extension. % of the matrix's sparsity pattern.

In what follows, it will be convenient to refer to particular principal submatrices of a sparse matrix, indexed by the maximal cliques of its sparsity graph. Given a clique $\mathcal{C}_k$ of $\mathcal{G}(\mathcal{V},\mathcal{E})$, we define a matrix $E_{\mathcal{C}_k} \in \mathbb{R}^{\vert  \mathcal{C}_k\vert  \times n}$ with entries
\begin{equation} \label{Eq:IndexMatrix}
    (E_{\mathcal{C}_k})_{ij} = \begin{cases} 1, \quad \text{if } {\mathcal{C}_k}(i) = j, \\ 0, \quad \text{otherwise}, \end{cases}
\end{equation}
where $\mathcal{C}_k(i)$ is the $i$-th vertex\footnote{The elements of $\mathcal{C}_k$ can be sorted in any convenient order. We implicitly use the natural ordering in this work, but using a different one simply amounts to a permutation of the columns of $E_{\mathcal{C}_k}$}.
Given $X \in \mathbb{S}^n$, the definition of $E_{\mathcal{C}_k}$ implies that the operation
$
    E_{\mathcal{C}_k}XE_{\mathcal{C}_k}^{\tr} \in \mathbb{S}^{\vert  \mathcal{C}_k\vert }
$
extracts the principal submatrix of $X$ indexed by the clique $\mathcal{C}_k$.
Conversely, the operation $E_{\mathcal{C}_k}^\tr YE_{\mathcal{C}_k}$ ``inflates'' a $\vert\mathcal{C}_k\vert \times \vert \mathcal{C}_k\vert $ matrix $Y$ into a sparse $n \times n$ symmetric matrix that has $Y$ as its principal submatrix indexed by $\mathcal{C}_k$, and is zero otherwise. For example, the chordal graph in~\cref{Fig:ChordalGraphs}(b) has a maximal clique $\mathcal{C}_1 = \{1,2,4\}$, and the corresponding matrix $E_{\mathcal{C}_1}$ is
\begin{equation*}
    E_{\mathcal{C}_1} = \begin{bmatrix} 1 & 0 & 0& 0 \\ 0 & 1 & 0 & 0 \\ 0 & 0& 0 & 1\end{bmatrix}.
\end{equation*}
For the sparse matrix $X\in\mathbb{S}^4$ in~\cref{Eq:Ch2ExampleMatrix} and any matrix $Y\in\mathbb{S}^3$, we have
$$
\begin{aligned}
% E_{\mathcal{C}_1} = \begin{bmatrix} 1 & 0 & 0& 0 \\ 0 & 1 & 0 & 0 \\ 0 & 0& 0 & 1\end{bmatrix},
% &\,\,
E_{\mathcal{C}_1}XE_{\mathcal{C}_1}^\tr &= \begin{bmatrix} X_{11} & X_{12} & X_{14} \\ X_{21} & X_{22} & X_{24} \\
    X_{41} & X_{42} & X_{44}\end{bmatrix},
\\
E_{\mathcal{C}_1}^\tr YE_{\mathcal{C}_1} &= \begin{bmatrix} Y_{11} & Y_{12} & 0 & Y_{13} \\ Y_{21} & Y_{22} & 0 & Y_{23}\\ 0 & 0 & 0 & 0 \\
    Y_{31} & Y_{32} & 0 & Y_{33}\end{bmatrix}.
\end{aligned}
$$

\subsubsection{Cone of sparse positive semidefinite matrices} \label{subsection:ch2PSDdecomposition}
\noindent
Denote the set of positive semidefinite matrices with sparsity pattern $\mathcal{E}$ by
\begin{equation*}
\mathbb{S}^n_+(\mathcal{E},0) := \mathbb{S}^n(\mathcal{E},0) \cap \mathbb{S}^n_+.
\end{equation*}
This set is a convex cone because it is the intersection of a subspace and a convex cone.
If $\mathcal{G}(\mathcal{V},\mathcal{E})$ is a chordal graph, $\mathbb{S}^n_+(\mathcal{E},0)$ can be represented using smaller but coupled convex cones, as stated in the following result (\citealt[Theorem 2.3]{agler1988positive}; \citealt[Theorem 4]{griewank1984existence}; \citealt[Theorem 1]{kakimura2010direct}).
\begin{theorem}
\label{T:ChordalDecompositionTheorem}
     Let $\mathcal{G}(\mathcal{V},\mathcal{E})$ be a chordal graph with maximal cliques $\mathcal{C}_1,\mathcal{C}_2, \ldots, \mathcal{C}_t$.
     Then, $Z\in\mathbb{S}^n_+(\mathcal{E},0)$ if and only if there exist matrices $Z_k \in \mathbb{S}^{\vert \mathcal{C}_k \vert}_+$ for $k=1,\,\ldots,\,t$ such that
    \begin{equation} \label{eq:ChordalMatrixDecomposition}
    Z = \sum_{k=1}^{t} E_{\mathcal{C}_k}^\tr Z_k E_{\mathcal{C}_k}.
    \end{equation}
\end{theorem}

The ``if'' part of \cref{T:ChordalDecompositionTheorem} is immediate, since a sum of positive semidefinite matrices is positive semidefinite. The ``only if'' part, instead, can be proven using the zero fill-in property of sparse Cholesky factorization for $Z \in \mathbb{S}^n_+(\mathcal{E},0)$~\cite[Section 9.2]{vandenberghe2015chordal}; see~\cref{appendix:zero_fillin,appendix:proof} for details. A similar elementary proof given by \cite{kakimura2010direct}, based on simple linear algebra and perfect elimination orderings for chordal graphs, reveals that one can impose a rank constraint in the decomposition \cref{eq:ChordalMatrixDecomposition}: there exist $Z_k$ with $\text{rank}(Z) = \sum_{k=1}^t \text{rank}(Z_k)$ such that~\cref{eq:ChordalMatrixDecomposition} holds.

\begin{remark}
    The chordality assumption in \cref{T:ChordalDecompositionTheorem} is necessary. For every nonchordal pattern $\mathcal{E}$, while particular matrices in $\mathbb{S}^n_+(\mathcal{E},0)$ admit the decomposition~\cref{eq:ChordalMatrixDecomposition}, there always exist matrices in $\mathbb{S}^n_+(\mathcal{E},0)$ that do not; see \citet[p.~342]{vandenberghe2015chordal} for an explicit example.  In addition, the decomposition~\cref{eq:ChordalMatrixDecomposition} generally requires all maximal cliques $\mathcal{C}_1, \ldots, \mathcal{C}_t$, even when a subset of maximal cliques has already covered the sparsity pattern $\mathcal{E}$ (that is $\mathcal{E} = \bigcup_{k \in \mathcal{I}} \mathcal{C}_k \times \mathcal{C}_k$ with $\mathcal{I} \subset \{1, \ldots, t\}$). An example of this is given in \Cref{app:maximal_cliques}. \markendexample
\end{remark}

\begin{example}
    %Similar to~\cref{eq:ex_M_decomposition},
    Consider the positive semidefinite matrix
    \begin{equation} \label{Eq:ExamplePSD1}
        Z = \begin{bmatrix} 2 & 1 & 0 \\ 1 & 1 & 1 \\ 0 & 1 & 2 \end{bmatrix},
    \end{equation}
    whose sparsity graph is a chordal chain graph with three vertices, edge set $\mathcal{E} = \{(1,1),\,(2,2),\,(1,2),\,(2,3)\}$, and maximal cliques $\mathcal{C}_1 = \{1,2\}$ and $\mathcal{C}_2 = \{2,3\}$. \Cref{T:ChordalDecompositionTheorem} guarantees that the decomposition~\cref{eq:ChordalMatrixDecomposition} exists. Indeed, we have
    $$
    E_{\mathcal{C}_1} = \begin{bmatrix} 1&0&0 \\0& 1& 0  \end{bmatrix}, \quad
    E_{\mathcal{C}_2} = \begin{bmatrix} 0&1&0 \\0& 0& 1  \end{bmatrix},
    $$
    and
    $$
         Z =  E_{\mathcal{C}_1}^\tr\underbrace{\begin{bmatrix} 2 & 1  \\ 1 & 0.5  \end{bmatrix}}_{Z_1\succeq 0}E_{\mathcal{C}_1}  + E_{\mathcal{C}_2}^\tr\underbrace{\begin{bmatrix} 0.5 & 1 \\  1 & 2 \end{bmatrix}}_{Z_2\succeq 0}E_{\mathcal{C}_2}.
    $$
    This decomposition satisfies the rank constraint mentioned above since $\text{rank}(Z) = 2$ and $\text{rank}(Z_1)=\text{rank}(Z_2)=1$. \markendexample
\end{example}

\begin{figure}[t]
    \centering
    \setlength{\abovecaptionskip}{0em}
    \setlength{\belowcaptionskip}{0em}
    \newcommand{\fighspace}{\hspace{0.1cm}}
    \subfigure[]
    { \label{Fig:LMIdecomposition_a}
      \includegraphics[scale=0.38]{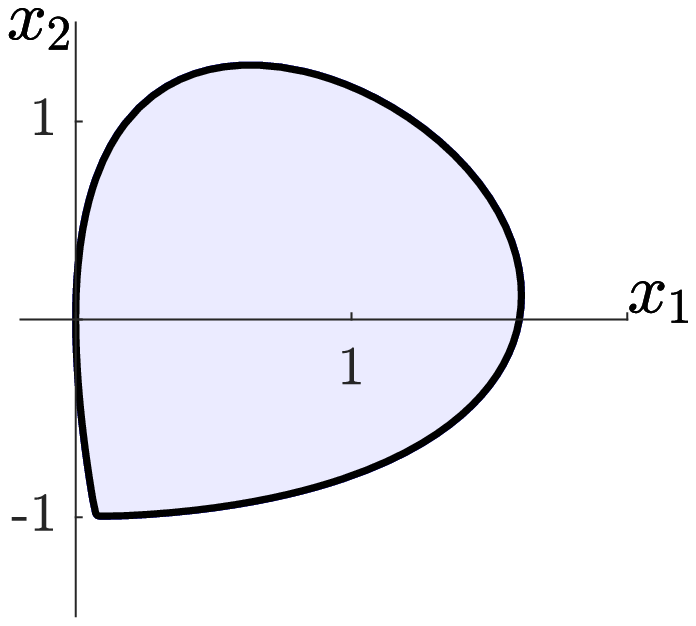}
    } %\fighspace
    \hfill
    \subfigure[]
    {
      \includegraphics[scale=0.38]{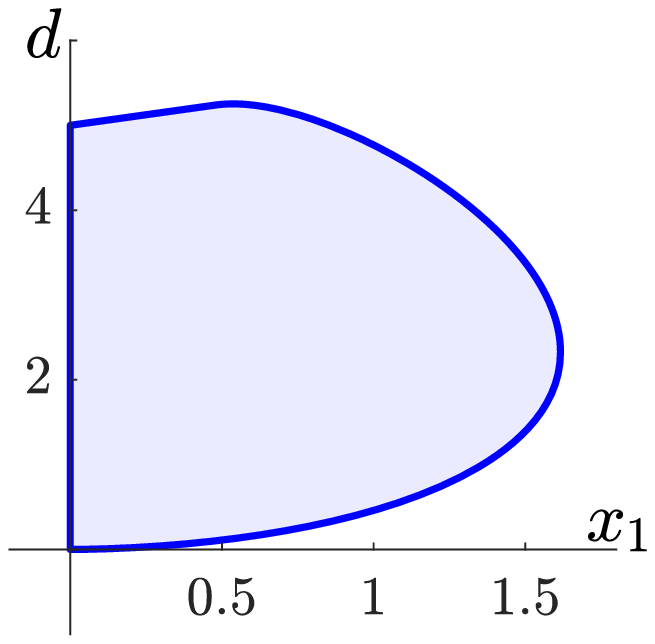}
    } %\fighspace
    \hfill
    \subfigure[]
    {
      \includegraphics[scale=0.38]{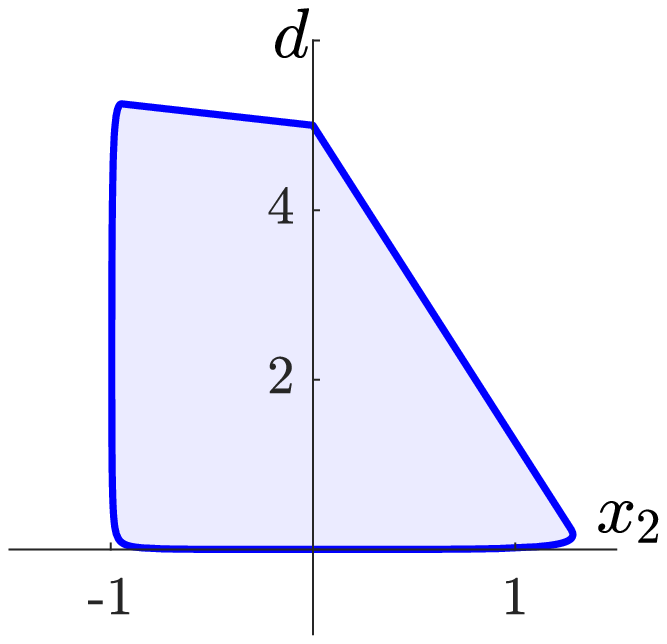}
    }
    \caption{Joint feasible set of the decomposed LMIs in~\cref{Eq:ExampleLMIDe}: \textit{(a)} projection onto the $(x_1,x_2)$ plane, \textit{(b)} projection onto the $(x_1, d)$ plane, \textit{(c)} projection onto the $(x_2,d)$ plane. Panel \textit{(a)} also shows the boundary of the feasible set of the original $3 \times 3$ LMI~\cref{Eq:ExampleLMI1}.}
    \label{Fig:LMIdecomposition}
\end{figure}

\begin{example} \label{Example:PSDdecomposition}
    Given a variable $x \in \mathbb{R}^2$, consider the $3 \times 3$ linear matrix inequality (LMI)
    \begin{equation}\label{Eq:ExampleLMI1}
        Z(x):=\begin{bmatrix} 2x_1 & x_1 + x_2 & 0 \\
                        x_1 + x_2 & 5 - x_1 - x_2 & x_1 \\
                        0 & x_1 & x_2 + 1 \end{bmatrix} \succeq 0.
    \end{equation}
    This LMI has the same chordal sparsity pattern as the matrix in~\cref{Eq:ExamplePSD1}. Consequently,~\cref{T:ChordalDecompositionTheorem} implies that~\cref{Eq:ExampleLMI1} holds if and only if there exist matrices
    $$
        Z_1:= \begin{bmatrix} a & b \\
                b & c \\  \end{bmatrix} \succeq 0
                \quad\text{and}\quad
         Z_2:= \begin{bmatrix} d & e \\
                e & f \\  \end{bmatrix} \succeq 0
    $$
    such that
    $$
        \begin{bmatrix} a & b & 0 \\
                       b & c + d &e\\
                        0 &e & f \end{bmatrix} = Z(x).
    $$
    After eliminating the variables $a$, $b$, $c$, $e$ and $f$ using this matching condition, we conclude that~\cref{Eq:ExampleLMI1} holds if and only if there exists $d$ such that
    \begin{equation} \label{Eq:ExampleLMIDe}
    \begin{aligned}
        \begin{bmatrix} 2x_1 & x_1 + x_2 \\
                x_1 + x_2 & 5 - x_1 - x_2  - d  \end{bmatrix} &\succeq 0,  \\
         \begin{bmatrix} d & x_1 \\
                x_1 & x_2 + 1 \\  \end{bmatrix} &\succeq 0.
        \end{aligned}
    \end{equation}
    \cref{Fig:LMIdecomposition} shows two-dimensional projections of the three-dimensional feasible set of the two LMIs in~\cref{Eq:ExampleLMIDe}. As expected, the projection on the $(x_1,x_2)$ plane coincides with the feasible set of LMI~\cref{Eq:ExampleLMI1}, which is contained inside the thick black line in \cref{Fig:LMIdecomposition}(a). This confirms that the LMIs in~\cref{Eq:ExampleLMIDe} are equivalent to the LMI~\cref{Eq:ExampleLMI1}. Therefore, we have decomposed a $3 \times 3$ LMI into two coupled LMIs of size $2 \times 2$.
    \markendexample
\end{example}

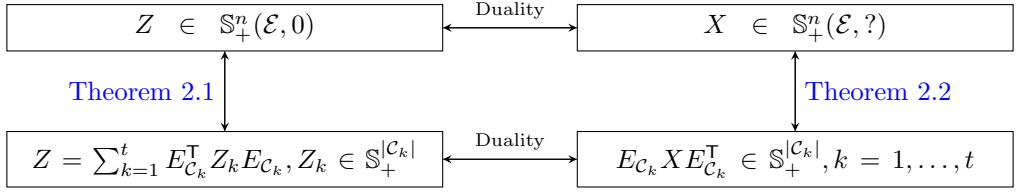
\begin{figure*}[t]
    \centering
    \setlength{\abovecaptionskip}{6pt}
    \setlength{\belowcaptionskip}{0em}
	%\footnotesize
	\begin{tikzpicture}
	  \matrix (m) [matrix of nodes,
	  		       row sep = 3.em,	
	  		       column sep = 5em,	
  			       nodes={rectangle, draw=black, align=center, text width=5.5cm}]
  	{
   	$Z \in \mathbb{S}^n_+(\mathcal{E},0)$ & $X \in \mathbb{S}^n_+(\mathcal{E},?)$ \\
   	$Z = \sum_{k=1}^t E_{\mathcal{C}_k}^\tr Z_k E_{\mathcal{C}_k}, Z_k \in \mathbb{S}^{|\mathcal{C}_k|}_+$  &
        $E_{\mathcal{C}_k} X E_{\mathcal{C}_k}^\tr \in \mathbb{S}^{|\mathcal{C}_k|}_+, k = 1, \ldots, t$ \\};
	\path[-stealth]
		(m-1-1) edge node [left, align=center]
            {\cref{T:ChordalDecompositionTheorem}}  (m-2-1)
        (m-2-1) edge (m-1-1)
		(m-1-2) edge node [right, align=center]
			{\cref{T:ChordalCompletionTheorem}}  (m-2-2)
        (m-2-2) edge (m-1-2)
		(m-1-1) edge node [above] {\scriptsize Duality} (m-1-2)
        (m-1-2) edge (m-1-1)
		(m-2-1) edge node [above] {\scriptsize Duality} (m-2-2)
        (m-2-2) edge (m-2-1);
	\end{tikzpicture}
    \caption{Summary of duality between $\mathbb{S}^n_+(\mathcal{E},0)$ and $\mathbb{S}^n_+(\mathcal{E},?)$ and duality between~\cref{T:ChordalDecompositionTheorem} and~\cref{T:ChordalCompletionTheorem} for a chordal graph $\mathcal{G}(\mathcal{V},\mathcal{E})$ with maximal cliques $\mathcal{C}_1, \ldots, \mathcal{C}_t$.}
    \label{Fig:SparseCone_Duality}
\end{figure*}

%%%%%%%%%%%%%%%%%%%%%%%%%%%%%%%%%%%%%%%%%%%%%%%%%%%%%%%%%%%%
\subsubsection{Cone of positive-semidefinite-completable matrices}
\noindent
A concept related to the matrix decomposition above is that of positive semidefinite matrix completion. Given a matrix $X \in \mathbb{S}^n$, let
\begin{equation}\label{e:matrix-sparse-projection}
    \mathbb{P}_{\mathbb{S}^n(\mathcal{E},0) }(X) = \begin{cases}
    X_{ij} &\text{if } (i,j) \in \mathcal{E},\\
    0 &\text{otherwise}
    \end{cases}
\end{equation}
be its projection onto the space of sparse matrices $\mathbb{S}^n(\mathcal{E},0)$ with respect to the Frobenius matrix norm. We define the cone
\begin{equation*}
\mathbb{S}^n_+(\mathcal{E},?) := \mathbb{P}_{\mathbb{S}^n(\mathcal{E},0) }( \mathbb{S}^n_+ ).
\end{equation*}
Using~\cref{e:matrix-sparse-projection}, it is not hard to see that a sparse matrix $X$ is in $\mathbb{S}^n_+(\mathcal{E},?) $ if and only if it has a positive semidefinite completion, meaning that some (or all) of the zero entries $X_{ij}$ with $(i,j) \notin \mathcal{E}$ can be replaced with nonzeros to obtain a positive semidefinite matrix $\overline{X}$. We call $\overline{X}$ the \textit{completion} of $X$ and refer to $\mathbb{S}^n_+(\mathcal{E},?)$ as the cone of positive-semidefinite-completable matrices.

\begin{remark}[Nonuniqueness of the positive semidefinite completion]\label{remark:completion-nonuniqueness}
    The positive semidefinite completion of a matrix $X \in \mathbb{S}_+^n(\mathcal{E},?)$ with sparsity pattern $\mathcal{E}$ is generally not unique. For a chordal sparsity pattern $\mathcal{E}$, two widely used and efficient strategies to compute a completion $\overline{X}$ are the maximum determinant completion~\citep[Chapter 10.2]{vandenberghe2015chordal}, which maximizes $\det \overline{X}$, and the minimum rank completion (see \citealt{dancis1992positive}; \citealt{jiang2017minimum}; \citealt[Chapter 3.3]{sun2015decomposition_thesis}), which minimizes $\mathrm{rank}(\overline{X})$. In particular, there exists a positive semidefinite completion $\overline{X}$ whose rank agrees with the maximum rank of the principal submatrices $E_{\mathcal{C}_i}XE_{\mathcal{C}_i}^\tr$~\cite[Theorem 1.5]{dancis1992positive}, i.e.,
    \begin{equation}\label{e:completion-rank-condition}
        \mathrm{rank}(\overline{X}) = \max_{k = 1, 2, \ldots, t} \mathrm{rank}(E_{\mathcal{C}_k}XE_{\mathcal{C}_k}^\tr).
    \end{equation}
    \markendexample
\end{remark}

For any undirected graph $\mathcal{G}(\mathcal{V},\mathcal{E})$, the cones $\mathbb{S}^n_{+}(\mathcal{E},?)$ and $\mathbb{S}_{+}^n(\mathcal{E},0)$ are dual to each other with respect to the trace inner product $\langle X,Z \rangle =\text{Trace}(XZ)$ in the space $\mathbb{S}^n(\mathcal{E},0)$~\cite[Chapter 10]{vandenberghe2015chordal}. To see this, observe that
$$
    \begin{aligned}
        (\mathbb{S}^n_{+}(\mathcal{E},?))^* &= \{Z \in \mathbb{S}^n(\mathcal{E},0) \mid \langle X,Z \rangle \geq 0, \forall X \in \mathbb{S}^n_{+}(\mathcal{E},?)\} \\
        &=\{Z \in \mathbb{S}^n(\mathcal{E},0) \mid \left\langle \mathbb{P}_{\mathbb{S}^n(\mathcal{E},0) }(X),Z \right\rangle \geq 0, \forall X \succeq 0\} \\
        & = \{Z \in \mathbb{S}^n(\mathcal{E},0) \mid \langle X,Z \rangle \geq 0, \forall X \succeq 0\} \\
        & = \{Z \in \mathbb{S}^n(\mathcal{E},0) \mid Z \succeq 0\} \\
        & = \mathbb{S}^n_+(\mathcal{E},0).
    \end{aligned}
$$
For a chordal matrix sparsity pattern, \cref{T:ChordalDecompositionTheorem} on the decomposition of the cone $\mathbb{S}^n_+(\mathcal{E},0)$ can be dualized to obtain the following characterization of $\mathbb{S}^n_{+}(\mathcal{E},?)$, first proved by~\citet[Theorem~7]{grone1984positive}.
\begin{theorem} \label{T:ChordalCompletionTheorem}
     Let $\mathcal{G}(\mathcal{V},\mathcal{E})$ be a chordal graph with maximal cliques $\mathcal{C}_1, \ldots, \mathcal{C}_t$.
     %and let $\{\mathcal{C}_1,\mathcal{C}_2, \ldots, \mathcal{C}_t\}$ be the set of its maximal cliques.
     Then, $X\in\mathbb{S}^n_+(\mathcal{E},?)$ if and only if
     %$E_{\mathcal{C}_k} X E_{\mathcal{C}_k}^\tr \in \mathbb{S}^{\vert \mathcal{C}_k \vert}_+$ for all $k=1,\,\ldots,\,t$.
     \begin{equation}
      E_{\mathcal{C}_k} X E_{\mathcal{C}_k}^\tr \in \mathbb{S}^{\vert \mathcal{C}_k \vert}_+ \qquad \forall k=1,\,\ldots,\,t.
     \end{equation}
\end{theorem}

The ``only if'' part of~\cref{T:ChordalCompletionTheorem} is immediate, since any principal submatrix of a positive semidefinite matrix is positive semidefinite. The ``if'' part, instead, relies on the properties of chordal graphs and, as mentioned above, can be proven by combining the duality between $\mathbb{S}^n_+(\mathcal{E},0)$ and $\mathbb{S}^n_+(\mathcal{E},?)$ with~\cref{T:ChordalDecompositionTheorem}~\cite[p.~357]{vandenberghe2015chordal}. Precisely,
\begin{equation*}
    \begin{aligned}
        X \in \mathbb{S}^n_+(\mathcal{E},?) &\Leftrightarrow \langle X,Z \rangle \geq 0\quad\forall Z \in \mathbb{S}^n_+(\mathcal{E},0), \\
        &\Leftrightarrow \bigg\langle X, \sum_{k=1}^{t} E_{\mathcal{C}_k}^\tr Z_k E_{\mathcal{C}_k} \bigg\rangle \geq 0\quad\forall Z_k \in \mathbb{S}^{|\mathcal{C}_k|}_+, \\
        &\Leftrightarrow \sum_{k=1}^{t} \left\langle E_{\mathcal{C}_k} X E_{\mathcal{C}_k}^\tr, Z_k  \right\rangle \geq 0\quad\forall Z_k \in \mathbb{S}^{|\mathcal{C}_k|}_+, \\
        &\Leftrightarrow  E_{\mathcal{C}_k} X E_{\mathcal{C}_k}^\tr \in \mathbb{S}^{|\mathcal{C}_k|}_+\quad\forall k = 1,\ldots, t.\\
    \end{aligned}
\end{equation*}
The first equivalence expresses the duality between $\mathbb{S}^n_+(\mathcal{E},0)$  and $\mathbb{S}^n_+(\mathcal{E},?)$, the second one follows from~\cref{T:ChordalDecompositionTheorem}, and the third one follows from the cyclic property of the trace operator: $\text{Trace}(MN) = \text{Trace}(NM)$ for any matrices $M,N$ of compatible dimensions.

\cref{Fig:SparseCone_Duality} illustrates how the duality between $\mathbb{S}^n_+(\mathcal{E},0)$ and $\mathbb{S}^n_+(\mathcal{E},?)$ is mirrored in the duality between~\cref{T:ChordalDecompositionTheorem} and~\cref{T:ChordalCompletionTheorem} for chordal graphs.

\begin{example}
    Consider the symmetric matrix
    \begin{equation*}% \label{E:ExPSDcompletion}
   X = \begin{bmatrix}
        2 & 1 & 0 \\
        1 & 0.5 & 1\\
        0 & 1 & 2
    \end{bmatrix},
\end{equation*}
whose sparsity pattern is the (by now usual) 3-node chordal chain graph with maximal cliques $\mathcal{C}_1 = \{1,2\}$ and $\mathcal{C}_2 = \{2,3\}$. It is easy to check that, while $X$ is not positive semidefinite, the principal submatrices indexed by the cliques $\mathcal{C}_1$ and $\mathcal{C}_2$ are.
Then,~\cref{T:ChordalCompletionTheorem} guarantees that $X \in \mathbb{S}^n_+(\mathcal{E},?)$, meaning that the zero entries may be replaced by nonzeros to obtain a positive semidefinite matrix $\overline{X}$. One possible positive semidefinite completion is
\begin{equation*}
    \overline{X} = \begin{bmatrix}
        2 & 1 & 2 \\
        1 & 0.5 & 1\\
        2 & 1 & 2
    \end{bmatrix}.
\end{equation*}
In fact, this is the minimum-rank completion whose rank, $\text{rank}(X) = 1$, coincides with the maximum rank of individual principal submatrices of $X$ (cf. \Cref{remark:completion-nonuniqueness}).  \markendexample
\end{example}

\begin{figure}[t]
    \centering
    \includegraphics[width=0.9\linewidth]{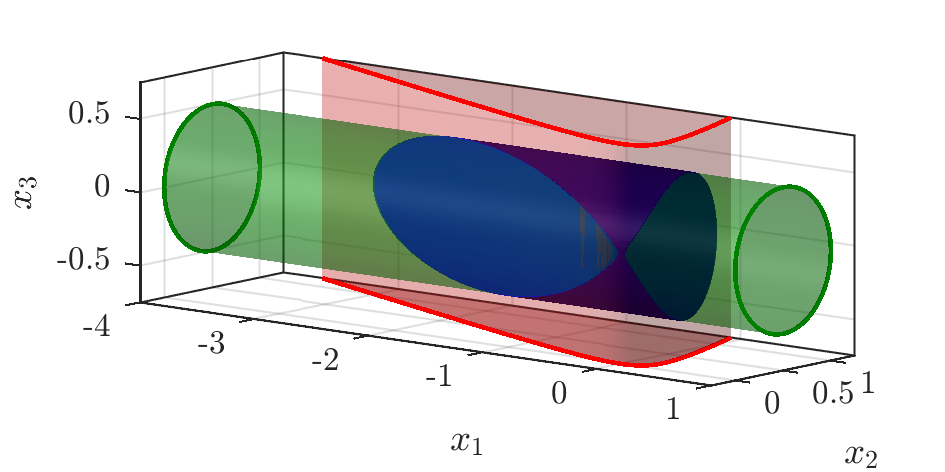}
    \caption{Region of $\mathbb{R}^3$ where the matrix $X(x)$ in~\cref{Eq:ExPSDcompletion} admits a positive semidefinite completion (blue shading). This region coincides with the intersection of the region of $\mathbb{R}^3$ where the first LMI in~\cref{Eq:ExPSDcompletionEq} is feasible (red shading; the region extends to infinity in the $x_3$ direction) and the cylindrical region of $\mathbb{R}^3$ where the second LMI in~\cref{Eq:ExPSDcompletionEq} is feasible (green shading; the region extends to infinity in the $x_1$ direction). Thick red and green lines highlight the cross section of these two regions.}
    \label{Fig:PSDcompletion}
\end{figure}

\begin{example}
    Consider the problem of finding a variable $x \in \mathbb{R}^3$ such that the matrix
    \begin{equation}\label{Eq:ExPSDcompletion}
       X(x) \!:= \!\begin{bmatrix} 1 - x_1 & x_1 + x_2 & 0 \\
        x_1 + x_2& x_2 & x_2 + x_3 \\
        0 & x_2 + x_3 & 2x_3 + 1
        \end{bmatrix} %\!\in \!\mathbb{S}^n_+(\mathcal{E},?).
    \end{equation}
    admits a positive semidefinite completion.
    This is equivalent to finding $x \in \mathbb{R}^3$ as well as a corresponding scalar $y \in \mathbb{R}$ such that
    \begin{equation}\label{Eq:ExPSDcompletion-LMI}
       \begin{bmatrix} 1 - x_1 & x_1 + x_2 & y \\
        x_1 + x_2& x_2 & x_2 + x_3 \\
        y & x_2 + x_3 & 2x_3 + 1
        \end{bmatrix} \succeq 0.
    \end{equation}
    %which involves an extra variable $y \in \mathbb{R}$.
    Since the sparsity graph of $X(x)$ is chordal, \cref{T:ChordalCompletionTheorem} implies that~\cref{Eq:ExPSDcompletion} is equivalent to the two LMIs
    \begin{equation} \label{Eq:ExPSDcompletionEq}
        \begin{bmatrix} 1-x_1 & x_1+x_2 \\x_1 + x_2 & x_2
        \end{bmatrix} \!\succeq \!0,
        \;
        \begin{bmatrix}
        x_2 & x_2 + x_3 \\
        x_2 + x_3 & 2x_3 + 1
        \end{bmatrix} \!\succeq\! 0.
    \end{equation}
    Feasible vectors $x$ for the first of these two LMIs
    %~\cref{Eq:ExPSDcompletionEq}
    can be found by imposing $1 - x_1 + x_2 \geq 0$ and $(1-x_1)x_2 - (x_1 + x_2)^2 \geq 0$, while feasible $x$ for the second LMI are found by requiring $x_2 + 2x_3 +1 \geq 0$ and $x_2(2x_3+1) - (x_2 + x_3)^2\geq 0$. The feasible sets obtained in each case are illustrated by the red and green regions in \Cref{Fig:PSDcompletion}, respectively. The blue region in the figure, instead, represents the three-dimensional set of feasible $x$ for~\cref{Eq:ExPSDcompletion-LMI}. As expected from \cref{T:ChordalCompletionTheorem}, this is exactly the intersection of the feasible regions for the two LMIs in~\cref{Eq:ExPSDcompletionEq}. Similar to~\cref{Example:PSDdecomposition}, one can therefore replace the original $3 \times 3$ completion constraint---which is equivalent to LMI~\cref{Eq:ExPSDcompletion-LMI}---with the two $2 \times 2$ LMIs in~\cref{Eq:ExPSDcompletionEq} without any loss of generality. \markendexample
\end{example}

\subsection{Block-partitioned matrices} \label{subsection:block-partition-matrices}
\noindent
\Cref{T:ChordalDecompositionTheorem,T:ChordalCompletionTheorem} can be extended to block-partitioned matrices characterized by block-sparsity. Such matrices arise, for example, when modeling network systems (cf. \Cref{subsection:control_applications}), where each block in the partition corresponds to an individual subsystem and sparsity in the network connectivity translates into block-sparsity. Block-partitioned matrices are also useful in extending factor-width decomposition that will be discussed in \Cref{section:factor-width-two}. 

\subsubsection{Sparse block matrices}\label{ss:block-matrices-notation}
\noindent
Given a positive integer $n$, any finite set of positive integers $\alpha = \{\alpha_1, \alpha_2, \ldots, \alpha_p\}$ is called a \textit{partition of $n$} if $\sum_{i=1}^p\alpha_i=n$. The set of all possible partitions of $n$ can be equipped with the following (partial) order relation.
\begin{definition}\label{definition:partition_order}
    Let $\alpha = \{\alpha_1, \ldots, \alpha_p\}$ and $\beta = \{\beta_1, \ldots, \beta_q\}$ be two partitions of an integer $n$ with $ p < q $. We say that $\beta$ is \emph{finer} than $\alpha$ (and $\alpha$ is \emph{coarser} than $\beta$), denoted by $\beta \sqsubset \alpha$, if there exist integers $\{m_1, m_2, \ldots, m_{p+1}\}$ with $m_1 = 1$, $m_{p+1} = q+1$ and $m_i < m_{i+1}$ for $i = 1, \ldots, p$ such that
    $\alpha_i = \sum_{j = m_i}^{m_{i+1} - 1} \beta_j$ for all $i = 1, \ldots, p.$
\end{definition}
Essentially, a finer partition $\beta$ breaks some entries of $\alpha$ into smaller ones (conversely,  a coarser partition $\alpha$ is obtained by merging some entries of $\beta$ into a bigger one). For example, the partitions $\alpha = \{4, 2\}, \beta = \{2, 2, 2\}$ and $\gamma = \{1, 1, 1, 1, 1, 1\}$ of $n=6$ satisfy
$\gamma \sqsubset \beta \sqsubset \alpha$.

Given any integer $n$ and any partition $\alpha =  \{\alpha_1, \ldots, \alpha_p\}$ of $n$, a matrix $M \in \mathbb{R}^{n \times n}$ can be written in the block form% is \emph{$\alpha$-partitioned} if
$$
    M = \begin{bmatrix} M_{11} & M_{12} & \ldots & M_{1p} \\
                        M_{12} & M_{22} & \ldots & M_{2p} \\
                        \vdots & \vdots & \ddots & \vdots \\
                        M_{p1} & M_{p2} & \ldots & M_{pp} \\  \end{bmatrix}
$$
with $M_{ij} \in \mathbb{R}^{\alpha_i \times \alpha_j}$ for all $i, j = 1, \ldots, p$. %pairs $(i,j) \in \{1, \ldots, n\}\times \{1, \ldots, n\}$. 
%An $N \times N$ matrix $M$ can obviously be partitioned in many different ways, and the trivial choice
For the finest partition $\alpha = \{1,\ldots,1\} = \mathbf{1}_n$, the block $M_{ij}$ reduces to the entry $(i,j)$ of $M$. As shown below and in \cref{section:blockFW2}, however, the freedom to consider a nontrivial partition offers considerable flexibility when devising decomposition strategies for a large matrix $M$. In particular, by refining or coarsening a partition one can in principle split a matrix into blocks of optimal size for the computational resources at one's disposal.

The block sparsity pattern of an $n \times n$ matrix $M$ whose blocks are defined by a partition $\alpha = \{\alpha_1,\ldots,\alpha_p\}$ of $n$ can be described using a graph $\mathcal{G}(\mathcal{V},\mathcal{E})$ with $\mathcal{V} = \{1, 2, \ldots, p\}$ and edge set such that $M_{ij}=0$ if $(i,j) \notin \mathcal{E}$, where $M_{ij}$ is the $(i,j)$-th block in $M$ and $0$ denotes a zero block of appropriate size. We call $\alpha$ a \textit{chordal partition} if the corresponding block sparsity graph $\mathcal{G}(\mathcal{V},\mathcal{E})$ is chordal. The linear space of sparse symmetric block matrices with a prescribed block sparsity pattern $\mathcal{E}$ is then given by
\begin{equation*} 
    \mathbb{S}^n_{\alpha}(\mathcal{E},0) := \{M \in \mathbb{S}^{n} | M_{ij} \!= \!0\; \text{if} \; (i,j) \!\notin \!\mathcal{E} \}.
\end{equation*}
The block-sparse positive semidefinite cone and the block-sparse positive-semidefinite-completable cone are simply
\begin{subequations}
\begin{align}
    \mathbb{S}^n_{\alpha,+}(\mathcal{E},0) &:= \mathbb{S}^n_{\alpha}(\mathcal{E},0) \cap \mathbb{S}^n_+,\\
    \mathbb{S}^n_{\alpha,+}(\mathcal{E},?) &:= \mathbb{P}_{\mathbb{S}^n_{\alpha}(\mathcal{E},0)}(\mathbb{S}^n_+).
\end{align}
\end{subequations}

\begin{figure}
    \centering
    \footnotesize
    \subfigure[][]{%
        \begin{tikzpicture}
            \matrix (m) [ampersand replacement=\&, matrix of nodes,
	  		             row sep = 1em,	
	  		             column sep = 1em,	
  			             nodes={circle, draw=black}] at (0,0)
      		{7 \& |[fill=blue!10]|8 \& |[fill=blue!10]|9 \\ |[fill=blue!10]|4 \& 5 \& |[fill=blue!10]|6\\ |[fill=blue!10]|1 \& |[fill=blue!10]|2 \& 3\\};
      		% Horizontal
    		\draw[ultra thick, blue] (m-3-1) -- (m-3-2);
    		\draw (m-3-2) -- (m-3-3);
    		\draw (m-2-1) -- (m-2-2);
    		\draw (m-2-2) -- (m-2-3);
            \draw (m-1-1) -- (m-1-2);
    		\draw[ultra thick, blue] (m-1-2) -- (m-1-3);
    		% Vertical
    		\draw (m-1-1) -- (m-2-1);
    		\draw[ultra thick, blue] (m-2-1) -- (m-3-1);
            \draw (m-1-2) -- (m-2-2);
    		\draw (m-2-2) -- (m-3-2);
    		\draw[ultra thick, blue] (m-1-3) -- (m-2-3);
    		\draw (m-2-3) -- (m-3-3);
    		% Diagonals
    		\draw (m-3-1) -- (m-2-2);
    		\draw[ultra thick, blue] (m-3-2) -- (m-2-3);
    		\draw[ultra thick, blue] (m-2-1) -- (m-1-2);
    		\draw (m-2-2) -- (m-1-3);
        \end{tikzpicture}
    }%
    \hspace{1.5em}
    \subfigure[][]{%
    	\tikzset{
    		myrectnodes/.style={rounded corners, draw=black},
    		mycircnodes/.style={circle, draw=black}
    	}
        \begin{tikzpicture}
            \matrix (m) [ampersand replacement=\&, matrix of nodes,
	  		             row sep = 1em,	
	  		             column sep = 1em] at (0,0)
      		{|[myrectnodes]|7, 8\& |[mycircnodes]|9 \\ |[mycircnodes]|4 \&|[myrectnodes]| 5, 6\\|[myrectnodes]| 1, 2\& |[mycircnodes]| 3\\};
      		% Links for node (1,2)
      		\draw (m-3-1) -- (m-3-2);
      		\draw (m-3-1) -- (m-2-1);
      		\draw (m-3-1) -- (m-2-2);
      		% Links for node (3)
      		\draw (m-3-2) -- (m-2-2);
      		% Links for node (4)
      		\draw (m-2-1) -- (m-2-2);
      		\draw (m-2-1) -- (m-1-1);
      		% Links for node (5,6)
      		\draw (m-2-2) -- (m-1-2);
      		\draw (m-2-2) -- (m-1-1);
      		% Links for node (7,8)
      		\draw (m-1-1) -- (m-1-2);
        \end{tikzpicture}
    }%
    \hspace{1.5em}
    \subfigure[][]{%
        \begin{tikzpicture}
            \matrix (m) [ampersand replacement=\&, matrix of nodes,
	  		             row sep = 1.25em,	
	  		             column sep = 1em,	
  			             nodes={rounded corners, draw=black}] at (0,0)
      		{7, 8, 9\\ 4, 5, 6\\ 1, 2, 3\\};
      		\draw (m-1-1) -- (m-2-1);
      		\draw (m-2-1) -- (m-3-1);
        \end{tikzpicture}
    }%
    \caption{
    \textit{(a)} Nonchordal sparsity graph of the $9\times 9$ matrix $M$ in \cref{remark:chordal-partition}. Blue nodes and edges form a cycle of length 6 with no chord.
    \textit{(b)} Chordal block sparsity graph of the same matrix with partition $\alpha_1=\{2,1,1,2,2,1\}$. 
    \textit{(c)} Chordal block sparsity graph of the same matrix with partition $\alpha_2 = \{3,3,3\}$.
	}
    \label{fig:chordal-partitions}
\end{figure}
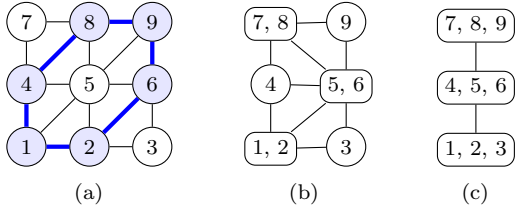

\begin{remark}[Chordal partitions and chordal extension]\label{remark:chordal-partition}
If $M$ is a sparse matrix with nonchordal sparsity pattern, it is often possible to find one or more chordal partitions $\alpha$. An example is the $9 \times 9$ symbolic matrix
\begin{equation*}
    \setlength\arraycolsep{0.5pt}
    \def\arraystretch{0.3}
    M = \left[
    \small
	\begin{array}{ccc ccc ccc}\\[0.1ex]
	\mycross{black} &\mycross{black} & & \mycross{black} &\mycross{black}&&&\\
	\mycross{black} &\mycross{black} &\mycross{black} & &\mycross{black} &\mycross{black} &\\
	 &\mycross{black} &\mycross{black} & & &\mycross{black} &\\
	\mycross{black} & & & \mycross{black}  &\mycross{black} & &\mycross{black}  &\mycross{black} \\
	\mycross{black} &\mycross{black} & &\mycross{black} &\mycross{black} &\mycross{black} &&\mycross{black} &\mycross{black} \\
	&\mycross{black} &\mycross{black} & &\mycross{black}  &\mycross{black} & & & \mycross{black} \\
	&&&\mycross{black} &&&\mycross{black} &\mycross{black} \\
	&&&\mycross{black} &\mycross{black} &&\mycross{black} &\mycross{black} &\mycross{black} \\
	&&&&\mycross{black} &\mycross{black} &&\mycross{black} &\mycross{black} \\[0.15ex]
	\end{array}
	\right]
	=
	\left[
	\small
	\begin{array}{cc | c | c | cc | cc |c}\\[0.1ex]
	\mycross{black} &\mycross{black} & & \mycross{black} &\mycross{black}&&&\\
	\mycross{black} &\mycross{black} &\mycross{black} & &\mycross{black} &\mycross{black} & &\\\hline
	 &\mycross{black} &\mycross{black} & & &\mycross{black} & & & \\\hline
	\mycross{black} & & & \mycross{black}  &\mycross{black} & &\mycross{black}  &\mycross{black} \\\hline
	\mycross{black} &\mycross{black} & &\mycross{black} &\mycross{black} &\mycross{black} &&\mycross{black} &\mycross{black} \\
	&\mycross{black} &\mycross{black} & &\mycross{black}  &\mycross{black} & & & \mycross{black} \\\hline
	&&&\mycross{black} &&&\mycross{black} &\mycross{black} \\
	&&&\mycross{black} &\mycross{black} &&\mycross{black} &\mycross{black} &\mycross{black} \\\hline
	&&&&\mycross{black} &\mycross{black} &&\mycross{black} &\mycross{black} \\[0.15ex]
	\end{array}
	\right]
	=
	\left[
	\small
	\begin{array}{ccc | ccc | ccc}\\[0.1ex]
	\mycross{black} &\mycross{black} & & \mycross{black} &\mycross{black}&&&\\
	\mycross{black} &\mycross{black} &\mycross{black} & &\mycross{black} &\mycross{black} &\\
	 &\mycross{black} &\mycross{black} & & &\mycross{black} &\\\hline
	\mycross{black} & & & \mycross{black}  &\mycross{black} & &\mycross{black}  &\mycross{black} \\
	\mycross{black} &\mycross{black} & &\mycross{black} &\mycross{black} &\mycross{black} &&\mycross{black} &\mycross{black} \\
	&\mycross{black} &\mycross{black} & &\mycross{black}  &\mycross{black} & & & \mycross{black} \\\hline
	&&&\mycross{black} &&&\mycross{black} &\mycross{black} \\
	&&&\mycross{black} &\mycross{black} &&\mycross{black} &\mycross{black} &\mycross{black} \\
	&&&&\mycross{black} &\mycross{black} &&\mycross{black} &\mycross{black} \\[0.15ex]
	\end{array}
	\right]
\end{equation*}
where the partitions $\alpha_1=\{2,1,1,2,2,1\}$ and $\alpha_2 = \{3,3,3\}$ are both chordal (the corresponding block sparsity graphs are illustrated in \cref{fig:chordal-partitions}). For a given chordal partition, in this example but also in general, completing all blocks of $M$ that are not identically zero results in a chordal extension of $M$. For instance, the chordal extension of the $9 \times 9$ matrix above resulting from the partitions $\alpha_1$ and $\alpha_2$ are, respectively,
\begin{equation*}
    \setlength\arraycolsep{0.5pt}
    \def\arraystretch{0.3}
   \left[
	\small
	\begin{array}{cc | c | c | cc | cc |c}\\[0.1ex]
	\mycross{black} &\mycross{black} &\mycross{matlabred} & \mycross{black} &\mycross{black}& \mycross{matlabred} &&\\
	\mycross{black} &\mycross{black} &\mycross{black} & \mycross{matlabred} &\mycross{black} &\mycross{black} & &\\\hline
	\mycross{matlabred} &\mycross{black} &\mycross{black} & & \mycross{matlabred} &\mycross{black} & & & \\\hline
	\mycross{black} &\mycross{matlabred} &  & \mycross{black}  &\mycross{black} & \mycross{matlabred} &\mycross{black}  &\mycross{black} & \\\hline
	\mycross{black} &\mycross{black} & \mycross{matlabred} &\mycross{black} &\mycross{black} &\mycross{black} & \mycross{matlabred}&\mycross{black} &\mycross{black} \\
	\mycross{matlabred} &\mycross{black} &\mycross{black} & \mycross{matlabred} &\mycross{black}  &\mycross{black} & \mycross{matlabred} & \mycross{matlabred} & \mycross{black} \\\hline
	&&&\mycross{black} &\mycross{matlabred}&\mycross{matlabred}&\mycross{black} &\mycross{black} &\mycross{matlabred} \\
	&&&\mycross{black} &\mycross{black} &\mycross{matlabred} &\mycross{black} &\mycross{black} &\mycross{black} \\\hline
	&&&&\mycross{black} &\mycross{black} & \mycross{matlabred}&\mycross{black} &\mycross{black} \\[0.5ex]
	\end{array}
	\normalsize
	\right] 
	\quad\text{and}\quad
	\left[
	\small
	\begin{array}{ccc | ccc | ccc}\\[0.1ex]
	\mycross{black} &\mycross{black} &\mycross{matlabred} & \mycross{black} &\mycross{black}& \mycross{matlabred} &&\\
	\mycross{black} &\mycross{black} &\mycross{black} & \mycross{matlabred} &\mycross{black} &\mycross{black} & &\\
	\mycross{matlabred} &\mycross{black} &\mycross{black} & \mycross{matlabred}& \mycross{matlabred} &\mycross{black} & & & \\\hline
	\mycross{black} &\mycross{matlabred} & \mycross{matlabred} & \mycross{black}  &\mycross{black} & \mycross{matlabred} &\mycross{black}  &\mycross{black} & \mycross{matlabred} \\
	\mycross{black} &\mycross{black} & \mycross{matlabred} &\mycross{black} &\mycross{black} &\mycross{black} & \mycross{matlabred}&\mycross{black} &\mycross{black} \\
	\mycross{matlabred} &\mycross{black} &\mycross{black} & \mycross{matlabred} &\mycross{black}  &\mycross{black} & \mycross{matlabred} & \mycross{matlabred} & \mycross{black} \\\hline
	&&&\mycross{black} &\mycross{matlabred}&\mycross{matlabred}&\mycross{black} &\mycross{black} &\mycross{matlabred} \\
	&&&\mycross{black} &\mycross{black} &\mycross{matlabred} &\mycross{black} &\mycross{black} &\mycross{black} \\
	&&&\mycross{matlabred}&\mycross{black} &\mycross{black} & \mycross{matlabred}&\mycross{black} &\mycross{black} \\[0.5ex]
	\end{array}
	\normalsize
	\right],
\end{equation*}
where entries colored in red have been added by the block-completion process.
Finding a chordal partition for a matrix, therefore, gives a way of performing a particular chordal extension of its sparsity pattern. The opposite, however, is not true: not all chordal extensions are obtained via a block-completion operation. One example for the $9\times 9$ matrix above is the chordal extension
\begin{equation*}
    \setlength\arraycolsep{0.5pt}
    \def\arraystretch{0.3}
   \left[
	\small
	\begin{array}{ccccccccc}\\[0.1ex]
	\mycross{black} &\mycross{black} & & \mycross{black} &\mycross{black}& \mycross{matlabred} &&&\mycross{matlabred}\\
	\mycross{black} &\mycross{black} &\mycross{black} & &\mycross{black} &\mycross{black} & &\\
	&\mycross{black} &\mycross{black} & & &\mycross{black} & & & \\
	\mycross{black} & &  & \mycross{black}  &\mycross{black} & &\mycross{black}  &\mycross{black} & \mycross{matlabred}\\
	\mycross{black} &\mycross{black} &  &\mycross{black} &\mycross{black} &\mycross{black} & &\mycross{black} &\mycross{black} \\
	\mycross{matlabred} &\mycross{black} &\mycross{black} & &\mycross{black}  &\mycross{black} &  & & \mycross{black} \\
	&&&\mycross{black} &&&\mycross{black} &\mycross{black} & \\
	&&&\mycross{black} &\mycross{black} & &\mycross{black} &\mycross{black} &\mycross{black} \\
	\mycross{matlabred}&&&\mycross{matlabred}&\mycross{black} &\mycross{black} & &\mycross{black} &\mycross{black} \\[0.15ex]
	\end{array}
	\normalsize
	\right],
\end{equation*}
which is obtained by a symbolic Cholesky factorization with approximately minimal degree ordering. \markendexample
\end{remark}

%Some notations

\subsubsection{Chordal decomposition of sparse block matrices}\label{ss:block-matrix-decomposition-theorems}
\noindent
As anticipated above, decomposition results similar to \cref{T:ChordalDecompositionTheorem,T:ChordalCompletionTheorem} hold for $\mathbb{S}^n_{\alpha,+}(\mathcal{E},?)$ and $\mathbb{S}^n_{\alpha,+}(\mathcal{E},0)$ when $\alpha$ is a chordal partition of $n$. Given a clique $\mathcal{C}_k$ of the chordal block sparsity graph $\mathcal{G}(\mathcal{V},\mathcal{E})$ subordinate to the chordal partition $\alpha$, we define the block matrix $E_{\mathcal{C}_k,\alpha} \in \mathbb{R}^{s(\alpha,k) \times n}$, where $s(\alpha,k) = \sum_{i \in \mathcal{C}_k}\alpha_i$, as
\begin{equation} \label{Eq:IndexMatrixBlock}
    (E_{\mathcal{C}_k,\alpha})_{ij} = \begin{cases} I_{\alpha_i}, & \text{if } \mathcal{C}_k(i) = j, \\ 0, & \text{otherwise}. \end{cases}
\end{equation}
Here, $I_{\alpha_i}$ is an identity matrix of dimension $\alpha_i$. %, and $\mathcal{C}_k$ is sorted with the natural ordering.   
When $\alpha = \{1,\ldots, 1\}$ is the trivial partition, $E_{\mathcal{C}_k,\alpha}$ reduces to the matrix $E_{\mathcal{C}_k}$ in~\cref{Eq:IndexMatrix}. 
% For a uniform partition $\alpha = \{a, a, \ldots, a\}$, we also have
% $
%     E_{\mathcal{C}_k,\alpha} = E_{\mathcal{C}_k} \otimes I_a, 
% $
% where $\otimes$ denotes the Kronecker product.
%
Similar to the case studied in~\cref{subsection:SparseMatrix}, the operation $E_{\mathcal{C}_k,\alpha}XE_{\mathcal{C}_k,\alpha}^\tr \in \mathbb{S}^{s(\alpha,k)}$ extracts the principal block-submatrix of $X$ whose blocks are indexed by $\mathcal{C}_k$, while $E_{\mathcal{C}_k,\alpha}^{\tr}YE_{\mathcal{C}_k,\alpha}$ ``inflates'' an $s(\alpha,k) \times s(\alpha,k)$ matrix into a sparse $n\times n$ block matrix.

We are now ready to %present the following two theorems~\cite[Theorems 2.17 \& 2.18]{zheng2019chordal} that 
extend~\cref{T:ChordalDecompositionTheorem,T:ChordalCompletionTheorem} to the case of sparse block matrices.

\begin{theorem} [Chordal block-decomposition]\label{T:GeneralDecompositionTheorem}
     Let $\mathcal{G}(\{1,\ldots,p\},\mathcal{E})$ be a chordal graph with maximal cliques $\mathcal{C}_1,\mathcal{C}_2, \ldots, \mathcal{C}_t$, and let $\alpha= \{\alpha_1,\ldots,\alpha_p\}$ be a partition of $n$. Then, $Z\in\mathbb{S}^n_{\alpha,+}(\mathcal{E},0)$ if and only if there exist matrices $Z_k \in \mathbb{S}^{s(\alpha,k)}_+$ for $k=1,\,\ldots,\,t$ such that
      \begin{equation} \label{E:DecompositionSparseCone}
        Z = \sum_{k=1}^{t} E_{\mathcal{C}_k,\alpha}^\tr Z_k E_{\mathcal{C}_k,\alpha}.
    \end{equation}
\end{theorem}

\begin{theorem} [Chordal block-completion]\label{T:GeneralCompletionTheorem}
     Let $\mathcal{G}(\{1,\ldots,p\},\mathcal{E})$ be a chordal graph with maximal cliques $\mathcal{C}_1,\mathcal{C}_2, \ldots, \mathcal{C}_t$, and let $\alpha= \{\alpha_1,\ldots,\alpha_p\}$ be a partition of $n$. Then, $X\in\mathbb{S}^n_{\alpha,+}(\mathcal{E},?)$ if and only if
     \begin{equation} \label{E:DecompositionCompletion}
        E_{\mathcal{C}_k,\alpha} X E_{\mathcal{C}_k,\alpha}^\tr \in \mathbb{S}^{s(\alpha,k)}_+ \quad \forall k=1,\,\ldots,\,t.
\end{equation}
\end{theorem}

The proofs of \cref{T:GeneralDecompositionTheorem,T:GeneralCompletionTheorem} rely on the fact that the block sparsity graph of $X \in \mathbb{S}^n_{\alpha}(\mathcal{E},0)$ induces a chordal extension of the standard sparsity graph of $X$ (cf. \Cref{remark:chordal-partition}) and, in fact, it is a \textit{hypergraph} of the latter. The normal chordal decomposition and completion from \cref{T:ChordalDecompositionTheorem,T:ChordalCompletionTheorem} can then be applied to the chordal extension of $X$, and the hypergraph structure implies the two results above. Interested readers are referred to~\citet[Chapter 2.4]{zheng2019chordal} for details.

%As discussed in~\Cref{remark:block_sparsity}, working with block-sparsity patterns directly can be advantageous. We can also ignore the sparsity in each block and detect maximal cliques in a small graph. The following example illustrates the block decomposition theorem. 

\begin{example}
Consider the $9 \times 9$ matrices
\begin{equation*}
    X = {\small\begin{pmatrix}
    2 & 1 & 0 & 1 & 1 & 0 & 0 & 0 & 0\\
    1 & 2 & 1 & 0 & 1 & 1 & 0 & 0 & 0\\
    0 & 1 & 2 & 0 & 0 & 1 & 0 & 0 & 0\\
    1 & 0 & 0 & 2 & 1 & 0 & 1 & 1 & 0\\
    1 & 1 & 0 & 1 & 2 & 1 & 0 & 1 & 1\\
    0 & 1 & 1 & 0 & 1 & 2 & 0 & 0 & 1\\
    0 & 0 & 0 & 1 & 0 & 0 & 2 & 1 & 0\\
    0 & 0 & 0 & 1 & 1 & 0 & 1 & 2 & 1\\
    0 & 0 & 0 & 0 & 1 & 1 & 0 & 1 & 2
    \end{pmatrix}}
\end{equation*}
and
\begin{equation*}
    Y = {\small\begin{pmatrix}
    1 & 1 & 0 & 1 & 1 & 0 & 0 & 0 & 0\\
    1 & 1 & 1 & 0 & 1 & 1 & 0 & 0 & 0\\
    0 & 1 & 1 & 0 & 0 & 1 & 0 & 0 & 0\\
    1 & 0 & 0 & 1 & 1 & 0 & 1 & 1 & 0\\
    1 & 1 & 0 & 1 & 1 & 1 & 0 & 1 & 1\\
    0 & 1 & 1 & 0 & 1 & 1 & 0 & 0 & 1\\
    0 & 0 & 0 & 1 & 0 & 0 & 1 & 1 & 0\\
    0 & 0 & 0 & 1 & 1 & 0 & 1 & 1 & 1\\
    0 & 0 & 0 & 0 & 1 & 1 & 0 & 1 & 1
    \end{pmatrix}},
\end{equation*}
which have the same nonchordal sparsity pattern as the symbolic matrix considered in \Cref{remark:chordal-partition}. Readers can easily check that $X$ is positive semidefinite, while $Y$ admits a positive semidefinite completion (e.g., replace all zero entries with ones to obtain $\overline{Y} = \boldsymbol{1}\boldsymbol{1}^\tr \succeq 0$). The partitions $\alpha_1=\{2,1,1,2,2,1\}$ and $\alpha_2 = \{3,3,3\}$ are both chordal, so while \cref{T:ChordalDecompositionTheorem} cannot be directly applied to decompose $X$, \cref{T:GeneralDecompositionTheorem} guarantees the existence of decompositions either in the symbolic form
\begin{multline*}
    \setlength\arraycolsep{0.5pt}
    \def\arraystretch{0.3}
    X = 
   \underbrace{\left[
	\small
	\begin{array}{cc|c|c|cc|cc|c}\\[0.1ex]
	\mycross{black} &\mycross{black} & \mycross{black} & \mycross{white} &\mycross{black} & \mycross{black} &\mycross{white} & \mycross{white}&\mycross{white}\\
	\mycross{black} &\mycross{black} & \mycross{black} & \mycross{white} &\mycross{black} & \mycross{black} &\mycross{white} & \mycross{white}&\mycross{white}\\\hline
	\mycross{black} &\mycross{black} & \mycross{black} & \mycross{white} &\mycross{black} & \mycross{black} & & &\\ \hline
	\mycross{white} &\mycross{white} & \mycross{white} & \mycross{white} &\mycross{white} & \mycross{white} &\mycross{white} & \mycross{white}&\mycross{white}\\\hline
	\mycross{black} &\mycross{black} & \mycross{black} & \mycross{white} &\mycross{black} & \mycross{black} & & &\\
	\mycross{black} &\mycross{black} & \mycross{black} & \mycross{white} &\mycross{black} & \mycross{black} & & &\\\hline
	\mycross{white}&&&&&&&&\\
	\mycross{white}&&&&&&&&\\\hline
	\mycross{white}&&&&&&&&\\[0.15ex]
	\end{array}
	\normalsize
	\right]}_{\succeq 0}
	+ 
    \underbrace{\left[
	\small
	\begin{array}{cc|c|c|cc|cc|c}\\[0.1ex]
	\mycross{black} &\mycross{black} & \mycross{white} & \mycross{black} &\mycross{black} & \mycross{black} &\mycross{white} & \mycross{white}&\mycross{white}\\
	\mycross{black} &\mycross{black} & \mycross{white} & \mycross{black} &\mycross{black} & \mycross{black} &\mycross{white} & \mycross{white}&\mycross{white}\\\hline
	\mycross{white} &\mycross{white} & \mycross{white} & \mycross{white} &\mycross{white} & \mycross{white} & & &\\ \hline
	\mycross{black} &\mycross{black} & \mycross{white} & \mycross{black} &\mycross{black} & \mycross{black} &\mycross{white} & \mycross{white}&\mycross{white}\\\hline
	\mycross{black} &\mycross{black} & \mycross{white} & \mycross{black} &\mycross{black} & \mycross{black} & & &\\
	\mycross{black} &\mycross{black} & \mycross{white} & \mycross{black} &\mycross{black} & \mycross{black} & & &\\\hline
	\mycross{white}&&&&&&&&\\
	\mycross{white}&&&&&&&&\\\hline
	\mycross{white}&&&&&&&&\\[0.15ex]
	\end{array}
	\normalsize
	\right]}_{\succeq 0}
	\\
	\setlength\arraycolsep{0.5pt}
    \def\arraystretch{0.3}
	+
	\underbrace{\left[
	\small
	\begin{array}{cc|c|c|cc|cc|c}\\[0.1ex]
	\mycross{white} &&&&&&&&\\
	\mycross{white} &&&&&&&&\\\hline
	\mycross{white} &&&&&&&&\\ \hline
	\mycross{white} &\mycross{white} & \mycross{white} & \mycross{black} &\mycross{black} & \mycross{black} &\mycross{black} & \mycross{black}&\mycross{white}\\\hline
	\mycross{white} &\mycross{white} & \mycross{white} & \mycross{black} &\mycross{black} & \mycross{black} & \mycross{black}& \mycross{black}&\\
	\mycross{white} &\mycross{white} & \mycross{white} & \mycross{black} &\mycross{black} & \mycross{black} & \mycross{black}& \mycross{black}&\\\hline
	\mycross{white} &\mycross{white} & \mycross{white} & \mycross{black} &\mycross{black} & \mycross{black} & \mycross{black}& \mycross{black}&\\
	\mycross{white} &\mycross{white} & \mycross{white} & \mycross{black} &\mycross{black} & \mycross{black} & \mycross{black}& \mycross{black}&\\\hline
	\mycross{white}&&&&&&&&\\[0.15ex]
	\end{array}
	\normalsize
	\right]}_{\succeq 0}
	+ 
    \underbrace{\left[
	\small
	\begin{array}{cc|c|c|cc|cc|c}\\[0.1ex]
	\mycross{white} &\mycross{white}&\mycross{white}&\mycross{white}&&&&&\\
	\mycross{white} &&&&&&&&\\\hline
	\mycross{white} &&&&&&&&\\\hline
	\mycross{white} &&&&&&&&\\\hline
	\mycross{white} &&&&\mycross{black} & \mycross{black} & \mycross{black}& \mycross{black}& \mycross{black}\\
	\mycross{white} &&&&\mycross{black} & \mycross{black} & \mycross{black}& \mycross{black}& \mycross{black}\\\hline
	\mycross{white} &&&&\mycross{black} &\mycross{black} & \mycross{black} & \mycross{black}& \mycross{black}\\
	\mycross{white} &&&&\mycross{black} & \mycross{black} & \mycross{black}& \mycross{black}& \mycross{black}\\\hline
	\mycross{white} &&&&\mycross{black} &\mycross{black} & \mycross{black} & \mycross{black}& \mycross{black}\\[0.15ex]
	\end{array}
	\normalsize
	\right]}_{\succeq 0}
\end{multline*}
(corresponding to the partition $\alpha_1$)
or in the symbolic form
\begin{equation*}
    X = \setlength\arraycolsep{0.5pt}
    \def\arraystretch{0.3}
   \left[
	\small
	\begin{array}{ccc|ccc|ccc}\\[0.1ex]
	\mycross{black} &\mycross{black} & \mycross{black} & \mycross{black} &\mycross{black} & \mycross{black} &\mycross{white} & \mycross{white}&\mycross{white}\\
	\mycross{black} &\mycross{black} & \mycross{black} & \mycross{black} &\mycross{black} & \mycross{black}\\
	\mycross{black} &\mycross{black} & \mycross{black} & \mycross{black} &\mycross{black} & \mycross{black}\\\hline
	\mycross{black} &\mycross{black} & \mycross{black} & \mycross{black} &\mycross{black} & \mycross{black}\\
	\mycross{black} &\mycross{black} & \mycross{black} & \mycross{black} &\mycross{black} & \mycross{black}\\
	\mycross{black} &\mycross{black} & \mycross{black} & \mycross{black} &\mycross{black} & \mycross{black}\\\hline
	\mycross{white}&&&&&&&&\\
	\mycross{white}&&&&&&&&\\
	\mycross{white}&&&&&&&&\\[0.15ex]
	\end{array}
	\normalsize
	\right]
	+ 
    \left[
	\small
	\begin{array}{ccc|ccc|ccc}\\[0.1ex]
    \mycross{white}&\mycross{white}&\mycross{white}&&&&&&\\
    \mycross{white}&\mycross{white}&\mycross{white}&&&&&&\\
    \mycross{white}&\mycross{white}&\mycross{white}&&&&&&\\\hline
    \mycross{white}&\mycross{white}&\mycross{white}&\mycross{black} &\mycross{black} & \mycross{black} & \mycross{black} &\mycross{black} & \mycross{black}\\
    &&&\mycross{black} &\mycross{black} & \mycross{black} & \mycross{black} &\mycross{black} & \mycross{black}\\
    &&&\mycross{black} &\mycross{black} & \mycross{black} & \mycross{black} &\mycross{black} & \mycross{black}\\\hline
    &&&\mycross{black} &\mycross{black} & \mycross{black} & \mycross{black} &\mycross{black} & \mycross{black}\\
    &&&\mycross{black} &\mycross{black} & \mycross{black} & \mycross{black} &\mycross{black} & \mycross{black}\\
    &&&\mycross{black} &\mycross{black} & \mycross{black} & \mycross{black} &\mycross{black} & \mycross{black}\\[0.15ex]
	\end{array}
	\normalsize
	\right]
\end{equation*}
(corresponding to the partition $\alpha_2$). These coincide with the classical chordal decompositions applied to the chordal extensions of $X$  from the chordal partitions $\alpha_1$ and $\alpha_2$. Thus, one can choose whether to decompose $X$ as a sum of four matrices with $5 \times 5$ nonzero principal submatrices, or as a sum of two matrices with $6 \times 6$ nonzero principal submatrices.
Similarly, one can apply \cref{T:GeneralCompletionTheorem} to verify that the matrix $Y$ admits a positive semidefinite completion by checking the positive semidefiniteness of either four $5\times 5$ principal submatrices, or two $6 \times 6$ ones.
\markendexample
\end{example}

%%%%%%%%%%%%%%%%%%%%%%%%%%%%%%%%%%%%%%%%%%%%%%%%%%%%
\section{Sparse semidefinite optimization}
\label{section:sparse-SDPs}
\noindent
The matrix decomposition and completion results in~\cref{T:ChordalDecompositionTheorem,T:ChordalCompletionTheorem} can be used to reduce the complexity of algorithms for sparse semidefinite optimization. A semidefinite program (SDP) in standard primal form takes the form
\begin{equation} \label{eq:SDP_primal}
    \begin{aligned}
        \min_{X} &\quad \langle C, X\rangle \\
        \text{subject to} & \quad \langle A_i, X \rangle = b_i, \quad i = 1, \ldots, m,\\
        &\quad X \in \mathbb{S}^n_+,
    \end{aligned}
\end{equation}
where $C, A_1, \ldots, A_m \in \mathbb{S}^n, b \in \mathbb{R}^m$ are the problem data. The dual problem to~\cref{eq:SDP_primal} is also an SDP,
\begin{equation} \label{eq:SDP_dual}
    \begin{aligned}
        \max_{y, Z} &\quad b^\tr y \\
        \text{subject to} & \quad Z + \sum_{i=1}^m y_i A_i = C,\\
        &\quad Z \in \mathbb{S}^n_+.
    \end{aligned}
\end{equation}

In this section, we describe decomposition techniques for SDPs that exploit the joint sparsity pattern of the coefficient matrices $C, A_1, \ldots, A_m$, called \textit{aggregate sparsity pattern}.
For simplicity, we assume that the matrices $A_1, \ldots, A_m$ are linearly independent and that there exist $X \succ 0$, $y \in \mathbb{R}$ and $Z \succ 0$ satisfying the equality constraints in~\cref{eq:SDP_primal,eq:SDP_dual}. This ensures that the primal and dual optimal values are finite, equal, and attained.
SDPs that are infeasible or have unbounded objective can be tackled using homogeneous self-dual embeddings \citep{ye1994nl,o2016conic,ye2011interior} or by analyzing the divergence of the iterates produced by solution algorithms \citep{liu2017new,banjac2019infeasibility}. Sparsity can be exploited within these frameworks, too, and we refer the interested reader to~\citet[Section 5]{ZFPGW2020chordal} and \cite{garstka2019cosmo} for details.

%%%%%%%%%%%%%%%%%%%%%%%%%%%%%%%%%%%%%%%%%%%%%%%%%%%%%%%%%%%%%%%%%%%%%%%%%%%%%%
\subsection{Aggregate sparsity} \label{subsection:aggregate_sparsity_pattern}
\noindent
The pair of SDPs~\cref{eq:SDP_primal}-\cref{eq:SDP_dual} is said to have aggregate sparsity graph $\mathcal{G}(\mathcal{V},\mathcal{E})$ if
\begin{equation} \label{eq:AggregateSparsity}
    C, A_1, \ldots, A_m \in \mathbb{S}^n(\mathcal{E},0). 
\end{equation}
Of course, if $\mathcal{E}'$ is an extension of $\mathcal{E}$, then $\mathcal{G}(\mathcal{V},\mathcal{E}')$ is also a suitable aggregate sparsity graph. The minimal one, therefore, is simply the union of the individual sparsity graphs of $C$, $A_1,\,\ldots,\,A_m$. Throughout this section, however, we consider a chordal extension of the minimal aggregate sparsity graph. We therefore assume from now on that the aggregate sparsity pattern $\mathcal{E}$ is chordal and has $t$ maximal cliques $\mathcal{C}_1,\,\ldots,\,\mathcal{C}_t$. 

It must be noted that an SDP may have a fully connected aggregate sparsity graph even if all coefficient matrices $C$ and $A_1,\ldots,A_m$ are very sparse; see \cite{zheng2018fast} for explicit examples. The decomposition methods described below cannot be applied to such problems. However, there are broad classes of SDPs for which the sparsity of the SDP data matrices can be expected to translate into very sparse aggregate sparsity graphs.

One such family consists of SDPs arising from relaxations of graph optimization problems and control problems over networks, which typically inherit the structure of the underlying network or graph. Notable examples include 
% "Abstract" applications
SDP relaxations of combinatorial graph optimization problems, such as Max-Cut \citep{goemans1995improved} and graph equipartition \citep{karisch1998semidefinite}, 
eigenvalue optimization problems over graphs \citep{boyd2004fastest},
% Control problems
analysis of linear networked systems \citep{mason2014chordal, ZMP2018Scalable,deroo2015distributed,ZKSP2018scalable}, 
sensor network localization \citep{kim2009exploiting,nie2009sum,so2007theory},  
neural network verification in machine learning \citep{raghunathan2018semidefinite,BPLZ2021neural}, 
and the optimal power flow problem in electricity networks \citep{bai2008semidefinite,andersen2014reduced,jabr2011exploiting}. 
We will briefly discuss some of these applications in~\Cref{section:applications}. 

Another source of SDPs with aggregate sparsity is the reformulation of intractable constraints (either convex or nonconvex) as tractable LMIs using auxiliary variables \citep{ben2001lectures,vandenberghe2015chordal}. For example, consider the uncountable family of ``uncertain'' convex quadratic constraints
$$
    x^\tr A^\tr Ax - 2b^\tr x -c \leq 0
$$
on a variable $x \in \mathbb{R}^q$, to be imposed for all matrices $A \in \mathbb{R}^{p\times q}$, vectors $b \in \mathbb{R}^q $ and scalars $c \in \mathbb{R}$ in the form
\begin{align*}
    A &= A_0 + \sum_{i=0}^r u_i A_i,&
    b &= b_0 + \sum_{i=0}^r u_i b_i,&
    c &= c_0 + \sum_{i=1}^r u_i c_i
\end{align*}
with $u^\tr u \leq 1$. Here, $A_0$, $b_0$ and $c_0$ are nominal reference values, and $\{A_i, b_i, c_i\}$ are fixed perturbations. \cite{andersen2010linear} showed that this family of constraints is equivalent to a sparse LMI in the form
\begin{equation} \label{eq:robust_constraint}
    \begin{bmatrix} 
    f(x)-t & (A_0x)^\tr & h(x)^\tr\\
    A_0x & I_p & G(x)^\tr\\
    h(x) & G(x) & t I_r
    \end{bmatrix} 
    \succeq 0,
\end{equation}
where $t \in \mathbb{R}$ while  $G:\mathbb{R}^q\to \mathbb{R}^{r\times p}$, $h:\mathbb{R}^q \to \mathbb{R}^r$ and $f:\mathbb{R}^q \to \mathbb{R}$ are known linear functions whose exact form is not important here. % \citep{andersen2010linear}. 
When $r\gg p$, this matrix has a ``block-arrow'' aggregate sparsity pattern analogous to that shown in \cref{fig:sparsematrix_b} (that figure is recovered exactly when $p=1$, $r=6$, and $q$ is arbitrary). This particular type of sparsity pattern is commonly encountered in robust optimization \citep{goldfarb2003robust,ben1998robust,andersen2010linear}.

\begin{remark}[Promoting aggregate sparsity]
    Sometimes, it is possible to reformulate SDPs with no aggregate sparsity as equivalent SDPs with very sparse aggregate sparsity graphs through a carefully chosen transformation of variables~(\citealp[Section 6]{fukuda2001exploiting}; \citealp[Chapter 14.1]{vandenberghe2015chordal}). For instance, the SDP relaxation of the graph equipartition problem studied by~\citet[Section 6]{fukuda2001exploiting} has sparse data matrices $C$ and $A_1,\ldots,A_{m-1}$, but the $m$-th constraint $\langle11^\tr, X\rangle = 0$ destroys the problem's aggregate sparsity because the matrix $A_m = 1 1^\tr$ is dense. However, any matrix $X \in \mathbb{S}^n_+$ satisfying $\langle11^\tr, X\rangle = 0$ can be expressed as $X = VYV^\tr$ for some matrix $Y \in \mathbb{S}^{n-1}_+$, where %$V$ being any $n \times (n - 1)$-matrix with columns that 
    $$
    {\def\arraystretch{0.8}
        \setlength\arraycolsep{3pt}
        V = \begin{bmatrix} 1 & 0 & 0 & \ldots & 0 & 0 \\
                            -1 & 1 & 0 & \ldots & 0 & 0 \\
                            0  & -1 & 1 & \ldots & 0 & 0 \\
                            \vdots & \vdots & \vdots& & \vdots& \vdots\\
                            0 & 0 & 0 & \ldots & -1 & 1 \\
                            0 & 0 & 0 & \ldots & 0 & -1\end{bmatrix} \in \mathbb{R}^{n\times (n - 1)}.
    }$$
    %By eliminating the constraint $\langle11^\tr, X\rangle = 0$ using the variable transformation $X = VYV^\tr$, we can reformulate 
    Thus, the original SDP can be reformulated as
    $$
        \begin{aligned}
            \min_{Y} &\quad \langle C', Y \rangle \\
            \text{subject to} & \quad \langle A_i', Y\rangle = b_i, \quad i = 1, \ldots, m-1,\\
            &\quad Y \in \mathbb{S}^{n-1}_+,
        \end{aligned}
    $$
    where $C' := V^\tr C V$ and $A_i' := V^\tr A_i V$.
    Since $V$ is a sparse basis matrix and the original data matrices are sparse, this new SDP is characterized by aggregate sparsity~\cite[Section 6]{fukuda2001exploiting}. Sparsity-promoting modeling strategies that generalize this example are discussed by \citet[Chapter 14.1]{vandenberghe2015chordal}.  
\end{remark}

%%%%%%%%%%%%%%%%%%%%%%%%%%%%%%%%%%%%%%%%%%%%%%%%%%%%%%%%%%%%%%%%%%%%%%%%%%%%%
\begin{table*}[t]
     \setlength{\abovecaptionskip}{0mm}
     \setlength{\belowcaptionskip}{0mm}
     \renewcommand\arraystretch{1.0}
     \caption{Comparison of first-order algorithms for solving SDPs. {``Chordal Sparsity'': whether the algorithm exploits chordal sparsity; ``SDP Type'': the types of SDP problems the algorithm considers; ``Algorithm'': the underlying first-order algorithm; ``infeas./unbounded'': whether the algorithm can detect infeasible or unbounded cases; ``Solver'': whether the code is open-source.}} 
     \begin{center}
     \small
  \begin{tabular}{c c c c c c }%{c C{15mm} C{20mm} C{24mm} C{18mm} c }
  \toprule
   \multicolumn{1}{c}{Reference}   & \makecell{Chordal Sparsity} &  SDP Type  & Algorithm  &   \makecell{Infeas./ Unbounded }    &  Solver \\%\makecell{open-source \\ solver}  \\
    \hline
    \cite{wen2010alternating} &  {\xmark} &  \cref{eq:SDP_dual} & ADMM & \xmark & \xmark \\
    \cite{zhao2010newton} & \xmark &\cref{eq:SDP_dual} & %Augmented Lagrangian 
    Augm. Lagrang. & \xmark & SDPNAL\\
    \cite{o2016conic} & \xmark & \cref{eq:SDP_primal}-\cref{eq:SDP_dual}  & ADMM & \checkmark & SCS \\
    \cite{yurtsever2021scalable}& \xmark & \cref{eq:SDP_primal}$^{1}$ & SketchyCGAL & \xmark & CGAL\\
    \hline
    \cite{lu2007large} & \checkmark &  \cref{eq:SDP_primal}  & Mirror-Prox & \xmark &\xmark \\
    \cite{lam2012distributed} & \checkmark &  OPF$^{2}$ & Primal-dual & \xmark &\xmark \\
    \cite{dall2013distributed} & \checkmark & OPF$^{2}$  & ADMM & \xmark &\xmark \\
    \cite{sun2014decomposition} & \checkmark & Special$^{3}$  &  Gradient proj. & \xmark &\xmark \\
    \cite{sun2015decomposition} & \checkmark &  \cref{eq:SDP_primal}-\cref{eq:SDP_dual} & Spingarn%'s method 
    & \xmark &\xmark \\
    \cite{kalbat2015fast} & \checkmark & Special$^{4}$  & ADMM & \xmark& \xmark \\
    \cite{madani2017low} & \checkmark &  General$^{5}$  &  ADMM & \xmark &\xmark \\
    \cite{ZFPGW2020chordal} &\checkmark  & \cref{eq:SDP_primal}-\cref{eq:SDP_dual} & ADMM & \checkmark& CDCS \\
    \cite{garstka2019cosmo} & \checkmark & Quad. SDP$^{6}$ & ADMM  & \checkmark& COSMO \\
    \bottomrule
  \end{tabular}
  \end{center}
  \raggedright
 {\footnotesize
 \begin{spacing}{1}
 Note: 1. It requires an explicit trace constraint on $X$; 2. Special SDPs from the optimal power flow (OPF) problem; 3. Special SDPs from the matrix nearness problem; 4. Special SDPs with decoupled affine constraints; 5. General SDPs with inequality constraints; 6. A dual SDP~\cref{eq:SDP_dual} with a quadratic objective function.
 \end{spacing}
}
\label{tab:comparision_SDP}
\end{table*}
%%%%%%%%%%%%%%%%%%%%%%%%%%%%%%%%%%%%%%%%%%%%%%%%%%%%%%%%%%%%%%%%%%%%%%%%%%%%%%%%%%

%%%%%%%%%%%%%%%%%%%%%%%%%%%%%%%%%%%%%%%%%%%%%%%%%%%%%%%%%%%%%%%%%%%%%%%
\subsection{Nonsymmetric formulation} \label{subsection:nonsymmetric}
\noindent
The aggregate sparsity of the primal-dual pair of SDPs~\cref{eq:SDP_primal}--\cref{eq:SDP_dual} can be exploited by reformulating them into a nonsymmetric pair of optimization problems, proposed by~\cite{fukuda2001exploiting} and later discussed extensively by~\cite{andersen2010implementation, sun2014decomposition, kim2011exploiting, ZFPGW2020chordal}.   

Consider first the dual-standard-form SDP~\cref{eq:SDP_dual}. Any feasible matrix $Z$ must be at least as sparse as the aggregate sparsity pattern of the SDP. We can therefore restrict $Z$ to the subspace $\mathbb{S}^n(\mathcal{E},0)$, where $\mathcal{E}$ is the edge set of the aggregate sparsity graph, and rewrite~\cref{eq:SDP_dual} as
\begin{equation}\label{eq:NonSymmetric_dual}
  \begin{aligned}
            \max_{y, Z} \quad & \langle b,y\rangle \\
            \text{subject to} \quad  & Z + \sum_{i=1}^m A_i\,y_i = C,\\
            & Z \in \mathbb{S}^n_{+}(\mathcal{E},0).
        \end{aligned}
\end{equation}

The primal-standard-form SDP~\cref{eq:SDP_primal}, instead, typically has a dense optimal matrix $X$. However, the value of the cost function and the equality constraints depend only on the entries $X_{ij}$ with %$i = j$ or 
$(i,j)\in\mathcal{E}$, while the remaining ones simply guarantee that {$X$ is positive semidefinite}. We can therefore pose~\cref{eq:SDP_primal} as an optimization problem over the cone $\mathbb{S}^n_{+}(\mathcal{E},?)$ of sparse matrix that admit a positive semidefinite completion,
\begin{equation}\label{eq:NonSymmetric_primal}
    \begin{aligned}
            \min_{X} \quad & \langle C, X \rangle \\
            \text{subject to} \quad & \langle A_i, X\rangle = b_i,
            \quad i = 1,\,\ldots,\,m, \\
                & X \in \mathbb{S}^n_{+}(\mathcal{E},?).
            \end{aligned}
\end{equation}

Problems~\cref{eq:NonSymmetric_primal} and~\cref{eq:NonSymmetric_dual} are a primal-dual pair of linear conic programs because the cones $\mathbb{S}^n_+(\mathcal{E}, ?)$ and $\mathbb{S}^n_+(\mathcal{E}, 0)$ are dual to each other (see~\cref{subsection:SparseMatrix} and~\cref{Fig:SparseCone_Duality}). Even though the sparse matrix cones $\mathbb{S}^n_+(\mathcal{E}, ?)$ and $\mathbb{S}^n_+(\mathcal{E}, 0)$ are not self-dual \citep{andersen2010implementation,andersen2011chordal}, so this sparse formulation is nonsymmetric, one can solve~\cref{eq:NonSymmetric_primal}, \cref{eq:NonSymmetric_dual}, or both problems simultaneously using a variety of first-order or interior-point algorithms. The next two subsections discuss some of them.

\begin{remark}
    A special type of aggregate sparsity arises when the data matrices $C,A_1,\,\ldots,\,A_m$ are block-diagonal. 
    In this case, any feasible matrix $X$ for~\cref{eq:NonSymmetric_primal} is automatically positive semidefinite and, consequently, can be restricted to $\mathbb{S}^n_+(\mathcal{E}, 0)$. Therefore, the nonsymmetric formulation described above becomes symmetric. In particular, problems~\cref{eq:NonSymmetric_primal,eq:NonSymmetric_dual} are simply SDPs with a Cartesian product $\mathbb{S}^{n_1}_+ \times \mathbb{S}^{n_2}_+ \times \ldots \times \mathbb{S}^{n_l}_+$ of semidefinite cones, where $n_i$ is the size of the $i$th diagonal block and $l$ is the number of blocks.  \markendexample
\end{remark}

%%%%%%%%%%%%%%%%%%%%%%%%%%%%%%%%%%%%%%%%%%%%%%%%%%%%%%%%%%%%%%%%%%%%%%%%%%%
\subsection{First-order algorithms} \label{subsection:first-order-methods}  
\noindent
First-order optimization algorithms rely only on gradient information and have iterations with low computational complexity, which can often  be implemented in a distributed manner \citep{boyd2011distributed,beck2017first}. For these reasons, the last decade has witnessed the development of a range of first-order methods to solve large-scale SDPs, many of which are listed in \cref{tab:comparision_SDP}. Some of these methods \citep{wen2010alternating, zhao2010newton, o2016conic, yurtsever2021scalable} focus on generic SDPs and do not exploit aggregate sparsity. Others, instead, tackle the sparsity-exploiting nonsymmetric formulations~\cref{eq:NonSymmetric_dual}--\cref{eq:NonSymmetric_primal} using so-called \textit{domain space} or \textit{range-space} conversion frameworks, which replace the matrix cones $\mathbb{S}^n_+(\mathcal{E}, ?)$ and $\mathbb{S}^n_+(\mathcal{E}, 0)$ with smaller positive semidefinite cones using the chordal decomposition and completion results in~\cref{T:ChordalDecompositionTheorem,T:ChordalCompletionTheorem}  (see, e.g., \citealp{lu2007large,lam2012distributed,dall2013distributed,sun2014decomposition,sun2015decomposition,kalbat2015fast,madani2017low,ZFPGW2020chordal,garstka2019cosmo}). Many of these works combine this strategy with additional separability assumptions for the equality constraints, which are satisfied in optimal power flow problems \citep{lam2012distributed,dall2013distributed,kalbat2015fast} and the matrix nearness problems \citep{sun2014decomposition} but not in general. To the best of our knowledge, the only first-order methods that can currently handle general SDPs with aggregate sparsity (including infeasible or unbounded ones) are those developed by~\cite{ZFPGW2020chordal} and \cite{garstka2019cosmo}.

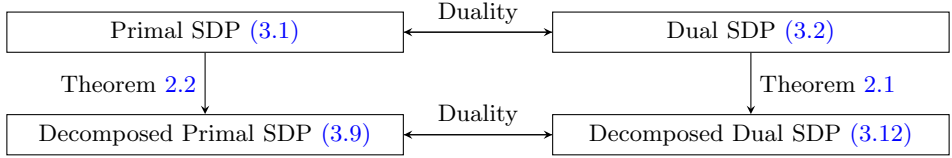
\begin{figure*}
    \centering
    \setlength{\abovecaptionskip}{6pt}
	\small
	\begin{tikzpicture}
	  \matrix (m) [matrix of nodes,
	  		       row sep = 2.5em,	
	  		       column sep = 6em,	
  			       nodes={rectangle, draw=black, align=center, text width=5cm}]
  	{
   	Primal SDP~\cref{eq:SDP_primal} & Dual SDP~\cref{eq:SDP_dual} \\
   	Decomposed Primal SDP~\cref{eq:DecomposedPrimalSDP} &
   	Decomposed Dual SDP~\cref{eq:DecomposedDualSDP} \\
   	%Algorithm~\ref{A:ADMMPrimal} & Algorithm~\ref{A:ADMMDual}\\
   	};
	\path[-stealth]
		(m-1-1) edge node [left, align=center]
			%{\scriptsize Grone's \scriptsize Theorem}  (m-2-1)
            { Theorem~\ref{T:ChordalCompletionTheorem}}  (m-2-1)
		(m-1-2) edge node [right, align=center]
			{ Theorem~\ref{T:ChordalDecompositionTheorem}}  (m-2-2)
		(m-1-1) edge node [above] { Duality} (m-1-2)
        (m-1-2) edge (m-1-1)
		(m-2-1) edge node [above] { Duality} (m-2-2)
        (m-2-2) edge (m-2-1);
	\end{tikzpicture}
    \caption{Duality between the original primal and dual SDPs, and the decomposed primal and dual SDPs.}
    \label{F:Duality}
\end{figure*}

%%%%%%%%%%%%%%%%%%%%%%%%%%%%%%%%%%%%%%%%%%%%%%%%%%%%%%%%%%%%%%%%%%%%%%%%%%%%%%%%%%%%%%%%%%%%
\subsubsection{Domain- and range-space conversion} \label{subsection:conversion_first_order}
\noindent
Consider problem~\cref{eq:NonSymmetric_primal}. When the aggregate sparsity graph is chordal and has maximal cliques $\mathcal{C}_1,\ldots,\mathcal{C}_t$, \cref{T:ChordalCompletionTheorem} allows one to replace the constraint $X \in \mathbb{S}^n_{+}(\mathcal{E},?)$ with
\begin{equation} \label{eq:coupling_cliques}
    E_{\mathcal{C}_k}XE_{\mathcal{C}_k}^\tr
    \in \mathbb{S}^{|\mathcal{C}_k|}_+
    \quad \forall k = 1,\, \ldots,\, t.
\end{equation}
These constraints are coupled in general because the matrices $E_{\mathcal{C}_p}XE_{\mathcal{C}_p}^\tr $ and $E_{\mathcal{C}_q}XE_{\mathcal{C}_q}^\tr $ depend on the same entries of $X$ if the cliques $\mathcal{C}_p$ and $\mathcal{C}_q$ overlap. The works referenced above differ primarily in how these couplings are handled and, as discussed in \cref{remark:conversion_first_order_method,remark:slack_varaibels_iterations} below, the choice of strategy can have a considerable impact on the overall complexity of the iterations in a first-order method. 

A simple but powerful strategy was proposed recently by~\cite{ZFPGW2020chordal}, who used ``slack'' matrices $X_1,\ldots,X_t$ to rewrite~\cref{eq:coupling_cliques} as
\begin{equation} \label{eq:PrimalSlackVariables}
    \begin{cases}
    		X_k = E_{\mathcal{C}_k}XE_{\mathcal{C}_k}^\tr  &\forall k = 1,\,\ldots,\,t,\\
    		X_k \in \mathbb{S}^{|\mathcal{C}_k|}_+  &\forall k = 1,\,\ldots,\,t.
    \end{cases}
\end{equation}
The primal SDP~\cref{eq:NonSymmetric_primal} is then equivalent to
\begin{equation} \label{eq:DecomposedPrimalSDP}
    \begin{aligned}
        \min_{X,X_1,\ldots,X_t} \quad & \langle C, X \rangle \\
        \text{subject to} \quad & \langle A_i, X\rangle = b_i, &i = 1, \ldots, m, \\
            & X_k = E_{\mathcal{C}_k}XE_{\mathcal{C}_k}^\tr ,  &k = 1, \ldots, t, \\
            & X_k \in \mathbb{S}^{|\mathcal{C}_k|}_+,  &k = 1, \ldots, t.
     \end{aligned}
\end{equation}
Following~\cite{fukuda2001exploiting} and \cite{ZFPGW2020chordal}, we refer to~\cref{eq:DecomposedPrimalSDP} as the \emph{domain-space} decomposition of the primal SDP~\cref{eq:SDP_primal}. 

A \emph{range-space} decomposition of the dual SDP~\cref{eq:SDP_dual} can be formulated in a very similar way. When the aggregate sparsity pattern $\mathcal{E}$ is chordal, \cref{T:ChordalDecompositionTheorem} implies that the constraint $Z\in \mathbb{S}^n_{+}(\mathcal{E},0)$ is equivalent to
\begin{equation}\label{e:range-space-constraint-initial}
    \begin{dcases}
    Z = \sum_{k=1}^{t} E_{\mathcal{C}_k}^\tr  Z_k E_{\mathcal{C}_k},\\[1ex]
    Z_k \in \mathbb{S}^{|\mathcal{C}_k|}_+  &\forall k = 1,\,\ldots,\,t.
    \end{dcases}
\end{equation}
Observe that, as before, the first of these conditions couples the positive semidefinite matrices $Z_p$ and $Z_q$ if the cliques $\mathcal{C}_p$ and $\mathcal{C}_q$ of the aggregate sparsity graph overlap. To decouple them, \cite{ZFPGW2020chordal} introduced slack variables $V_1,\ldots,V_t$ and reformulated~\cref{e:range-space-constraint-initial} as
\begin{equation}
    \begin{dcases}
    Z = \sum_{k=1}^{t} E_{\mathcal{C}_k}^\tr  V_k E_{\mathcal{C}_k},\\[1ex]
    V_k = Z_k & \forall k = 1,\,\ldots,\,t,\\
    Z_k \in \mathbb{S}^{|\mathcal{C}_k|}_+  &\forall k = 1,\,\ldots,\,t.
    \end{dcases}
\end{equation}
Using this to eliminate $Z$ from~\cref{eq:NonSymmetric_dual} yields the range-space decomposition
\begin{equation} \label{eq:DecomposedDualSDP}
    \begin{aligned}
            \max_{y,Z_1,\ldots, Z_t, V_1,\ldots,V_t} \quad & \langle b,y\rangle \\
            \text{subject to} \quad  & \sum_{i=1}^m A_i\,y_i + \sum_{k=1}^{t} E_{\mathcal{C}_k}^\tr  V_k E_{\mathcal{C}_k}= C,\\
            & Z_k - V_k = 0, \; k = 1, \ldots, t, \\
            & Z_k \in  \mathbb{S}^{|\mathcal{C}_k|}_+, \quad \; k = 1, \ldots, t.
        \end{aligned}
\end{equation}

While the domain- and range-space decompositions \cref{eq:DecomposedPrimalSDP} and \cref{eq:DecomposedDualSDP} have been derived independently, it is not difficult to verify that they are a primal-dual pair of SDPs. The duality between the original SDPs~\cref{eq:SDP_primal} and~\cref{eq:SDP_dual} is thus inherited by the decomposed SDPs~\cref{eq:DecomposedPrimalSDP} and~\cref{eq:DecomposedDualSDP} by virtue of the duality between~\Cref{T:ChordalCompletionTheorem} and~\Cref{T:ChordalDecompositionTheorem}. This elegant picture is illustrated in~\cref{F:Duality}.

\begin{remark} \label{remark:conversion_first_order_method}
The introduction of variables $X_k$ and $V_k$ leads to redundancies in the affine constraints of~\cref{eq:DecomposedPrimalSDP,eq:DecomposedDualSDP}, but is essential to obtain a decomposition framework that is suitable for the development of fast first-order SDP solvers. For example, as explained in~\cref{subsection:ADMM} below, applying the alternating direction method of multipliers (ADMM) to~\cref{eq:DecomposedPrimalSDP} leads to an algorithm whose iterations have closed-form update rules that can be implemented efficiently. The same is usually not true if the redundant constraints in~\cref{eq:DecomposedPrimalSDP} are used to eliminate the matrix $X$: the iterations of the first-order method proposed by~\cite{sun2014decomposition}, for instance, require the solution of a further SDP with quadratic objective function, which limits its scalability. 
However, the matrix $X$ may be eliminated from~\cref{eq:DecomposedPrimalSDP} without compromising efficiency if the original primal SDP~\cref{eq:SDP_primal} has separable affine constraints. This observation was exploited to solve sparse SDPs arising from optimal power flow problems \citep{kalbat2015fast,dall2013distributed} and matrix nearness problems \citep{sun2015decomposition}. 
Similar observations apply to the seemingly redundant matrices $V_k$ in the range-space decomposed SDP~\cref{eq:DecomposedDualSDP}.  \markendexample
\end{remark}

%%%%%%%%%%%%%%%%%%%%%%%%%%%%%%%%%%%%%%%%%%%%%%%%%%%%%%%%%%%%%%%%%%%%%%%%%%%
\subsubsection{ADMM for decomposed SDPs} \label{subsection:ADMM}
\noindent
The alternating direction method of multipliers (ADMM) is a first-order operator-splitting method developed in the mid-1970s \citep{gabay1976dual,glowinski1975approximation} %\red{[Give original references?]} 
to solve general optimization problems in the form
\begin{equation}
\label{e:ADMM_base_prob}
    \begin{aligned}
        \min_{\substack{\mathcal{X} \in \mathbb{X}\\\mathcal{Y} \in \mathbb{Y}}} \quad & f(\mathcal{X})+g(\mathcal{Y}) \\
        \text{subject to} \quad & \mathcal{A}(\mathcal{X}) + \mathcal{B}(\mathcal{Y}) = \mathcal{C},
    \end{aligned}
\end{equation}
where $f$ and $g$ are proper convex (but not necessarily smooth) functions on finite-dimensional normed vector spaces $\mathbb{X}$ and $\mathbb{Y}$,  $\mathcal{A}$ and $\mathcal{B}$ are given linear operators from $\mathbb{X}$ and $\mathbb{Y}$ into a finite-dimensional normed vector space $\mathbb{Z}$, and $\mathcal{C} \in \mathbb{Z}$ is given. 
Given a penalty parameter $\rho>0$ and a dual variable $\mathcal{Z} \in \mathbb{Z}$ that acts as a Lagrange multilier for the equality constraint, ADMM finds a saddle point of the (scaled) augmented Lagrangian %\red{[the scaled formulation is much more convenient]}
\begin{equation*}
\mathcal{L}_{\rho}(\mathcal{X},\mathcal{Y},\mathcal{Z}) := f(\mathcal{X}) + g(\mathcal{Y}) + \frac{\rho}{2} \left\|\mathcal{A}(\mathcal{X}) + \mathcal{B}(\mathcal{Y}) - \mathcal{C} + \mathcal{Z}\right\|^2
\end{equation*}
by updating the primal variables $\mathcal{X}$, $\mathcal{Y}$ and the dual variable $\mathcal{Z}$ according to the following rules:
\begin{subequations}\label{eq:ADMM}
    \begin{align}
        \mathcal{X}^{(q+1)} & = \text{arg} \min_{\mathcal{X}} \; \mathcal{L}_{\rho}(\mathcal{X},\mathcal{Y}^{(q)},\mathcal{Z}^{(q)}),
        \label{eq:ADMM_S1}\\
        \mathcal{Y}^{(q+1)} & = \text{arg} \min_{\mathcal{Y}} \; \mathcal{L}_{\rho}(\mathcal{X}^{(q+1)},\mathcal{Y},\mathcal{Z}^{(q)}),
        \label{eq:ADMM_S2}\\
        \mathcal{Z}^{(q+1)} &= \mathcal{Z}^{(n)} + \mathcal{A}(\mathcal{X}^{(q+1)}) + \mathcal{B}(\mathcal{Y}^{(q+1)}) - \mathcal{C}. \label{eq:ADMM_S3}
    \end{align}
\end{subequations}
The superscript $(q)$ indicates that a variable is fixed to its value at the $q$-th iteration. Under mild technical conditions \cite[Section 3.2]{boyd2011distributed}, the method converges to  %solution
an $\epsilon$-approximate solution of~\cref{e:ADMM_base_prob} using at most $\mathcal{O}(1/\epsilon)$ iterations.  
%with a rate $\mathcal{O}(\frac{1}{q})$. 

Given its slow convergence rate, ADMM is suitable only when~\cref{eq:ADMM_S1} and~\cref{eq:ADMM_S2} have closed-form expressions and/or can be solved efficiently. Below, we show that this is true when the method is applied to the decomposed SDPs~\cref{eq:DecomposedPrimalSDP,eq:DecomposedDualSDP}.

\paragraph{Domain-space decomposition}
Consider the domain-space decomposition~\cref{eq:DecomposedPrimalSDP}. Let $\chi_{\mathcal{K}}(x)$ denote the characteristic function of a set $\mathcal{K}$, {\it i.e.},
$$
    \chi_{\mathcal{K}}(x) :=
    \begin{cases}
    0 &\text{if } x \in \mathcal{K}, \\
    + \infty &\text{otherwise}.
    \end{cases}
$$
For simplicity, we write $\chi_{0}$ when $\mathcal{K}\equiv \{0\}$. Problem~\cref{eq:DecomposedPrimalSDP} is equivalent to
\begin{subequations}
    \begin{align}
    \min_{X,X_1,\ldots,X_t} \quad &\langle C,X \rangle + 
    \sum_{i=1}^m \left[ \chi_0\left( \langle A_i, X \rangle -b_i \right) + \chi_{\mathbb{S}^{|\mathcal{C}_k|}_+}(X_k) \right]\nonumber
    \\
    \text{subject to} \quad &X_k = E_{\mathcal{C}_k}XE_{\mathcal{C}_k}^\tr,  \quad k = 1, \ldots, t. \nonumber  %\label{eq:ADMMPrimal} 
    \end{align}
\end{subequations}

Upon letting $\mathcal{X}:= \{X\}$ and $\mathcal{Y}:= \{X_1,\,\ldots,\,X_t\}$, this problem may be written in the standard form~\cref{e:ADMM_base_prob} over the spaces $\mathbb{X}=\mathbb{S}^n$ and $\mathbb{Y}=\mathbb{Z}=\mathbb{S}^{\abs{\mathcal{C}_1}} \times \cdots \times \mathbb{S}^{\abs{\mathcal{C}_t}}$, and can therefore be solved using ADMM. 
Introducing a penalty parameter $\rho>0$ and a dual variable $\mathcal{Z}:= \{\Lambda_1,\,\ldots,\,\Lambda_t\}$, where each $\Lambda_k \in \mathbb{S}^{\abs{\mathcal{C}_k}}$ acts as a Lagrange multiplier for the corresponding constraint $X_k = E_{\mathcal{C}_k}XE_{\mathcal{C}_k}^\tr$, it is not difficult to check that the ADMM step~\cref{eq:ADMM_S1} reduces to an equality-constrained quadratic program,
%
% \begin{subequations}
    \begin{multline*}
        X^{(q+1)} = \argmin_{\substack{\langle A_i,X\rangle = b_i\\ i=1,\ldots,m}}
        \bigg\{ \frac{\rho}{2} \sum_{k=1}^{t} 
        \left\| X_k^{(q)} - E_{\mathcal{C}_k}XE_{\mathcal{C}_k}^\tr + \Lambda_k^{(q)} \right\|_F^2
        \\[-2ex]
        +  \langle C,X \rangle \bigg\}.
        % \label{eq:MinXblockPrimal}
    \end{multline*}
Step~\cref{eq:ADMM_S2}, instead, reduces to $t$ independent positive semidefinite projections of the form for $k = 1, \ldots, t$
\begin{equation*}
    X_k^{(q+1)} = \argmin_{X_k \in \mathbb{S}_+^{|\mathcal{C}_k|} } \left\| X_k - E_{\mathcal{C}_k}X^{(q+1)}E_{\mathcal{C}_k}^\tr  + \Lambda_k^{(q)} \right\|_F^2.
    % \label{eq:MinYblockPrimal}
\end{equation*}
Finally, step~\cref{eq:ADMM_S3} updates the multipliers $\Lambda_1,\ldots,\Lambda_t$ according to
\begin{equation*}
% \label{eq:LambdakUpdatePrimal}
    \Lambda_k^{(q+1)} =
    \Lambda_k^{(q)} + X_k^{(q+1)} - E_{\mathcal{C}_k} X^{(q+1)} E_{\mathcal{C}_k}^\tr. %, %k = 1, \ldots, t.
\end{equation*}
% \end{subequations}

These three steps have efficient closed-form solutions and can be implemented efficiently~\cite[Section 4.1]{ZFPGW2020chordal}. In particular, the $t$ independent projections onto the cones $\mathbb{S}^{\vert \mathcal{C}_k\vert}_{+}$ required to compute $X_1^{(q+1)},\ldots,X_t^{(q+1)}$ can be computed through an eigenvalue decomposition with complexity of $O(\abs{\mathcal{C}_k}^3)$ floating-point operations. This is not expensive when all cliques $\mathcal{C}_1,\ldots,\mathcal{C}_t$ of the aggregate sparsity graph $\mathcal{E}$ are small, which is often true in many applications. In contrast, the first-order algorithms for generic SDPs developed in~\cite{wen2010alternating,o2016conic} require a projection onto the semidefinite cone $\mathbb{S}^n_+$ at each iteration, which becomes a bottleneck when $n \gg 1$. It is therefore clear that exploiting sparsity via chordal decomposition can bring significant computational savings in ADMM algorithms.

%%%%%%%%%%%%%%%%%%%%%%%%%%%%%%%%%%%%%%
\paragraph{Range-space decomposition:} The range-domain decomposition~\cref{eq:DecomposedDualSDP} of the dual-standard-form SDP~\cref{eq:SDP_dual} can  be solved using an ADMM algorithm very similar to that presented above. First, observe that~\cref{eq:DecomposedDualSDP} is equivalent to
\begin{subequations}
    \begin{align}
    \min_{y,V_k,Z_k} \quad
    &-\langle b,y \rangle + \chi_0\!\bigg( C - \sum_{i=1}^m A_i\,y_i - \sum_{k=1}^{t} E_{\mathcal{C}_k}^\tr  V_k E_{\mathcal{C}_k} \bigg) \nonumber\\
    &\qquad \qquad \qquad \qquad\qquad \qquad + \sum_{k=1}^{t} \chi_{\mathbb{S}^{|\mathcal{C}_k|}_+}(Z_k)  \nonumber\\
    \text{subject to} \quad &
    Z_k - V_k = 0, \quad k=1,\,\ldots,\,t. \nonumber%\label{E:DualDecomposedSDPvec}
    \end{align}
\end{subequations}

Grouping the variables as $\mathcal{X}:= \{y,V_1,\,\ldots,\,V_t\}$ and
$\mathcal{Y}:= \{Z_1,\,\ldots,\,Z_t\}$, this problem can be written in the general form~\cref{e:ADMM_base_prob} over the spaces $\mathbb{X}=\mathbb{R}^m \times \mathbb{S}^{\abs{\mathcal{C}_1}} \times \cdots \times \mathbb{S}^{\abs{\mathcal{C}_t}}$ and $\mathbb{Y} = \mathbb{Z} = \mathbb{S}^{\abs{\mathcal{C}_1}} \times \cdots \times \mathbb{S}^{\abs{\mathcal{C}_t}}$
Given a penalty parameter $\rho>0$ and a dual variable $\mathcal{Z}:= \{\Lambda_1,\,\ldots,\,\Lambda_t\}$, where each $\Lambda_k$ acts as a Lagrange multiplier for the corresponding constraint $Z_k - V_k=0$, one can easily verify that the ADMM step~\cref{eq:ADMM_S1} reduces to solving the equality-constrained quadratic program
\begin{subequations}
\begin{align*}
    \min_{y,V_1,\ldots,V_t} \quad &-\langle b,y\rangle + \frac{\rho}{2}\sum_{k=0}^{t} \left\| Z_k^{(q)} - V_k + \Lambda_k^{(q)} \right\|_F^2 \nonumber \\
    \text{subject to} \quad &  C - \sum_{i=1}^m A_i\,y_i - \sum_{k=1}^{t} E_{\mathcal{C}_k}^\tr  V_k E_{\mathcal{C}_k} = 0. 
    % \label{eq:MinXblockDual}
\end{align*}
Step~\cref{eq:ADMM_S2}, instead, reduces to $t$ independent positive semidefinite projections of the form
\begin{equation*}
   Z_k^{(q+1)} = \argmin_{Z_k \in \mathbb{S}_+^{|\mathcal{C}_k|}} \left\| Z_k - V_k^{(q+1)}  + \Lambda_k^{(q)} \right\|_F^2.
\end{equation*}
Finally, the dual variables $\Lambda_1,\ldots,\Lambda_t$ are updated through step~\cref{eq:ADMM_S3} as
\begin{equation*}
    % \label{eq:LambdakUpdateDual}
    \Lambda_k^{(q+1)} =
    \Lambda_k^{(q)} + Z_k^{(q+1)} - V_k^{(q+1)}.
\end{equation*}
\end{subequations}

Again, these iterations admit inexpensive closed-loop expressions. Moreover, it is not difficult to see that the ADMM iterations for the range-space decomposition~\cref{eq:DecomposedDualSDP} and for the domain-space decomposition~\cref{eq:DecomposedPrimalSDP} have similar leading-order complexity. In fact, \citet[Section 4.3]{ZFPGW2020chordal} showed that the ADMM algorithms for the primal and dual decomposed SDPs are scaled versions of each other. This shows that the duality picture of \cref{F:Duality,Fig:SparseCone_Duality} is reflected also at the algorithmic level. 

\begin{remark} 
For all fixed penalty $\rho>0$, the primal and dual ADMM algorithms outlined above converge to a solution of~\cref{eq:DecomposedPrimalSDP} and~\cref{eq:DecomposedDualSDP}, respectively, provided that strict primal-dual feasibility conditions are satisfied~\cite[Section 3.2]{boyd2011distributed}. An efficient ADMM algorithm that can handle primal or dual infeasible problems was developed by \citet[Section 5]{ZFPGW2020chordal}, who considered the homogeneous self-dual embedding \citep{ye1994nl,o2016conic} of the domain-space decomposition~\cref{eq:DecomposedPrimalSDP} and the range-space decomposition~\cref{eq:DecomposedDualSDP}. \markendexample
\end{remark}

\begin{remark} \label{remark:slack_varaibels_iterations}
As anticipated in~\Cref{remark:conversion_first_order_method}, considering the variables $X_k$ and the constraints $X_k = E_{\mathcal{C}_k}XE_{\mathcal{C}_k}^\tr$ without eliminating any redundant variables is essential to obtain efficient ADMM iterations. This is because the conic constraints separate completely from the affine ones in~\cref{eq:DecomposedPrimalSDP} when applying the splitting strategy of ADMM, making it easy to update each $X_k$ via simple projections onto positive semidefinite cones. Similarly, the redundant variables $V_k$ and the constraints $Z_k = V_k $ in~\cref{eq:DecomposedDualSDP} are essential to decouple the conic constraints from the affine ones, which enables one to handle positive semidefinite constraints via simple projections. \markendexample
\end{remark}

%%%%%%%%%%%%%%%%%%%%%%%%%%%%%%%%%%%%%%%%%%%%%%%%%%%%%%%%%%%%%%%%%%%%%%%%%%%%%%%%%
\subsection{Interior-point algorithms} \label{subsection:interior-point-methods}
\noindent
Interior-point algorithms for convex optimization problems with equality and inequality constraints employ Newton's method to solve a sequence of modified equality-constrained problems, obtained by replacing any inequality constraints with barrier functions in the objective~\citep{nesterov2003introductory,ye2011interior}. These barrier functions approximate the characteristic function of the set defined by the original inequality constraints and ensure that the optimal solution of each modified problem is strictly feasible for the original problem, meaning that it is an interior point of the original feasible set.

Since Newton's method relies on second-order (Hessian) information, interior-point algorithms do not share the slow convergence of first-order methods. Instead, they 
%Unlike the slow convergence of first-order methods, interior-point methods 
converge to an $\epsilon$-approximate solution using at most $\mathcal{O}(\log(1/\epsilon))$ Newton iterations \citep{nesterov2003introductory,ye2011interior}. In practice, convergence often occurs within tens of iterations. Therefore, interior-point methods are typically preferred when solving~\cref{eq:SDP_primal}-\cref{eq:SDP_dual} to high accuracy. The general-purpose SDP solvers \software{SeDuMi} \citep{sturm1999using}, \software{SDPT3} \citep{tutuncu2003solving}, \software{SDPA} \citep{yamashita2012latest}, and \software{MOSEK} \citep{mosek2015mosek} are all based on primal-dual interior-point methods, and they can very reliably solve small and medium-sized SDPs (e.g., when $n$ is less than a few hundreds and $m$ is less than a few thousands in~\cref{eq:SDP_primal}-\cref{eq:SDP_dual}) on regular computers. However, they become impractical for large SDPs because the CPU time and memory requirements for each interior-point iteration increase as $\mathcal{O}(n^3m + n^2m^2 + m^3)$ and $\mathcal{O}(n^2 + m^2)$, respectively \citep[Section 4.3.3]{nesterov2003introductory}. 
%$\mathcal{O}(n^6)$ time and $\mathcal{O}(n^4)$ memory when $m=\mathcal{O}(n^2)$; see, e.g.,~\cite[Section 4.3.3]{nesterov2003introductory}).

Chordal graph techniques can be exploited to improve the efficiency of interior-point methods when solving large-scale SDPs with chordal aggregate sparsity \citep{fukuda2001exploiting,andersen2011chordal,de2010exploiting}. This section reviews two general approaches for doing so. The first one, similar to the conversion methods in~\cref{subsection:conversion_first_order}, reformulates problems~\cref{eq:NonSymmetric_primal,eq:NonSymmetric_dual}  as SDPs with small positive semidefinite cones, which are often easier to solve with general-purpose interior-point solvers \citep{fukuda2001exploiting,nakata2003exploiting,kim2011exploiting, zhang2020sparse}. The second approach, instead, directly solves~\cref{eq:NonSymmetric_primal}-\cref{eq:NonSymmetric_dual} using an interior-point method for nonsymmetric conic optimization \citep{andersen2010implementation,nesterov2012towards,skajaa2015homogeneous,coey2020towards}. For other ways to exploit chordal sparsity in the computation of interior-point search directions, we refer the reader to the works by \citep{benson2000solving}, \cite{pakazad2017distributed} and \citet[Section 5]{fukuda2001exploiting}.

\subsubsection{Conversion methods} \label{subsection:conversion_IPM}
\noindent
Starting from the domain-space decomposed SDP~\cref{eq:DecomposedPrimalSDP}, \cite{fukuda2001exploiting} and \cite{kim2011exploiting} suggested to eliminate the global matrix $X$ and rewrite the SDP~\cref{eq:NonSymmetric_primal} only in terms of variables $X_k \in \mathbb{S}^{|\mathcal{C}_k|}_+, k = 1, \ldots t$. To rewrite the cost function and the first set of equality constraints, one must choose matrices $C_k$ and $A_{ik}$ that satisfy 
$$
    \sum_{k=1}^t \langle  C_k, X_k \rangle = \langle C, X \rangle
$$
and
$$
     \sum_{k=1}^t \langle  A_{ik}, X_k \rangle = \langle A_i, X \rangle, \quad i = 1, \ldots, m. 
$$
These affine relations do not usually determine $C_k$ and $A_{ik}$ uniquely, and some choices may be more convenient than others from the point of view of computations (\citealp[Section 3.1]{sun2014decomposition}, \citealp[Section 6]{zhang2020sparse}). The second set of constraints in~\cref{eq:DecomposedPrimalSDP}, instead, can be enforced via consistency constraints on the entries of $X_1,\ldots,X_t$ that correspond to the same elements of $X$. Such consistency constraints can be formulated as 
\begin{equation}\label{eq:verlappingEquality}
\begin{aligned}
            &E_{\mathcal{C}_j \cap \mathcal{C}_k}
            \left(
            		E_{\mathcal{C}_k}^\tr  X_k E_{\mathcal{C}_k}
            	  - E_{\mathcal{C}_j}^\tr  X_j E_{\mathcal{C}_j}
            	  \right)
            	 E_{\mathcal{C}_j \cap \mathcal{C}_k}^\tr  =0 \\
            &\qquad \qquad \qquad \qquad  \forall j,k :\; \mathcal{C}_j \cap \mathcal{C}_k \neq \emptyset.
\end{aligned}
\end{equation}
The primal SDP~\cref{eq:NonSymmetric_primal} can therefore be rewritten as 
\begin{equation} \label{eq:conversion_primal_IPM}
    \begin{aligned}
        \min_{X_1, \ldots, X_t} & \quad \sum_{k=1}^t \langle  C_k, X_k \rangle  \\
        \text{subject to} &\quad \sum_{k=1}^t \langle A_{ik}, X_k \rangle = b_i, \; i = 1, \ldots, m, \\
        & \quad\cref{eq:verlappingEquality}, X_k \in \mathbb{S}^{|\mathcal{C}_k|}_+,\;\; k = 1, \ldots t.
    \end{aligned}
\end{equation}

This conversion process, first proposed in~\cite{fukuda2001exploiting}, is known as the \emph{domain-space} decomposition \citep{kim2011exploiting}. The reformulated problem~\cref{eq:conversion_primal_IPM} has more variables and constraints than the original SDP~\cref{eq:SDP_primal}, but the large matrix constraint $X \in \mathbb{S}^n_+$ is replaced by $t$ smaller ones, $X_k \in \mathbb{S}^{|\mathcal{C}_k|}_+$ for $k = 1, \ldots, t$. In certain cases, the decomposed problem~\cref{eq:conversion_primal_IPM} is easier to solve than the original SDP~\cref{eq:SDP_primal} using general-purpose interior-point solvers; see~\cite{nakata2003exploiting} and \cite{fujisawa2009user} for numerical examples. 
Three other variants of this conversion method, including \emph{range-space} decompositions, have been studied by~\cite{kim2011exploiting}. 

The main drawback of these conversion methods is that, sometimes, the additional consistency constraints~\cref{eq:verlappingEquality} significantly increase the size of the Schur complement system that needs to be solved at each interior-point iteration. This can offset the benefits of the clique-based matrix decomposition. As shown recently by \cite{zhang2020sparse}, this issue can be mitigated using a dualization technique \citep{lofberg2009dualize}.

\begin{remark}[Removing redundant constraints]
Since the maximal cliques in a chordal graph satisfy the running intersection property \citep{blair1993introduction,fukuda2001exploiting} (see also \cref{app:maximal_cliques}), it is in fact sufficient to enforce the consistency between pairs $\mathcal{C}_j, \mathcal{C}_k$ that correspond to the parent-child pairs in a clique tree. Redundant constraints in~\cref{eq:verlappingEquality} can therefore be removed using the running intersection property. Interested readers are referred to~\cite{kim2011exploiting} and \cite{vandenberghe2015chordal} for details. \markendexample
\end{remark}

\begin{remark}[Dropping or fixing consistency constraints]
In some applications, the SDP~\cref{eq:SDP_primal} comes from a semidefinite relaxation of a nonconvex optimization problem. Dropping some consistency constraints in~\cref{eq:verlappingEquality} leads to a valid weaker relaxation with a lower computational complexity. This idea was successfully applied to semidefinite relaxations for optimal power flow problems \citep{andersen2014reduced} and neural network verification \citep{BPLZ2021neural}.  
Other times, one can enforce some of the consistency conditions \emph{a priori} and look for feasible (but suboptimal) points for an SDP at a low computational cost. This idea was used  in~\cite{ZMP2018Scalable} to develop a scalable approach for solving distributed control problems.  \markendexample
\end{remark}

\subsubsection{Nonsymmetric interior-point algorithms} \label{subsection:nonsymmetric_algorithms}
\noindent
Chordal graph techniques can also be exploited to speed up interior-point methods for the nonsymmetric pair of sparse SDPs~\cref{eq:NonSymmetric_primal}-\cref{eq:NonSymmetric_dual} without appealing to the matrix decomposition and conversion frameworks described above. Since the cones $\mathbb{S}^n_+(\mathcal{E},?)$ and $\mathbb{S}^n_+(\mathcal{E},0)$ are not self-dual, such sparsity-exploiting methods cannot enjoy a complete primal-dual symmetry \citep{andersen2010implementation}. Instead, one must resort to purely primal, purely dual, or nonsymmetric primal-dual path-following methods \citep{andersen2010implementation,nesterov2012towards,skajaa2015homogeneous,burer2003semidefinite,coey2020towards}.

To construct nonsymmetric interior-point methods, \cite{dahl2008covariance} and \cite{andersen2010implementation} introduced barrier functions $\phi: \mathbb{S}^n(\mathcal{E},0) \rightarrow \mathbb{R}$ and $\phi_*: \mathbb{S}^n(\mathcal{E},0) \rightarrow \mathbb{R}$ for the cones $\mathbb{S}^n_+(\mathcal{E},0)$ and $\mathbb{S}_+^n(\mathcal{E},?)$, defined as
\begin{subequations} \label{eq:barrier_functions}
    \begin{gather} 
    \phi(Z) = \begin{cases}
    - \log \det Z & Z \in {\rm int}(\mathbb{S}^n_+(\mathcal{E},0)),\\
    +\infty &\text{otherwise},
    \end{cases}  \label{eq:dual_barrier}\\
    \intertext{and}
    \phi_*(X) = \sup_{Z \in \mathbb{S}^n(\mathcal{E},0)} (-\langle X, Z \rangle - \phi(Z)). \label{eq:primal_barrier}
    \end{gather}
\end{subequations}
Note that $\phi$ (resp. $\phi_*$) is finite only on the interior of $\mathbb{S}^n_+(\mathcal{E},0)$ (resp. $\mathbb{S}_+^n(\mathcal{E},?)$) and tends to $+\infty$ as $Z$ (resp. $X$) approaches the boundary of this cone. Observe also that $\phi_*$ is simply the Legendre transform of $\phi$ evaluated at $-X$.

Thanks to the properties of the barrier functions, a minimizing sequence $\{X^\mu\}_{\mu>0}$ for~\cref{eq:NonSymmetric_primal} can be computed by solving the regularized primal problem
\begin{equation}\label{eq:SDP_primal_barrier}
    \begin{aligned}
            \min_{X} \quad & \langle C, X \rangle + \mu \phi_*(X) \\
            \text{subject to} \quad & \langle A_i, X\rangle = b_i,
            \quad i = 1,\,\ldots,\,m,
    \end{aligned}
\end{equation}
and letting $\mu \to 0$. Similarly, a minimizing sequence $\{y^\mu,Z^\mu\}_{\mu>0}$ for~\cref{eq:NonSymmetric_dual} is found upon solving the regularized dual problem
\begin{equation}\label{eq:SDP_dual_barrier}
  \begin{aligned}
            \max_{y, Z} \quad & \langle b,y\rangle - \mu \phi(Z) \\
            \text{subject to} \quad  & Z + \sum_{i=1}^m A_i\,y_i = C
        \end{aligned}
\end{equation}
for $\mu \to 0$. Solutions of the regularized problems for fixed finite $\mu$ are usually found using Newton's method, leading to so-called \emph{primal scaling} and \emph{dual scaling} interior-point methods. Other methods can also be used; for instance, \cite{jiang2021bregman} recently suggested solving~\cref{eq:SDP_primal_barrier} with a Bregman first-order method, where the complexity of evaluating the Bregman proximal operator can be reduced using a sparse Cholesky factorization. 
%
%\red{[Summarize the advantages of this in a sentence?]}

When Newton's method is applied to~\cref{eq:SDP_primal_barrier}, the KKT optimality conditions are
\begin{subequations} \label{eq:KKT_primal_barrier}
    \begin{align}
        \langle A_i, X^\mu \rangle &= b_i, \quad i = 1, \ldots, m, \\
        \sum_{i=1}^m y_i A_i + Z &= C, \\
        \mu \nabla \phi_*(X^\mu) + Z&= 0, \label{eq:KKT_primal_barrier_s3}
    \end{align}
\end{subequations}
where $y \in \mathbb{R}^m$ is a Lagrange multiplier for the equality constraint in~\cref{eq:SDP_primal_barrier} and $Z$ is an auxiliary variable arising from the definition of $\phi_*$ via the Legendre transform. Solutions $X^\mu \in \mathbb{S}^n_+(\mathcal{E},?)$ as $\mu$ is varied define the so-called \textit{central path} for~\cref{eq:NonSymmetric_primal}. Similarly, the KKT optimality conditions for~\cref{eq:SDP_dual_barrier} are
\begin{subequations} \label{eq:KKT_dual_barrier}
    \begin{align}
        \langle A_i, X\rangle &= b_i, \quad i = 1, \ldots, m, \\
        \sum_{i=1}^m y^\mu_i A_i + Z^\mu &= C, \\
        \mu \nabla \phi(Z^\mu) +X&= 0, \label{eq:KKT_dual_barrier_s3}
    \end{align}
\end{subequations}
where $X$ is a Lagrange multiplier for the equality constraint in~\cref{eq:SDP_dual_barrier}. Solutions $\{y^\mu, Z^\mu\} \in \mathbb{R}^m \times \mathbb{S}^n_+(\mathcal{E},0)$ as $\mu$ is varied define the central path for~\cref{eq:NonSymmetric_dual}. It is possible to show that~\cref{eq:KKT_primal_barrier} and~\cref{eq:KKT_dual_barrier} are equivalent~\cite[Chapter 3]{andersen2011chordal}, so the set of points $\{X^\mu, y^\mu, Z^\mu\}_{\mu>0}$ in $\mathbb{S}^n_+(\mathcal{E},?) \times \mathbb{R}^m \times \mathbb{S}^n_+(\mathcal{E},0)$ define a primal-dual central path. 

The rest of this section briefly outlines how the chordality of the sparsity pattern $\mathcal{E}$ can be exploited in the context of dual-scaling interior point methods. Similar ideas can be used to formulate primal-scaling methods, and we refer interested readers to the work by~\citet[Section 4.2]{andersen2010implementation} for details.

\paragraph{Dual-scaling interior-point methods} Search directions in a dual-scaling interior-point method are obtained by linearizing~\cref{eq:KKT_dual_barrier} around the current interior iterate $X \in {\rm int}(\mathbb{S}^n_+(\mathcal{E},?))$, $y\in \mathbb{R}^m$ and $Z \in {\rm int}(\mathbb{S}^n_+(\mathcal{E},0))$. Replacing $X$, $y$ and $Z$ with $X + \Delta X$, $y + \Delta y$, $Z + \Delta Z$ in~\cref{eq:KKT_dual_barrier}, linearizing~\cref{eq:KKT_dual_barrier_s3}, and eliminating $\Delta Z$ yields the Newton equations %\red{[GF, YZ: should be $1/\mu$?]}
\begin{equation} \label{eq:dual_newton}
    \begin{aligned}
        \langle A_i, \Delta X \rangle &= r_i, \, i = 1, \ldots, m, \\
        \sum_{i=1}^m \Delta y_i A_i - \frac{1}{\mu} \nabla^2 \phi(Z)^{-1}[\Delta X] &= R,
    \end{aligned}
\end{equation}  
where $\nabla^2 \phi(Z)^{-1}$ is the inverse Hessian of $\phi$ at $Z$, $r_i = b_i - \langle A_i, X \rangle$ and
%
%\begin{align*}
    %, &&
 $R = C - \sum_{i = 1}^m y_i A_i - 2Z + \frac{1}{\mu} \nabla^2 \phi(Z)^{-1}[X].$
%\end{align*}
%
Further elimination of $\Delta X$ leads to the Schur complement equation
\begin{equation} \label{eq:SchurComplement}
    H \Delta y = g,
\end{equation}
where $g \in \mathbb{R}^m$ is a vector and $H$ is an $m \times m$ positive definite matrix, both depending only on the current (known) iterates $X$, $y$ and $Z$. Explicit expression for these quantities are given by~\citet[Section 4.3]{andersen2010implementation}. 
%\red{[GF: Is it obvious that these Newton/Schur equations guarantee the symmetry of the seach directions $\Delta X$ and $\Delta Z$? This is an issue in standard interior-point methods, and one must ``symmetrize'' the Newton system to preserve symmetry. I am not familiar enough with interior-point methods to be able to comment on this. \textbf{UPDATE:} I checked and one must be careful. I think that everything is fine because the hessians $\nabla^2 \phi$ and $\nabla^2\phi_*$ should map symmetric matrices onto symmetric matrices. This is not immediately obvious and must be checked, for instance using the explicit expression for $\nabla^2 \phi$ given below...]}

Finding the dual-scaling search direction  $\Delta X, \Delta y, \Delta Z$ requires solving the Newton equation~\cref{eq:dual_newton} or the Schur complement equation~\cref{eq:SchurComplement}. To do this using a direct method, one must first calculate the Hessian and inverse Hessian of the barrier function $\phi(X)$ in~\cref{eq:dual_barrier}, and then form and factorize the matrix $H$. This is the most computationally expensive part of any interior-point method. It is in this computation that one can exploit the chordality of the sparsity pattern $\mathcal{E}$ \citep{andersen2010implementation}.

\paragraph{Fast calculations involving the barrier functions} 
The value, gradient, Hessian, and inverse Hessian of the %the primal barrier $\phi_*(X)$ and 
dual barrier $\phi(Z)$ in~\cref{eq:dual_barrier} can be computed efficiently if the sparsity pattern $\mathcal{E}$ of $Z$ is chordal. Similar fast algorithms exist for the primal barrier $\phi_*(X)$ in~\cref{eq:primal_barrier}, but we do not review them here and refer interested readers to~\citet[Section 3.2]{andersen2010implementation} for details.  

The key ingredient of these efficient algorithms is a sparse Cholesky factorization with zero fill-in \citep{rose1970triangulated,blair1993introduction,vandenberghe2015chordal}: as reviewed in \cref{appendix:zero_fillin}, for any positive definite matrix $Z$ in ${\rm int}(\mathbb{S}^n_+(\mathcal{E},0))$ with chordal sparsity there exists a permutation matrix $P$ and a lower triangular matrix $L$ such that
\begin{equation} \label{eq:sparseCholesky}
    P Z P^\tr  = LL^\tr, \qquad P^\tr (L + L^\tr)P \in \mathbb{S}^n(\mathcal{E},0).
\end{equation}
%As highlighted in~\cref{subsection:SparseMatrix}, this is one central result for proving the chordal decomposition in~\cref{T:ChordalDecompositionTheorem,T:ChordalCompletionTheorem}. 
%This sparse factorization is the basis of the efficient algorithms below. 
%For sparse chordal patterns, 
This factorization can be computed efficiently by following a recursion on a clique tree~\cite[Chapter 9.3]{vandenberghe2015chordal}. 

Now, to evaluate $\phi(Z)$ it suffices to substitute $Z = P^\tr L L^\tr P$ into~\cref{eq:dual_barrier} and observe that
$$
    \phi(Z) = -2 \sum_{i=1}^n \log L_{ii}
$$
because determinants distribute over products and permutation matrices have unit determinant. Thus, $\phi(Z)$ can be evaluated efficiently once the Cholesky factorization~\cref{eq:sparseCholesky} has been computed.

The gradient of $\phi(Z)$, instead, is given by the following negative projected inverse
$$
    \nabla \phi(Z) = - \mathbb{P}_{\mathbb{S}^n(\mathcal{E},0)}(Z^{-1}). 
$$
Despite the fact that $Z^{-1}$ is in general dense, the projection onto $\mathbb{S}^n(\mathcal{E},0)$ can be computed from its sparse Cholesky factorization~\cref{eq:sparseCholesky} without computing any other entries of $Z^{-1}$~\cite[Chapter 9.5]{vandenberghe2015chordal}. 

The Hessian of $\phi$ at $Z$ applied to a matrix $Y \in \mathbb{S}^n(\mathcal{E},0)$ is computed as 
$$
    \nabla^2 \phi(Z)[Y] = \frac{d}{dt} \nabla \phi(Z + tY) \mid_{t = 0} = \mathbb{P}_{\mathbb{S}^n(\mathcal{E},0)}(Z^{-1}YZ^{-1}).
$$
Again, this quantity can be evaluated knowing only the sparse Cholesky factorization of $Z$ and its projected inverse $\mathbb{P}_{\mathbb{S}^n(\mathcal{E},0)}(Z^{-1})$, without explicitly computing the inverse $Z^{-1}$ or the matrix product $Z^{-1}YZ^{-1}$ \citep{andersen2010implementation,andersen2013logarithmic}. 

Finally, thanks to the chordal structure, solving the linear equation $\nabla^2 \phi(Z)[U] = Y$ for $U$ in order to evaluate the inverse Hessian $\nabla^2 \phi(Z)^{-1}[Y]$ has the same cost as the
evaluating the Hessian $\nabla^2 \phi(Z)[Y]$; see \citet[Section 3.2]{andersen2010implementation} and \cite{andersen2013logarithmic}.

\subsection{Algorithm implementations} \label{subsection:implementations_SDP}
\noindent
We conclude this section by providing a list of numerical packages that implement some of the approaches reviewed above. This list is not exhaustive, and the goal here is to give the interested reader a starting point for numerical experiments. First-order solvers based on augmented Lagrangian methods and ADMM for generic SDPs include \software{SDPNAL}/\software{SDPNAL+} \citep{sun2020sdpnal+, zhao2010newton} and \software{SCS} \citep{scs}. \software{CDCS} \citep{CDCS} and \software{COSMO} \citep{garstka2019cosmo} are two open-source first-order solvers that exploit chordal sparsity in SDPs. The MATLAB package \software{CDCS} implements the algorithms described in~\Cref{subsection:ADMM} and has interfaces with the optimization toolboxes \software{YALMIP} \citep{lofberg2004yalmip} and \software{SOSTOOLS} \citep{prajna2002introducing}. The Julia package \software{COSMO} solves SDPs with quadratic objective functions.

The conversion methods in~\Cref{subsection:conversion_IPM} are implemented in \software{SparseCoLO} \citep{fujisawa2009user} and \software{CHOMPACK} \citep{andersen2015chompack}. We note that \software{CHOMPACK} also provides useful implementation of many other chordal matrix computations, including maximum determinant positive definite completion and minimum rank positive semidefinite completion. Another MATLAB package \software{Dual-CTC} \citep{dual_ctc} implements a dualized clique tree conversion \citep{zhang2020sparse}. The reformulated SDPs after conversion can be solved using general-purpose interior-point solvers, such as \software{SeDuMi} \citep{sturm1999using}, \software{SDPT3} \citep{tutuncu2003solving}, \software{SDPA} \citep{yamashita2012latest}, and \software{MOSEK} \citep{mosek2015mosek}. \software{SMCP} \citep{SMCP} is a nonsymmetric interior-point solver that provides a Python implementation of the algorithms in~\Cref{subsection:nonsymmetric_algorithms}. Finally, \software{SDPA-C} \citep{fujisawa2004sdpa} is a primal-dual interior-point solver that exploits chordal sparsity using the maximum-determinant positive definite completion. 

%%%%%%%%%%%%%%%%%%%%%%%%%%%%%%%%%%%%%%%%%%%%%%%%%%%%
\section{Sparse polynomial optimization}
\label{section:polynomial_optimization}
\noindent
We have seen in \cref{section:sparse-SDPs} that the chordal decomposition of large semidefinite matrices allows for significant efficiency gains in the solution of sparse SDPs. The same ideas can often be leveraged to replace SDP relaxations of intractable optimization problems, which generally have no inherent sparsity or other computationally advantageous structure, with SDPs that do.

This section describes how sparsity (primarily chordal, but also nonchordal) can be exploited in the context of sum-of-squares (SOS) relaxation techniques for polynomial optimization. As mentioned in the introduction, SOS methods are at the heart of many recent tractable frameworks for the analysis and optimal control of nonlinear systems with polynomial dynamics; see \cite{Lasserre2008,valmorbida2017region,henrion2014convex,Lasagna2016sos,Jones2019reachable_sets,Majumdar2014occupation_measures,Han2018control_om,fantuzzi2016bounds,Fantuzzi2020siads,Goluskin2020attractors,Korda2021invariant_measures,Miller2021peak,Ahmadi2018robust,papachristodoulou2005tutorial,prajna2004nonlinear} to name but a few contributions. %\red{[Add more to balance out the French group?]}

Our goal is not to offer an exhaustive review of all sparsity-exploiting methods that have been proposed in this field, but rather to introduce the key ideas underpinning most of these methods from a general perspective, in the hope that this can guide further developments. For this reason, we concentrate mainly on two basic problems. The first, discussed in \cref{s:pop-sparse-global}, is to prove that an $n$-variate polynomial of even degree $2d$ is a sum of squares and, therefore, globally nonnegative. In this case, we seek to exploit the structure of polynomials that depend only a small subset of all possible degree-$2d$ monomials---a property often referred to as \textit{term sparsity}. The second problem, discussed in~\cref{s:pop-sparse-matrix}, is to check whether a sparse and symmetric $n$-variate polynomial matrix $P(x)$ is SOS, and therefore positive semidefinite for all $x \in \mathbb{R}^n$. In this case, our goal is to leverage the \textit{structural sparsity} of $P$, meaning that many of its entries are zero.

Although we focus only on global nonnegativity, all of the sparsity-exploiting techniques discussed in this section can be extended to prove polynomial (matrix) nonnegativity locally on basic semialgebraic sets. Such extensions, which have been studied extensively in order to build hierarchies of sparse SDP relaxations for polynomial optimization problems~\citep{waki2006sums,waki2008algorithm,lasserre2006convergent,Wang2020tssos,Wang2020chordal-tssos,Wang2020cs-tssos,zheng2020sum}, require some careful technical adjustments, but the underlying strategy is the same as for the global nonnegativity setting. We outline some of these adjustments in \cref{s:pop-sparse-constrained,ss:local-matrix-sos}, and refer readers to the excellent literature on this topic for full details.

%%%%%%%%%%%%%%%%%%%%%%%%%%%%%%%%%
\subsection{Background}\label{s:pop-preliminaries}
\noindent
Let $\mathbb{R}[x]_{n,d}$ be the $\binom{n+d}{d}$-dimensional space of polynomials with independent variables $x=(x_1,\ldots,x_n)$ and degree no larger than $d$. The $n$-variate monomial with exponent $\beta = (\beta_1,\ldots,\beta_n) \in \mathbb{N}^n$ and degree $\abs{\beta} = \beta_1 + \cdots +\beta_n$ is denoted by $x^\beta= x_1^{\beta_1}x_2^{\beta_2}\cdots x_n^{\beta_n}$. Given a finite set of exponents $\mathbb{B} \subset \mathbb{N}^n$, we write $x^\mathbb{B} = (x^\beta)_{\beta \in \mathbb{B}}$ for the (column) vector of monomials with exponents in $\mathbb{B}$. The cardinality of $\mathbb{B}$ is denoted by $\abs{\mathbb{B}}$. We also define
\begin{subequations}
	\begin{align}
	\mathbb{B}+ \mathbb{B} &:= \{\beta+\gamma: \; \beta, \gamma \in \mathbb{B}\},\\
	2\mathbb{B} &:= \{ 2\beta : \; \beta \in \mathbb{B} \}.
	\end{align}
\end{subequations}

If $\mathbb{N}^n_d = \{\beta\in\mathbb{N}^n: \abs{\beta} \leq d\}$ is the set of all $n$-variate exponents of degree $d$ or less, the vector $x^{\mathbb{N}^n_d}$ is a basis for $\mathbb{R}[x]_{n,d}$ and any polynomial $f \in \mathbb{R}[x]_{n,d}$ can be written as $f(x) = \sum_{\beta \in \mathbb{N}^n_d} f_\beta x^\beta$ for some coefficients $f_\beta \in \mathbb{R}$. The set of exponents with nonzero coefficient,
\begin{equation}
\supp(f) = \{\beta \in \mathbb{N}^n_d: f_\beta \neq 0\},
\end{equation}
is called the \textit{support} of $f$.
Its convex hull is called the \textit{Newton polytope} of $f$ and is denoted by $\New(f)$.

\subsubsection{SOS polyonomials and SDPs}
\noindent
A polynomial $f \in \mathbb{R}[x]_{n,2d}$ of even degree $2d$ is SOS if there exist degree-$d$ polynomials $f_1,\ldots,f_k \in \mathbb{R}[x]_{n,d}$ such that
\begin{equation}\label{e:sos-definition}
f = f_1^2 + \cdots + f_k^2.
\end{equation}
The set of $n$-variate degree-$2d$ SOS polynomials, denoted by $\Sigma_{n,2d}$, is a proper cone in $\mathbb{R}[x]_{n,2d}$~\cite[Theorem 3.26]{blekherman2012semidefinite}. Given an exponent set $\mathbb{A} \subseteq \mathbb{N}^n_{2d}$, we define the subcone of SOS polynomials supported on $\mathbb{A}$ as
\begin{equation}
\Sigma[\mathbb{A}] := \{f \in \Sigma_{n,2d}:\; \supp(f) \subseteq \mathbb{A}  \}.
\end{equation}

It is well known (see, e.g., \citealp{parrilo2003semidefinite,parrilo2013semidefinite}) that a polynomial $f \in \mathbb{R}[x]_{n,2d}$ is SOS if and only if there exist a set of exponents $\mathbb{B} \subseteq \mathbb{N}^n_d$ and a positive semidefinite matrix $Q \in \mathbb{S}^{\abs{\mathbb{B}}}_+$ such that
\begin{equation}\label{e:sos-gram}
f(x) = (x^{\mathbb{B}} )^\tr \, Q \, x^{\mathbb{B}}.
\end{equation}
In particular, if $f$ is SOS, this so-called \textit{Gram matrix representation}~\cref{e:sos-gram} is guaranteed to exist with \citep{reznick1978extremal}
\begin{equation}\label{e:newton-basis}
\mathbb{B} = \frac12 \New(f) \cap \mathbb{N}^n_d.
\end{equation}
The exponent set obtained with this Newton polytope reduction can be simplified further using more general \textit{facial reduction} techniques \citep{lofberg2009pre,permenter2014basis,permenter2014partial,waki2010facial}. These techniques analyze the support of $f$ in order to remove redundant elements from $\mathbb{B}$, and construct a smaller exponent set for which~\cref{e:sos-gram} is guaranteed to hold as long as $f$ is SOS.

It is clear that SOS polynomials are nonnegative globally. The converse is true only for univariate polynomials ($n=1$, $d$ arbitrary), quadratic polynomials ($d=1$, $n$ arbitrary), and bivariate quartics ($n=2$, $d=2$) \citep{Hilbert1888}. In general, therefore, being SOS is only a sufficient condition for global nonnegativity, and there are well-known examples of nonnegative polynomials that are not SOS, such the Motzkin polynomial \citep{Motzkin1967}. However, while verifying polynomial nonnegativity is an NP-hard problem \citep{Murty1987}, checking whether a polynomial $f$ is SOS can be done in polynomial time by solving an SDP. Specifically, for each exponent $\alpha \in \mathbb{B}+ \mathbb{B}$, let $A_\alpha \in \mathbb{S}^{\abs{\mathbb{B}}}$ be the symmetric binary matrix satisfying
\begin{equation}
[A_\alpha]_{\beta,\gamma} := \begin{cases}
1, &\beta+\gamma=\alpha,\\
0, &\text{otherwise},
\end{cases}
\end{equation}
and observe that
\begin{equation}
(x^{\mathbb{B}} )^\tr \, Q \, x^{\mathbb{B}} 
= \langle Q, x^{\mathbb{B}} (x^{\mathbb{B}} )^\tr\rangle 
= \sum_{\alpha \in \mathbb{B}+\mathbb{B}} \langle Q, A_\alpha \rangle x^\alpha.
\end{equation}
Then, condition~\cref{e:sos-gram} holds if and only if $\langle Q, A_\alpha\rangle = f_\alpha$ for all $\alpha \in \mathbb{B}+ \mathbb{B}$ and we conclude that
\begin{equation}\label{e:sos-sdp}
f \in \Sigma_{n,2d}
\;\;\iff\;\;
\begin{cases}
\exists Q \in \mathbb{S}^{\abs{\mathbb{B}}}_+ \text{ such that}\\
\langle Q, A_\alpha\rangle = f_\alpha \quad \forall \alpha \in \mathbb{B}+ \mathbb{B}.
\end{cases} 
\end{equation}
The condition on the right-hand side defines an SDP, so a positive semidefinite Gram matrix $Q$ certifying that $f$ is SOS can (in principle) be constructed in polynomial time.

%%%%%%%%%%%%%%%%%%%%%%%%%%%%%%%%%%%%%%%%%%%%%%%%
\subsubsection{SOS polynomial matrices and SDPs}\label{s:matrix-sos}
\noindent
Let $\mathbb{R}[x]_{n,d}^{r \times s}$ be the space of $r \times s$ matrices whose entries are $n$-variate polynomials of degree $d$. We say that a symmetric polynomial matrix $P \in \mathbb{R}[x]_{n,2d}^{r \times r}$ is positive semidefinite (resp. definite) globally if $P(x) \succeq 0$ (resp. $P(x) \succ 0$) for all $x \in \mathbb{R}^n$. We also say that $P$ is positive semidefinite locally on a set $\mathbb{K}$ if the same conditions hold for $x \in \mathbb{K}$, but not necessarily otherwise.

A symmetric polynomial matrix $P \in \mathbb{R}[x]_{n,2d}^{r \times r}$ is called SOS if there exists an integer $s$ and a polynomial matrix $M \in \mathbb{R}[x]_{n,d}^{s \times r}$ such that
\begin{equation}\label{e:matrix-sos-def}
P(x) = M(x)^\tr M(x).
\end{equation}
The set of $r \times r$ SOS polynomial matrices with entries in $\mathbb{R}[x]_{n,2d}$ will be denoted by $\Sigma^{r}_{n,2d}$. All SOS polynomial matrices are clearly positive semidefinite globally, and the converse is true in the univariate case ($n=1$); see~\cite{aylward2007explicit} for a recent proof.

It is well known (see, e.g., \citealp{kojima2003sums,parrilo2013semidefinite,gatermann2004symmetry}) that a symmetric polynomial matrix $P \in \mathbb{R}[x]_{n,2d}^{r \times r}$ is SOS if and only if it admits a Gram matrix representation in the form
\begin{equation}\label{e:matrix-sos-gram}
P(x) = (I_r \otimes x^\mathbb{B})^\tr \, Q \, (I_r \otimes x^\mathbb{B})
\end{equation}
for some exponent set $\mathbb{B} \subseteq \mathbb{N}^n_d$ and some positive semidefinite symmetric matrix $Q \in \mathbb{S}_+^{r\abs{\mathbb{B}}}$. One may always take $\mathbb{B} = \mathbb{N}^n_d$, and smaller exponent sets can be constructed with the same reduction techniques used for SOS polynomials. As in the scalar case ($r=1$), condition \cref{e:matrix-sos-gram} defines a set of affine constraints on $Q$, so verifying that a polynomial matrix is SOS amounts to solving an SDP.

%%%%%%%%%%%%%%%%%%%%%%%%%%%%%%%%%%%%%%%%%
\subsection{Sparse SOS decompositions}\label{s:pop-sparse-global}
\noindent
A major obstacle to constructing SOS certificates of global polynomial nonnegativity via semidefinite programming is that the matrix $Q$ is both dense and very large. If $f \in \mathbb{R}[x]_{n,2d}$ has dense support $\supp(f)=\mathbb{N}^n_{2d}$, then one must take $\mathbb{B}=\mathbb{N}^n_d$ 
%to be the full set of exponents of degree $d$ or less, 
and $Q$ is a $\binom{n+d}{d} \times \binom{n+d}{d}$ dense matrix. Often, however, the support of $f$ is small, i.e.,
$\abs{\supp(f)}$ is much smaller than $\binom{n + 2d}{2d}$.
% \begin{equation}
% %\max\left\{ \abs{\supp(f)}, \abs{\supp(g_1)},\ldots, \abs{\supp(g_m)} \right\} \ll \binom{n + 2\omega_0}{2\omega_0}.
% \abs{\supp(f)} \ll \binom{n + 2d}{2d}.
% \end{equation}
%
This property, called \textit{term sparsity} \citep{Wang2019sos-term-sparsity,Wang2020cs-tssos,Wang2020tssos,Wang2020chordal-tssos}, 
can be exploited in various ways to reduce the computational complexity of the SDP in~\cref{e:sos-sdp}.

The facial reduction techniques mentioned above, which replace the full exponent set $\mathbb{N}^n_d$ with a (sometimes significantly) smaller subset, are arguably the simplest way to exploit term sparsity. However, as the next example demonstrates, they are often not sufficient.

\begin{example}\label{example:facial-reduction}
	Fix $n=50$ and $d=2$. The support of %the quartic polynomial
	\begin{equation}\label{e:facial-reduction-quartic}
	f(x) = \sum_{i=2}^{49} (x_{i-1} + x_{i} + x_{i+1})^4
	\end{equation}
	contains only $485$ out of the $\binom{50+4}{4}=316\,251$ possible monomials, so $f$ is term sparse. However, it is not hard to check that the Newton polytope $\New(f)$ consists of all points $\xi \in \mathbb{R}^{50}_+$ with $\|\xi \|_1 = 4$, so the Newton-reduced exponent set $\mathbb{B}=\frac12\New(f) \cap \mathbb{N}^{50}_2 = \mathbb{N}^{50}_2 \setminus \mathbb{N}^{50}_1$ contains \textit{all} homogeneous exponents of degree $2$. Therefore, Newton polytope reduction removes only $\binom{50+1}{1} = 51$ of the possible $\binom{50+2}{2}=1326$ in the full set $\mathbb{N}^{50}_2$, and the SDP in~\cref{e:sos-sdp} still involves a $1275 \times 1275$ Gram matrix $Q$. \markendexample
\end{example}

Techniques to exploit term sparsity beyond what can be achieved with facial reduction methods alone are clearly desirable.
\Cref{s:sparse-sos-general-strategy} describes a general strategy to search for \textit{sparse} SOS decompositions, which is based on the same matrix decomposition approach used to tackle large-scale sparse SDPs in \cref{section:sparse-SDPs}. 
%, that often brings considerable reduction in computational complexity. 
\Cref{s:sos-csp-unconstrained,s:tssos} show that different types of sparse SOS decompositions proposed in the literature are particular cases of this general approach. 
\Cref{s:pop-sparse-constrained} outlines how these methods can be extended to prove polynomial nonnegativity on basic semialgebraic sets, rather than globally. 
Throughout, $\mathbb{B}$ will denote a fixed set of candidate exponents for the SOS decomposition of a polynomial $f$, generated from $\mathbb{N}^n_d$ using facial reduction or any other exponent selection technique.

\subsubsection{General approach}
\label{s:sparse-sos-general-strategy}
\noindent
Let $\mathbb{A}$ be a small subset of $\mathbb{N}^n_{2d}$ and $f$ be a term-sparse polynomial supported on $\mathbb{A}$. To reduce the cost of testing if $f$ is SOS, a natural idea is to check whether $f$ belongs to a subset of the sparse SOS cone $\Sigma[\mathbb{A}]$ that admits a semidefinite representation with low computational complexity. Such a subset can be constructed using a simple strategy: \textit{prescribe a sparsity graph $\mathcal{G}(\mathbb{B},\mathcal{E})$ for the Gram matrix $Q$ and impose its positive semidefiniteness through matrix decomposition.} 

Precisely, let $\mathcal{G}(\mathbb{B},\mathcal{E})$ be a graph with maximal cliques $\mathcal{C}_{1},\ldots,\mathcal{C}_{t}$ and with edge set $\mathcal{E} \subseteq \mathbb{B} \times \mathbb{B}$ satisfying
\begin{equation}\label{e:sos-sparsity-graph-basic-ass}
\mathbb{A} \subseteq \{ \beta + \gamma:\; (\beta,\gamma) \in \mathcal{E} \}.
\end{equation}
Consider the cone of sparse SOS polynomial whose Gram matrix $Q$ has sparsity graph $\mathcal{G}$ and admits the clique-based positive semidefinite decomposition
\begin{equation}\label{e:pop-sparse-matrix-decomposition}
Q = \sum_{k=1}^{t} E_{\mathcal{C}_{k}}^\tr S_{k} E_{\mathcal{C}_{k}}, \qquad
S_{k} \in \mathbb{S}_+^{\abs{\mathcal{C}_{k}}}.
\end{equation}
We denote this cone by
\begin{multline} \label{eq:sparse_sos_general}
\Sigma[\mathbb{A};\mathcal{E}] := \big\{
f \in \Sigma[\mathbb{A}]:\;\;
f(x)=(x^{\mathbb{B}} )^\tr \, Q \, x^{\mathbb{B}},\\
 Q \text{ satisfies~\cref{e:pop-sparse-matrix-decomposition}}
\big\}.
\end{multline}
Conditions~\cref{e:sos-sparsity-graph-basic-ass} and~\cref{e:pop-sparse-matrix-decomposition} imply that $\Sigma[\mathbb{A};\mathcal{E}] \subseteq \Sigma[\mathbb{A}]$. 
Moreover, inserting the clique-based decomposition~\cref{e:pop-sparse-matrix-decomposition} of $Q$ into~\cref{e:sos-sdp} one finds that $f \in \Sigma[\mathbb{A};\mathcal{E}]$ if and only if
\begin{equation}\label{e:sos-sparsified-sdp-unconstrained}
%\quad\iff\quad
%\begin{dcases}
\begin{aligned}
&\exists S_{1} \in \mathbb{S}_+^{\abs{\mathcal{C}_{1}}},\,\ldots,\, S_{t} \in \mathbb{S}_+^{\abs{\mathcal{C}_{t}}}\;\text{ such that}\\
&\sum_{k=1}^t\langle  S_{k} , E_{\mathcal{C}_{k}} A_\alpha E_{\mathcal{C}_{k}}^\tr\rangle = f_\alpha 
\quad \forall \alpha \in \mathbb{B}+ \mathbb{B}.
\end{aligned}
%\end{dcases}
\end{equation}
If the cliques of the prescribed sparsity graph are small, the right-hand side is an SDP with small semidefinite cones and can be solved more efficiently than~\cref{e:sos-sdp}.

\begin{remark}[Chordality of the sparsity graph] \label{remark:chordality_term_sparsity}
	The Gram matrix decomposition~\cref{e:pop-sparse-matrix-decomposition} is motivated by the chordal decomposition result in~\cref{T:ChordalDecompositionTheorem}. However, we do not assume here that the sparsity graph $\mathcal{G}(\mathbb{B},\mathcal{E})$ is chordal, so~\cref{e:pop-sparse-matrix-decomposition} is generally \textit{not} equivalent to requiring $Q \in \mathbb{S}_+^{\abs{\mathbb{B}}}(\mathcal{E},0)$. The lack of chordality makes searching for the maximal cliques $\mathcal{C}_1,\ldots,\mathcal{C}_t$ an NP-hard problem \citep{tomita2006worst}. Allowing for nonchordal graphs with small cliques that can be determined analytically, however, can be extremely useful when a chordal extension leads to unacceptably large cliques even if it is approximately minimal. Examples of this situation can be found in works by \cite{nie2009sparse} and \cite{kovcvara2020decomposition}. \markendexample
\end{remark}

\begin{remark}[Sparse SOS decompositions]\label{remark:sparse-sos-decomposition}
	Given a sparsity graph $\mathcal{G}(\mathbb{B},\mathcal{E})$, the cone $\Sigma[\mathbb{A};\mathcal{E}] \subset \Sigma[\mathbb{A}]$ contains special SOS polynomials that admit a \textit{sparse SOS decomposition}, i.e., a decomposition into a sum of sparse SOS polynomials. Indeed, substituting~\cref{e:pop-sparse-matrix-decomposition} into~\cref{e:sos-gram} yields
	\begin{align}
	f(x)
	&= \sum_{k=1}^{t} (x^{\mathbb{B}})^\tr E_{\mathcal{C}_k}^\tr \, S_{k} \, E_{\mathcal{C}_k} x^{\mathbb{B}} 
	\nonumber \\
	&= \sum_{k=1}^{t} \underbrace{(E_{\mathcal{C}_k} x^{\mathbb{B}})^\tr \,S_{k} \,(E_{\mathcal{C}_k} x^{\mathbb{B}})}_{=:\sigma_k(x)}.
	\label{e:sparse-sos-algebraic}
	\end{align}
	Each polynomial $\sigma_{k}(x)$ is SOS because $S_k$ is positive semidefinite, and is sparse because the operation $E_{\mathcal{C}_{k}} x^{\mathbb{B}}$ extracts a subset of the full monomial vector $x^{\mathbb{B}}$. \markendexample
\end{remark}

It is important to observe that the reduction in computational complexity granted by the clique-based decomposition~\cref{e:pop-sparse-matrix-decomposition} usually comes at the expense of conservatism. %generality. 
This is because sparsity in the support set $\mathbb{A}$ does not guarantee the existence of a sparse Gram matrix $Q$. For a given support set $\mathbb{A}$, special choices of the sparsity graph $\mathcal{G}(\mathbb{B},\mathcal{E})$ may ensure that $\Sigma[\mathbb{A};\mathcal{E}]=\Sigma[\mathbb{A}]$
(\citealp[Corollaries~4.1 \& 4.2]{zheng2020sum};
 \citealp[Theorem~2.1]{mai2020sparse}; 
 \citealp[Theorem~4.1]{Wang2019sos-term-sparsity};
 \citealp[Theorem~3.3]{Wang2020tssos}). 
In general, however, sparse SOS polynomials need not admit a sparse SOS decomposition, so the inclusion $\Sigma[\mathbb{A};\mathcal{E}] \subset \Sigma[\mathbb{A}]$ is \textit{strict}. The next example illustrates this.

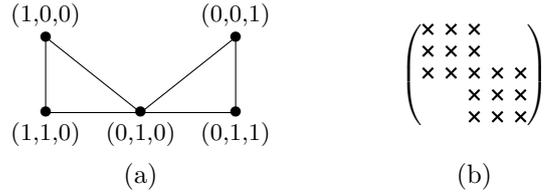
\begin{figure}
	\centering
	\begin{tikzpicture}
	\draw (-1.25,0) -- (-1.25,1);
	\draw (-1.25,0) -- (0,0);
	\draw (0,0) -- (-1.25,1);
	\draw (1.25,0) -- (1.25,1);
	\draw (1.25,0) -- (0,0);
	\draw (0,0) -- (1.25,1);
	\node at (-1.25,0) {$\bullet$};
	\node at (-1.25,1) {$\bullet$};
	\node at (0,0) {$\bullet$};
	\node at (1.25,0) {$\bullet$};
	\node at (1.25,1) {$\bullet$};
	\node[anchor=south] at (-1.25,1) {\small(1,0,0)};
	\node[anchor=north] at (-1.25,0) {\small(1,1,0)};
	\node[anchor=north] at (0,0) {\small(0,1,0)};
	\node[anchor=north] at (1.25,0) {\small(0,1,1)};
	\node[anchor=south] at (1.25,1) {\small(0,0,1)};
	\node at (0,-0.85) {(a)};
	\end{tikzpicture}
	\hspace{1cm}
	\begingroup % keep the change local
	\setlength\arraycolsep{2pt}
	\def\arraystretch{0.75}
	\begin{tikzpicture}
	\node at (0,0) {\small$\begin{pmatrix}
		\mycross{black} & \mycross{black} & \mycross{black} \\
		\mycross{black} & \mycross{black} & \mycross{black} \\
		\mycross{black} & \mycross{black} & \mycross{black} & \mycross{black} & \mycross{black} \\ 
		&&\mycross{black} & \mycross{black} & \mycross{black} \\
		&&\mycross{black} & \mycross{black} & \mycross{black}
		\end{pmatrix}$};
	\node at (0,-1.35) {(b)};
	\end{tikzpicture}
	\endgroup
	\vskip-1ex
	\caption{
		\textit{(a)}~Sparsity graph $\mathcal{G}(\mathbb{B},\mathcal{E})$ for \cref{ex:chordal-sos-failure}. The vertices $\mathbb{B}$ are such that $x^{\mathbb{B}}=(x_1,x_1x_2,x_2,x_3,x_2x_3)$.
		\textit{(b)}~Sparsity pattern of the Gram matrix of SOS polynomials in $\Sigma[\mathbb{A};\mathcal{E}]$.
	}
	\label{f:example-csp-inexact}
\end{figure}

\begin{example}{\cite[Lemma 5.2]{klep2019sparse}}\label{ex:chordal-sos-failure}
	Consider the polynomial $f(x) = x_1^2 - 2x_1x_2 + 3x_2^2 - 2x_1^2x_2 + 2x_1^2x_2^2 - 2x_2x_3 + 6x_3^2 + 18x_2^2x_3 - 54x_2x_3^2 + 142x_2^2x_3^2,$ and set $\mathbb{A} = \supp(f)$. Let $\mathbb{B} \subset \mathbb{N}^3_2$ be the exponent set such that $x^{\mathbb{B}}=(x_1,x_1x_2,x_2,x_3,x_2x_3)$, which is obtained via Newton polytope reduction. Consider also the (chordal) sparsity graph $\mathcal{G}(\mathbb{B},\mathcal{E})$ shown in \cref{f:example-csp-inexact}, which satisfies~\cref{e:sos-sparsity-graph-basic-ass}. We claim that $f$ belongs to $\Sigma[\mathbb{A}]$ but not to $\Sigma[\mathbb{A};\mathcal{E}]$. To see this, observe that any Gram matrix representation of $f$ must take the form
	\begin{equation*}
	f(x) = \begin{pmatrix}x_1\\x_1x_2\\x_2\\x_3\\x_2x_3\end{pmatrix}^\tr
	\underbrace{\small\begin{pmatrix}
		1 & -1 & -1 & 0 & \alpha \\
		-1 & 2 & 0 & -\alpha & 0 \\
		-1 & 0 & 3 & -1 & 9\\
		0 & -\alpha & -1 & 6 & -27 \\
		\alpha & 0 & 9 & -27 & 142
		\end{pmatrix}}_{Q}
	\begin{pmatrix}x_1\\x_1x_2\\x_2\\x_3\\x_2x_3\end{pmatrix},
	\end{equation*}
	where $\alpha \in \mathbb{R}$ can be chosen arbitrarily. Setting $\alpha=1$ makes the Gram matrix $Q$ positive semidefinite, so $f \in \Sigma[\mathbb{A}]$. However, $f$ cannot be in $\Sigma[\mathbb{A};\mathcal{E}]$ because this would require $\alpha=0$, for which $Q$ is not positive semidefinite.  \markendexample
	%This shows that sparsity of the support $\mathbb{A}$ does not imply the existence of a sparse Gram matrix.
\end{example}

\begin{figure*}
	\centering
	\begin{tikzpicture}
	\footnotesize
	\matrix (m3) [matrix of nodes,
	row sep = 3em,	
	column sep = 3em,	
	nodes={circle, draw=black}] at (-3.5,1)
	{1  & 2\\ 4 &  3\\};
	\draw[ultra thick,color=matlabblue] (m3-1-1) -- (m3-1-2);
	\draw[ultra thick,color=matlabred] (m3-1-2) -- (m3-2-2);
	\draw[ultra thick,color=matlabyellow] (m3-2-2) -- (m3-2-1);
	\draw[ultra thick,color=matlabpurple] (m3-2-1) -- (m3-1-1);
	\node at (4,1) {\includegraphics[scale=0.95]{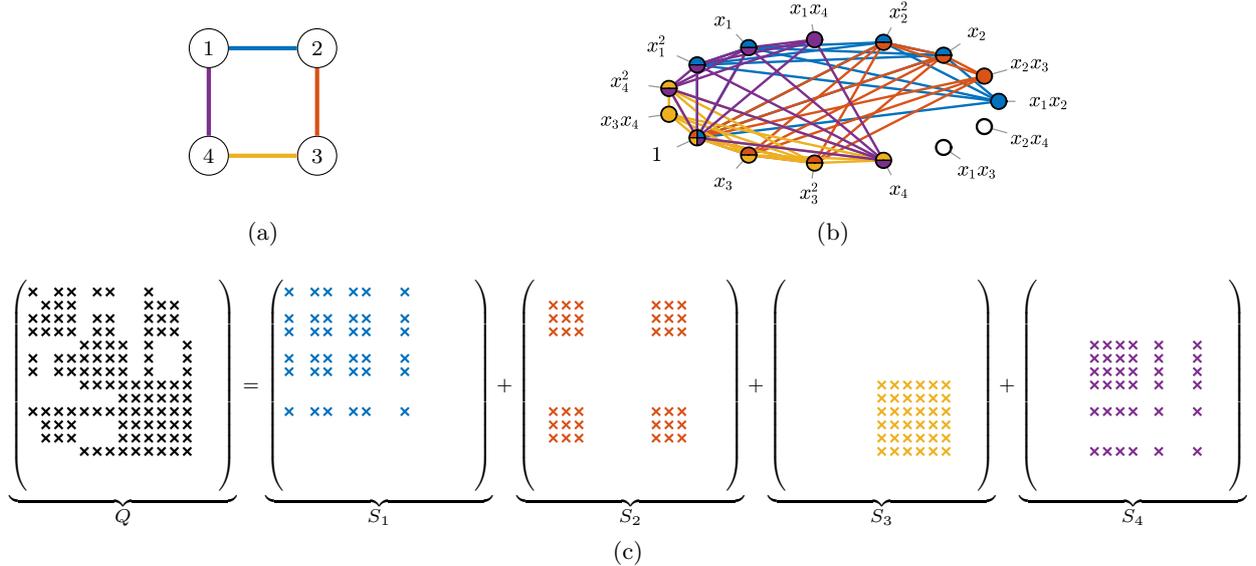}};
	\node at (-3.5,-0.75) {\small(a)};
	\node at (4,-0.75) {\small(b)};
	%\node at (0,-3) {\small(c)};
	\end{tikzpicture}
	\\[1ex]
	\input{./figures/pops/csp-example-arc.tex}
	\caption{
		\textit{(a)}~Correlative sparsity graph of the polynomial in \cref{ex:csp-example}. 
		\textit{(b)}~The corresponding sparsity graph $\mathcal{G}(\mathbb{B},\mathcal{E}_{\csp})$. Graph vertices are labelled by monomials in $x^\mathbb{B}$ instead of the corresponding exponents in $\mathbb{B}$ to ease the visualization. Filled vertices have a self-loop (not shown), empty ones do not. Colors mark the maximal cliques 
		$\mathcal{C}_1$~({\color{matlabblue}\solidrule}), 
		$\mathcal{C}_2$~({\color{matlabred}\solidrule}), 
		$\mathcal{C}_3$~({\color{matlabyellow}\solidrule}) and
		$\mathcal{C}_4$~({\color{matlabpurple}\solidrule}).
		\textit{(c)}~Sparsity pattern of the Gram matrix $Q$ induced by $\mathcal{G}(\mathbb{B},\mathcal{E}_{\csp})$ and its clique-based decomposition. The vertices of $\mathcal{G}(\mathbb{B},\mathcal{E}_{\csp})$ are ordered anticlockwise starting from $x_1x_2$. The last two rows and columns of all matrices are empty.
	}
	\label{f:sos-csp-example}
\end{figure*}

%%%%%%%%%%%%%%%%%%%%%%%%%%%%%%%%%%%%%%%%%%%%%%%%%%%%%%%
\subsubsection{Correlative sparsity}\label{s:sos-csp-unconstrained}
\noindent
The general approach presented in~\Cref{s:sparse-sos-general-strategy} %the previous section 
requires specifying the sparsity graph for the Gram matrix $Q$ in~\cref{e:sos-gram}.
A natural strategy to do this, pioneered by \cite{waki2006sums} and \cite{lasserre2006convergent}, is to consider the couplings between any two independent variables $x_i$ and $x_j$ in a polynomial $f$ supported on $\mathbb{A}$. Two variables $x_i$ and $x_j$ are considered coupled if a monomial in the vector $x^\mathbb{A}$ depends on both simultaneously, i.e., if there exists $\alpha \in \mathbb{A}$ with $\alpha_i \alpha_j>0$. These couplings can be described using the \textit{correlative sparsity (csp) graph} of the support set $\mathbb{A}$ (or, alternatively, of the polynomial $f$), which has vertices $\{1,\ldots,n\}$ and edge set
\begin{equation}
\mathcal{S}_{\csp}(\mathbb{A}):=\{(i,j):\;\exists \alpha \in \mathbb{A} \,\text{ with } \,\alpha_i\alpha_j > 0\}.
\end{equation}

Correlatively sparse SOS decompositions are obtained upon imposing that the entry $Q_{\gamma,\beta}$ of the Gram matrix in~\cref{e:sos-gram} vanishes if the monomial $x^{\beta+\gamma}$ introduces couplings between variables that are not consistent with the csp graph of $f$. This amounts to requiring that $Q$ has sparsity graph $\mathcal{G}_{\csp}(\mathbb{B},\mathcal{E}_{\rm csp})$ with edge set
\begin{multline}\label{e:csp-edge-set}
\mathcal{E}_{\rm csp} := \{ (\beta,\gamma) \in \mathbb{B}\times \mathbb{B}:\; 
\\(\beta_i+\gamma_i)(\beta_j+\gamma_j)> 0 
\Rightarrow (i,j) \in \mathcal{S}_{\csp}(\mathbb{A}) \}.
\end{multline}

One may consider $\mathcal{G}(\mathbb{B},\mathcal{E}_{\csp})$ a ``hypergraph'' with $|\mathbb{B}|$ nodes, built from the csp graph of $f$ (which has $n$ nodes) to ensure that polynomials $(x^\mathbb{B})^\tr \, Q \,x^\mathbb{B}$ with $Q \in \mathbb{S}^{|\mathbb{B}|}(\mathcal{E}_{\csp},0)$ inherit the correlative sparsity of the original support set~$\mathbb{A}$. Unsurprisingly, therefore, the properties of $\mathcal{G}(\mathbb{B},\mathcal{E}_{\csp})$ can be inferred from those of the (usually much smaller) csp graph. In the following statement, which can be proved using arguments similar to those given by \citet[Section~2.4.3]{zheng2019chordal}, $\nnz(\beta)$ denotes the indices of the nonzero entries of an exponent $\beta$.

\begin{proposition}\label{th:sos-csp-cliques}
	Suppose that the csp graph of the support set $\mathbb{A}$ has maximal cliques $\mathcal{J}_1,\ldots,\mathcal{J}_t$. Then, $\mathcal{G}(\mathbb{B},\mathcal{E}_{\csp})$ has maximal cliques $\mathcal{C}_k = \{\beta \in \mathbb{B}: \nnz(\beta) \subseteq \mathcal{J}_k\}$ for $k=1,\ldots,t$.
	Moreover, if the csp graph of $\mathbb{A}$ is chordal, then so is $\mathcal{G}(\mathbb{B},\mathcal{E}_{\csp})$.
\end{proposition}

\Cref{th:sos-csp-cliques} considerably simplifies the construction of the ``inflation'' matrices $E_{\mathcal{C}_{k}}$ in~\cref{e:pop-sparse-matrix-decomposition}, because it suffices to find the maximal cliques of the csp graph of $\mathbb{A}$ without building the (much larger) graph $\mathcal{G}(\mathbb{B},\mathcal{E}_{\csp})$. In addition, it is not difficult to check that the operation $E_{\mathcal{C}_{k}} x^{\mathbb{B}}$ extracts monomials that depend only on variables indexed by $\mathcal{J}_k$. Using \cref{e:sparse-sos-algebraic}, one concludes  that exploiting correlative sparsity amounts to searching for a sparse SOS decomposition in the form
\begin{equation}\label{e:csp-sos-algebraic}
f(x) = \sum_{k=1}^t  \sigma_{k}(x_{\mathcal{J}_k}), \qquad \sigma_k \text{ is SOS},
\end{equation}
where $x_{\mathcal{J}_k}$ denotes the subset of variables $x$ indexed by~$\mathcal{J}_k$ (cf. \citealp[Theorem 2]{zheng2018sparse}). 

\begin{remark}\label{remark:csp-failure}
	\Cref{ex:chordal-sos-failure} shows that correlatively sparse SOS polynomials need not admit the sparse SOS decomposition~\cref{e:csp-sos-algebraic}, even if the csp graph is chordal.
	%(which is a line of three nodes in~\Cref{ex:chordal-sos-failure}). 
	Thus, the inclusion $\Sigma[\mathbb{A};\mathcal{E}_{\csp}] \subset \Sigma[\mathbb{A}]$ is generally strict. For further discussion on the existence of sparse SOS decompositions for polynomials with chordal correlative sparsity, see \cite{mai2020sparse} and \cite{zheng2020sum}. \markendexample
\end{remark}

\begin{example}
	The quartic polynomial $f$ in~\cref{e:facial-reduction-quartic} is correlatively sparse, and the csp graph of its support is chordal with maximal cliques $\mathcal{J}_i = \{i, i+1, i + 2\}$ for $i = 1, \ldots, n-2$. It is clear that $f$ admits a sparse SOS decomposition~\cref{e:csp-sos-algebraic} and this can be searched for by solving the SDP in~\cref{e:sos-sparsified-sdp-unconstrained}. Since, for each clique $\mathcal{J}_i$, only six elements in $\mathbb{B}=\mathbb{N}^n_2 \setminus \mathbb{N}^n_1$ can be multiplied together without introducing spurious couplings to different cliques, this SDP has semidefinite matrix variables $S_1,\ldots,S_{n-2} \in \mathbb{S}^6_+$. Its computational complexity is clearly much lower than the corresponding dense formulation in~\cref{example:facial-reduction}, and a sparse SOS decomposition for $f$ can be found in less than one second on a standard laptop. \markendexample
\end{example}

\begin{example}\label{ex:csp-example}
		Consider the quartic polynomial
		\begin{multline*}
		f(x) = 2 + x_1^2x_{4}^2 (x_1^2x_{4}^2 - 1) - x_1^2 + x_1^4
		\\+ \textstyle\sum_{i=2}^4 \left(x_i^2x_{i-1}^2 (x_i^2x_{i-1}^2 - 1) - x_i^2 + x_i^4\right).
		\end{multline*}
		 Its csp graph, shown in \cref{f:sos-csp-example}(a), is nonchordal and has maximal cliques $\mathcal{J}_1=\{1,2\}$, $\mathcal{J}_2=\{2,3\}$, $\mathcal{J}_1=\{3,4\}$ and $\mathcal{J}_1=\{4,1\}$. The corresponding graph $\mathcal{G}(\mathbb{B},\mathcal{E}_{\csp})$, where the set of exponents obtained with Newton polytope reduction is $\mathbb{B}=\mathbb{N}^4_2$, is shown in \cref{f:sos-csp-example}(b) and has cliques $\mathcal{C}_1,\ldots,\mathcal{C}_4$ containing $6$ elements each, which are determined using \cref{th:sos-csp-cliques}. 
		The sparsity pattern of the Gram matrix $Q$ induced by $\mathcal{G}(\mathbb{B},\mathcal{E}_{\csp})$ and the clique-based matrix decomposition in~\cref{e:pop-sparse-matrix-decomposition}, also illustrated in the figure, replaces a $15\times 15$ positive semidefinite contraint on $Q$ with four  semidefinite constraints on $6\times 6$ matrices $S_1,\ldots,S_4$. According to~\cref{e:csp-sos-algebraic}, searching for these matrices is equivalent to looking for a sparse SOS decomposition $f = \sigma_1(x_1,x_2) + \sigma_2(x_2,x_3) + \sigma_3(x_3,x_4) + \sigma_4(x_4,x_1)$. Such a decomposition is not guaranteed to exist even if $f$ were SOS, but it does for this example with
		%
		%\begin{subequations}
			\begin{align*}
			\sigma_1(x_1,x_2) 
			&= \tfrac12 \left(x_1^2 - \tfrac12 \right)^2 + \left(x_1 x_2 - \tfrac12 \right)^2 + \tfrac12 \left(x_2^2-\tfrac12\right)^2
			\\
			\sigma_2(x_2,x_3) 
			&= \tfrac12 \left(x_2^2 - \tfrac12 \right)^2 + \left(x_2 x_3 - \tfrac12 \right)^2 + \tfrac12 \left(x_3^2-\tfrac12\right)^2
			\\
			\sigma_3(x_3,x_4) 
			&= \tfrac12 \left(x_3^2 - \tfrac12 \right)^2 + \left(x_3 x_4 - \tfrac12 \right)^2 + \tfrac12 \left(x_4^2-\tfrac12\right)^2
			\\
			\sigma_4(x_4,x_1) 
			&= \tfrac12 \left(x_4^2 - \tfrac12 \right)^2 + \left(x_4 x_1 - \tfrac12 \right)^2 + \tfrac12 \left(x_1^2-\tfrac12\right)^2.
			\end{align*}
		%\end{subequations}
		%
		This proves that $f \in \Sigma[\supp(f);\mathcal{E}_{\csp}]$. \markendexample
\end{example}

%%%%%%%%%%%%%%%%%%%%%%%%%%%%%%%%%%%%%%%%%%%%%%%%%%%%%%%
\subsubsection{TSSOS, chordal-TSSOS and related hierarchies}\label{s:tssos}
\noindent
Fix an exponent set $\mathbb{A} \subset \mathbb{N}^n_{2d}$ and a polynomial $f$ with $\supp(f)\subseteq \mathbb{A}$. Correlative sparsity exploits only the sparse couplings between variables as encoded by the csp graph of $\mathbb{A}$, but does not take into account any further structure of $\mathbb{A}$. This is not efficient when $\abs{\mathbb{A}}$ is much smaller than $\binom{n+2d}{2d}$, so $f$ is term-sparse, but the csp graph is fully connected or nearly so. 

For this reason, \cite{Wang2019sos-term-sparsity,Wang2020tssos,Wang2020chordal-tssos} introduced the term-sparse-SOS (TSSOS) and the chordal-TSSOS decomposition hierarchies, which exploit term sparsity irrespective of whether $f$ is correlatively sparse. These are two particular examples of a broader family of possible sparsity-exploiting SOS decomposition hierarchies, each of which is obtained upon imposing the clique-based Gram matrix decomposition~\cref{e:pop-sparse-matrix-decomposition} for a sequence $\{\mathcal{G}(\mathbb{B},\mathcal{E}_k)\}_{k\geq 1}$ of increasingly connected sparsity graphs ($\mathcal{E}_k \subseteq \mathcal{E}_{k+1}$). 

Irrespective of the particular hierarchy being considered (TSSOS, chordal-TSSOS, or another), the construction of such sparsity graphs begins with the observation that, in order to ensure~\cref{e:sos-sparsity-graph-basic-ass}, each edge set $\mathcal{E}_k$ should contain at least all edges $(\beta,\gamma)$ with $\beta + \gamma \in \mathbb{A}$. This guarantees that $\mathbb{A} \subseteq \supp( (x^\mathbb{B})^\tr \,Q \,x^\mathbb{B} )$ for any Gram matrix $Q$ defined via the clique-based decomposition~\cref{e:pop-sparse-matrix-decomposition}, which is necessary for the feasibility of the SDP in~\cref{e:sos-sparsified-sdp-unconstrained}. One should also not force diagonal entries $Q_{\beta\beta}$ of the Gram matrix to vanish, because this amounts to saying that the monomial $x^\beta$ is redundant and $\beta$ could be removed from the exponent set $\mathbb{B}$. For these reasons, we define an initial exponent set $\mathbb{B}_0$ and an initial edge set $\mathcal{E}_0$ as
\begin{subequations}
    \begin{gather}\label{e:tssos-initialization}
    %\mathcal{S}_0 
    \mathbb{B}_0:= 2 \mathbb{B} \cup \mathbb{A},\\
    \mathcal{E}_0 = \{(\beta, \gamma) \in \mathbb{B}\times \mathbb{B}: \beta+\gamma \in \mathbb{B}_{0}\}.
\end{gather}
\end{subequations}
%and the initial edge set $\mathcal{E}_0 = \{(\beta, \gamma) \in \mathbb{B}\times \mathbb{B}: \beta+\gamma \in \mathbb{B}_{0}\}$.

Next, consider an \textit{extension operator}
$\mathbb{E}: \mathbb{B}\times \mathbb{B} \to \mathbb{B}\times \mathbb{B},$
%\begin{equation}
%\begin{aligned}
%\mathbb{E}: \mathbb{B}\times \mathbb{B} &\to \mathbb{B}\times \mathbb{B}\\
%\mathcal{E} &\mapsto \mathbb{E}(\mathcal{E})
%\end{aligned}
%\end{equation}
which extends a given edge set $\mathcal{E} \subset \mathbb{B} \times \mathbb{B}$ according to a given rule. The edge sets $\mathcal{E}_1 \subseteq \mathcal{E}_2 \subseteq \cdots  \subseteq \mathcal{E}_k \subseteq \cdots$ and their corresponding support sets $\mathbb{B}_k$ are defined using the iterative rule
\begin{subequations}
	\label{e:term-sparsity-iterations}
	\begin{align}
	\label{e:tssos-extension}
	\mathcal{E}_k &:= \mathbb{E}\left(
	\{(\beta,\gamma) \in \mathbb{B}\times \mathbb{B}: \beta+\gamma \in \mathbb{B}_{k-1} \}
	\right),\\
	\mathbb{B}_k &:= \{\beta+\gamma: (\beta,\gamma) \in \mathcal{E}_k \}.
	\end{align}
\end{subequations}
Note that $\mathcal{E}_k \subseteq \{(\beta,\gamma) \in \mathbb{B}\times \mathbb{B}: \beta+\gamma \in \mathbb{B}_{k}\}$, so the extension operator guarantees that $\mathcal{E}_{k} \subseteq \mathcal{E}_{k+1}$. Moreover, the sequence $\{\mathcal{E}_k\}_{k\geq 1}$ must converge to an edge set $\mathcal{E}^*$ in a finite number of iterations because $\mathcal{E}_k$ cannot be extended beyond the complete edge set $\mathbb{B}\times \mathbb{B}$.
The sequence of sparsity graphs $\{\mathcal{G}(\mathbb{B},\mathcal{E}_k)\}_{k\geq 1}$ obtained in this way is therefore finite, and yields the (finite) hierarchy of nested sparse SOS cones
\begin{equation}
\Sigma[\mathbb{A};\mathcal{E}_1] \subseteq \Sigma[\mathbb{A};\mathcal{E}_2] \subseteq \cdots \subseteq \Sigma[\mathbb{A};\mathcal{E}^*] \subseteq \Sigma[\mathbb{A}].
\end{equation}
Here, $\Sigma[\mathbb{A};\mathcal{E}]$ is as defined in \cref{eq:sparse_sos_general} and all inclusions are strict in general.

Different extension operators produce different types of sparse SOS decomposition hierarchies. In particular:
\begin{itemize}
	\item If $\mathbb{E}$ is a \textit{block-completion} operator that completes all connected components of the edge set $\{(\beta,\gamma) \in \mathbb{B}\times \mathbb{B}: \beta+\gamma \in \mathbb{B}_{k-1} \}$, one recovers the TSSOS hierarchy \citep{Wang2019sos-term-sparsity,Wang2020tssos}. 
	At each step of the hierarchy, $Q$ has chordal sparsity (specifically, a block-diagonal structure) and~\cref{e:pop-sparse-matrix-decomposition} is equivalent to imposing $Q \in \mathbb{S}^{\abs{\mathbb{B}}}_+(\mathcal{E}_k,0)$.
	\item If $\mathbb{E}$ is an \textit{approximately minimal chordal extension} operator that extends the edge sets $\{(\beta,\gamma) \in \mathbb{B}\times \mathbb{B}: \beta+\gamma \in \mathbb{B}_{k-1} \}$ such that $\mathcal{G}(\mathbb{B},\mathcal{E}_k)$ is chordal, one recovers the chordal-TSSOS hierarchy \citep{Wang2020chordal-tssos}. 
	At each step of the hierarchy, $Q$ has chordal sparsity and~\cref{e:pop-sparse-matrix-decomposition} is equivalent to requiring $Q \in \mathbb{S}^{\abs{\mathbb{B}}}_+(\mathcal{E}_k,0)$.
\end{itemize} 
%At each step of the TSSOS and chordal-TSSOS hierarchies, $Q$ has chordal sparsity (in fact, for the TSSOS hierarchy, it is block-diagonal) and by \cref{T:ChordalDecompositionTheorem} condition~\cref{e:pop-sparse-matrix-decomposition} is equivalent to requiring $Q \in \mathbb{S}^{\abs{\mathbb{B}}}_+(\mathcal{E}_k)$.
In both cases, the edge extensions are performed on a graph with $|\mathbb{B}|$ nodes and the maximal cliques of $\mathcal{G}(\mathbb{B},\mathcal{E}_k)$ must be found at each iteration. This is unlike the correlative sparsity strategy in~\cref{s:sos-csp-unconstrained}, where the maximal cliques of $\mathcal{G}(\mathbb{B},\mathcal{E}_{\csp})$ are built from those in the csp graph of $\mathbb{A}$, which has only $n$ nodes (cf. \cref{th:sos-csp-cliques}).

It is also clear that the choice of extension operator determines the computational complexity of the resulting sparse SOS decomposition hierarchy, as well as the gap between $\Sigma[\mathbb{A},\mathcal{E}^*]$ and $\Sigma[\mathbb{A}]$.
For example, the chordal-TSSOS hierarchy has a lower complexity than the TSSOS one in general, as its sparsity graphs have fewer edges (see \citealp{Wang2020tssos,Wang2020chordal-tssos} for detailed complexity estimates). However, the TSSOS hierarchy has a higher representation power because  $\Sigma[\mathbb{A},\mathcal{E}^*]=\Sigma[\mathbb{A}]$, which is generally not true for the chordal-TSSOS hierarchy.

\begin{theorem}[\citealp{Wang2020tssos}]\label{th:tssos-sign-symm}
	%Let $\mathbb{E}$ be the block-completion operator and 
	If $\mathcal{E}^*_{\text{\sc tssos}}$ is the stabilized edge set of the TSSOS hierarchy, then $\Sigma[\mathbb{A},\mathcal{E}^*_{\text{\sc tssos}}]=\Sigma[\mathbb{A}]$, i.e., $f$ is SOS if and only if $f \in \Sigma[\mathbb{A},\mathcal{E}^*_{\text{\sc tssos}}]$.
\end{theorem}
\begin{remark}\label{remark:tssos-sign-symmetries}
\Cref{th:tssos-sign-symm} follows from a stronger result~\cite[ Theorem~6.5]{Wang2020tssos} which reveals that the constraint $Q \in \mathbb{S}^{\abs{\mathbb{B}}}_+(\mathcal{E}^*_{\text{\sc tssos}},0)$ imposes the well-known block-diagonal structure implied by the \textit{sign symmetries} of $f$ (see, e.g.,~\citealp{lofberg2009pre}). \markendexample
\end{remark}

\begin{figure}[t]
	\centering
	\input{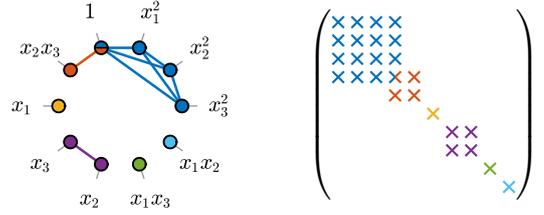}
	\caption{Sparsity graph (left) and corresponding matrix sparsity pattern (right) for the chordal-TSSOS hierarchy in \cref{ex:tssos}, which stabilizes at the first iteration ($k=1$). Graph vertices are labelled by monomials in $x^\mathbb{B}$ instead of the corresponding exponents in $\mathbb{B}$ to ease the visualization. Colours mark the maximal cliques; multicoloured vertices and matrix entries belong to multiple cliques.
		\label{f:ex-chordal-tssos}
	}
\end{figure}

\begin{example}\label{ex:tssos}
	The trivariate quartic polynomial
	$$f(x) = 1 + x_1^4 + x_2^4 + x_3^4 + x_1^2 x_2^2 + x_1^2 x_3^2 + x_2^2 x_3^2 + x_2 x_3$$
	is term sparse but not correlatively sparse, since its csp graph is a complete graph with three nodes. The candidate exponent set to search for an SOS decomposition of $f$ is $\mathbb{B}=\mathbb{N}^3_2$, as Newton polytope reduction removes no terms. For convenience, we order $\mathbb{B}$ such that 
	$$x^\mathbb{B} = (x_3^2, x_2^2, x_1^2, x_2x_3, 1, x_1, x_1x_3, x_1x_2, x_3, x_2)^\tr.$$
	
	\begin{figure}[t]
		\centering
		\input{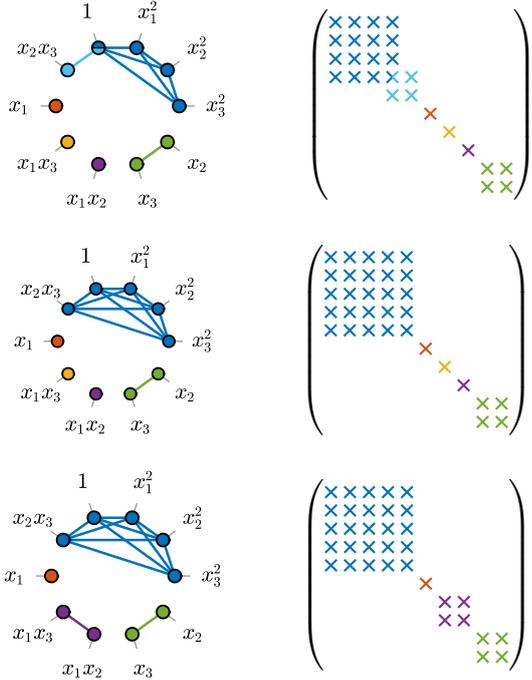}
		\caption{Sparsity graphs and corresponding matrix sparsity patterns for the TSSOS hierarchy in \cref{ex:tssos} at initialization (top; edge set $\mathcal{E}_0$), at the first iteration (middle; edge set $\mathcal{E}_1$) and at the second iteration (bottom; edge set $\mathcal{E}_2$). After that, the hierarchy stabilizes. Graph vertices are labelled by monomials in $x^\mathbb{B}$ instead of the corresponding exponents in $\mathbb{B}$ to ease the visualization. Colors mark the maximal cliques. 
			\label{f:ex-tssos}
		}
	\end{figure}
	
	The TSSOS hierarchy yields the sparsity graphs shown in \cref{f:ex-tssos}, which stabilize at the second iteration ($k=2$). The corresponding sparsity patterns of the Gram matrix $Q$ are also shown in that figure. Observe how the connected components of the initial graph $\mathcal{G}(\mathbb{B},\mathcal{E}_0)$ are completed at the first iteration to obtain the graph $\mathcal{G}(\mathbb{B},\mathcal{E}_1)$.
	%At the fist iteration, the constraint $Q \in \mathbb{S}^{\abs{\mathbb{B}}}(\mathcal{E}_1,0)$ decomposes into five semidefinite constraints on matrices of size $5\times 5$, $1\times 1$, $2 \times 2$, $1 \times 1$, and $1 \times 1$. 
	As discusses in \Cref{remark:tssos-sign-symmetries}, the stabilized block-diagonal structure of $Q$ coincides with the partition of $x^\mathbb{B}$ into the groups $\{x_3^2, x_2^2, x_1^2, 1, x_2x_3\}$, $\{x_1\}$, $\{x_3, x_2\}$ and $\{x_1x_2, x_1x_3\}$ implied by the sign symmetries of $f$, which is invariant under the transformations $(x_1,x_2,x_3) \mapsto (-x_1,-x_2,-x_3)$ and $(x_1,x_2,x_3) \mapsto (x_1,-x_2,-x_3)$ (the four groups of monomials are invariant under both, the first, the second, and none of these transformations). In this example, the SDP in~\cref{e:sos-sparsified-sdp-unconstrained} is feasible at all iterations of the TSSOS hierarchy because $f$ admits the positive semidefinite Gram matrix representation
	\begin{equation*}%\label{e:tsos-example-explicit-decomp}
	f(x) = \frac18
	%\begin{pmatrix}x_3^2\\x_2^2\\x_1^2\\1\\x_2 x_3\\x_1\\x_1x_3\\x_1x_2\\x_3\\x_2\end{pmatrix}^\tr
	(x^\mathbb{B})^\tr
	{\small
	\begin{pmatrix}
	8 & 3 & 4 & 0 & 0 & 0 & 0 & 0 & 0 & 0\\
	3 & 8 & 4 & 0 & 0 & 0 & 0 & 0 & 0 & 0\\
	4 & 4 & 8 & 0 & 0 & 0 & 0 & 0 & 0 & 0\\
	0 & 0 & 0 & 2 & 4 & 0 & 0 & 0 & 0 & 0\\
	0 & 0 & 0 & 4 & 8 & 0 & 0 & 0 & 0 & 0\\
	0 & 0 & 0 & 0 & 0 & 0 & 0 & 0 & 0 & 0\\
	0 & 0 & 0 & 0 & 0 & 0 & 0 & 0 & 0 & 0\\
	0 & 0 & 0 & 0 & 0 & 0 & 0 & 0 & 0 & 0\\
	0 & 0 & 0 & 0 & 0 & 0 & 0 & 0 & 0 & 0\\
	0 & 0 & 0 & 0 & 0 & 0 & 0 & 0 & 0 & 0
	\end{pmatrix}
	}
	%\begin{pmatrix}x_3^2\\x_2^2\\x_1^2\\1\\x_2 x_3\\x_1\\x_1x_3\\x_1x_2\\x_3\\x_2\end{pmatrix}
	x^\mathbb{B}
	\end{equation*}
	and the Gram matrix is consistent with the sparsity graphs in \cref{f:ex-tssos}. Thus, all steps of the TSSOS hierarchy are able to prove that $f$ is SOS. Note that this can be guaranteed a priori only for the last step by virtue of \cref{th:tssos-sign-symm}.
	
	For the same polynomial $f$, the chordal-TSSOS hierarchy stabilizes at the first iteration ($k=1$) and yields the sparsity graph shown in \cref{f:ex-chordal-tssos}. The corresponding Gram matrix $Q$ is sparser than those encountered in the TSSOS hierarchy, leading to smaller semidefinite constraints in~\cref{e:sos-sparsified-sdp-unconstrained}. Again, this SDP is feasible in light of the Gram matrix decomposition given above, so the chordal-TSSOS hierarchy is able to prove that $f$ is SOS. This, however, cannot be guaranteed a priori. \markendexample
\end{example}

\begin{remark}
	The explicit Gram matrix decomposition in \cref{ex:tssos} reveals that the smaller monomial basis $x^\mathbb{B} = (x_3^2, x_2^2, x_1^2, x_2 x_3, 1)$ would suffice to construct an SOS decomposition of $f$. It remains to be seen whether this reduced basis can be identified using strategies that are more sophisticated than the Newton polytope reduction.
\end{remark}

%%%%%%%%%%%%%%%%%%%%%%%%%%%%
\subsubsection{Correlatively term-sparse hierarchies}\label{s:cs-tssos}
\noindent
The sparse SOS decomposition hierarchies described in \cref{s:tssos} can be combined with the correlative sparsity techniques outlined in \cref{s:sos-csp-unconstrained} in a natural way. Let $\mathcal{E}_k$ be the edge set obtained using the iterations in~\cref{e:term-sparsity-iterations} for a given extension operator, and let $\mathcal{E}_{\csp}$ be the edge set in~\cref{e:csp-edge-set} constructed using correlative sparsity. Then, the sequence of sparsity graphs
\begin{equation}\label{e:corr-term-sp-graphs}
\mathcal{G}(\mathbb{B},\mathcal{E}_k\cap \mathcal{E}_{\csp}), \quad k\geq 1
\end{equation}
yields a hierarchy of ``correlatively term-sparse'' SOS decompositions, which exploit simultaneously term and correlative sparsity. Since $\mathcal{E}_k \subseteq \mathcal{E}_{k+1}$ by construction, and since the sequence $\{\mathcal{E}_k\}$ stabilizes onto an edge set $\mathcal{E}^*$ in a finite number of steps, the sparse SOS cones corresponding to this hierarchy satisfy
\begin{multline*}
\Sigma[\mathbb{A};\mathcal{E}_1\cap\mathcal{E}_{\csp}] \subseteq 
\Sigma[\mathbb{A};\mathcal{E}_2\cap\mathcal{E}_{\csp}] \subseteq
\cdots 
\\
\cdots\subseteq
\Sigma[\mathbb{A};\mathcal{E}^*\cap\mathcal{E}_{\csp}]
\subseteq
\begin{cases}
\Sigma[\mathbb{A};\mathcal{E}^*]\\
\Sigma[\mathbb{A};\mathcal{E}_{\csp}],
\end{cases}
\end{multline*}
and all inclusions are generally strict.
Note also that one may remove from $\mathbb{B}$ all exponents that violate the correlative sparsity \textit{before} constructing the edge sets $\mathcal{E}_k$, because the intersection with $\mathcal{E}_{\csp}$ eliminates all edges between such exponents (including self-loops). 

When the extension operator used to build $\mathcal{E}_k$ is the block-completion operation used in the TSSOS hierarchy, the sparsity graphs in~\cref{e:corr-term-sp-graphs} yield exactly the CS-TSSOS hierarchy introduced by \cite{Wang2020cs-tssos}. In this case, by \cref{th:tssos-sign-symm}, the stabilized sparsity graph $\mathcal{G}(\mathbb{B},\mathcal{E}^*\cap \mathcal{E}_{\csp})$ simply encodes sign symmetries and correlative sparsity. Since exploiting sign symmetries in SOS decompositions brings no conservatism \citep{lofberg2009pre}, one immediately obtains the following corollary.

\begin{proposition}
	If $\mathcal{E}^*_{\text{\sc tssos}}$ is the stabilized edge set of the TSSOS hierarchy, then $\Sigma[\mathbb{A};\mathcal{E}^*_{\text{\sc tssos}}\cap\mathcal{E}_{\csp}] = \Sigma[\mathbb{A};\mathcal{E}_{\csp}]$ for any exponent set $\mathbb{A}\subseteq \mathbb{N}^n_{2d}$.
\end{proposition}

\begin{figure}[t]
	\centering
	\input{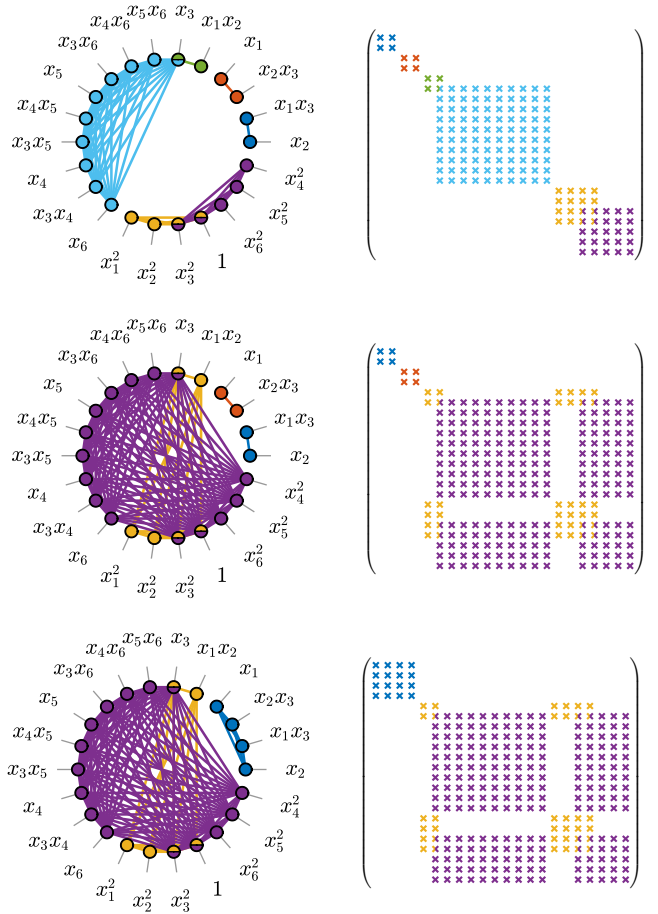}
	\caption{Sparsity graphs and corresponding matrix sparsity patterns for the first ({top}),  second ({middle}) and third ({bottom}) iterations of the CS-TSSOS hierarchy in \cref{ex:cs-tssos}. After that, the hierarchy stabilizes. Graph vertices are labelled by monomials in $x^\mathbb{B}$ instead of the corresponding exponents in $\mathbb{B}$ to ease the visualization. Colors mark the maximal cliques; multicolor vertices and matrix entries belong to multiple cliques.
		\label{f:example-cs-tssos}}
\end{figure}

\begin{example}{\citep[Example~3.4]{Wang2020cs-tssos}}\label{ex:cs-tssos}
	Let%Consider the polynomial
	\begin{multline}\label{e:cs-tssos-example}
	f(x) = 1 + x_1x_2x_3 + x_3x_4x_5 + x_3x_4x_6 \\[-2ex]
	+ x_3x_5x_6 + x_4x_5x_6 + \sum_{i=1}^6 x_i^4.
	\end{multline}
	This polynomial is both term and correlatively sparse, and its csp graph has two maximal cliques $\mathcal{J}_1 = \{1,2,3\}$ and $\mathcal{J}_2 = \{3,4,5,6\}$. It is also invariant under the sign symmetry transformation $(x_1,x_2,x_3,x_4,x_5,x_6) \mapsto (-x_1,-x_2,x_3,x_4,x_5,x_6)$. To search for an SOS decomposition of $f$ using the CS-TSSOS hierarchy, we let $\mathbb{B}$ be the exponent set defining the monomial vector % (monomials violating the correlative sparsity are removed)
	\begin{multline*}
	x^\mathbb{B} = (x_2, x_1x_3, x_2x_3,x_1,x_1x_2,x_3,x_5x_6,x_4x_6,x_3x_6,\\x_5,x_4x_5,x_3x_5,x_4,x_3x_4,x_6,x_1^2,x_2^2,x_3^2,1,x_6^2,x_5^2,x_4^2).
	\end{multline*}
	This is obtained upon removing from the full basis $x^{\mathbb{N}^6_2}$ all monomials that violate the correlative sparsity of $f$ (these would be removed anyway by the CS-TSSOS hierarchy).
	
	The CS-TSSOS sparsity graphs obtained with~\cref{e:corr-term-sp-graphs} and the corresponding sparsity patterns for the Gram matrix $Q$ of $f$, illustrated in \cref{f:example-cs-tssos}, are chordal. The hierarchy stabilizes after three steps. At the first step, the clique-based decomposition~\cref{e:pop-sparse-matrix-decomposition} replaces the $22\times22$ semidefinite constraint on the Gram matrix $Q$ with six semidefinite constraints of size 2, 2, 2, 10, 4 and 5. The first step of the TSSOS hierarchy, instead, leads to an SDP with five semidefinite constraints of size 2, 2, 2, 10, 7~\cite[Example~3.4]{Wang2020cs-tssos}. The second iteration of the CS-TSSOS hierarchy produces significant fill-in, and the size of the largest semidefinite constraint increases to 15. The third iteration brings only minimal additional fill-in. At this final stage, the connected components of the sparsity graph correspond to a partition of the monomials $x^\mathbb{B}$ according to the sign symmetry of $f$ (the first four monomials in $x^\mathbb{B}$ are not invariant under the symmetry transformation, while the rest are), but the correlative sparsity prevents the completion of the largest connected component, i.e., of the bottom-right connected matrix block in \cref{f:example-cs-tssos}\textit{(c)}.
	
	Numerical solution of the SDP~\cref{e:sos-sparsified-sdp-unconstrained} shows that all steps of the CS-TSSOS hierarchy are feasible, so an SOS decomposition of the polynomial $f$ in~\cref{e:cs-tssos-example} can be constructed at a lower computational cost than any other hierarchy discussed in this work. Note that feasibility cannot be guaranteed a priori at any step of the hierarchy, even the last (stabilized) one, due to the conservative nature of correlatively sparse SOS decomposition (see \Cref{remark:csp-failure}).  \markendexample
\end{example}

%%%%%%%%%%%%%%%%%%%%%%%%%%%%%%%%%%%%%%%%%
\subsubsection{Sparse SOS decompositions on semialgebraic sets}\label{s:pop-sparse-constrained}
\noindent
%Brief outline of what is needed to extend sparsity exploitation to the constrained case: take constraints into account carefully. The notation here is going to be quite painful -- perhaps only describe the main idea for the correlative sparsity case (choose sparsity to preserve couplings between independent variables)?
The sparsity-exploiting methods to construct SOS decompositions described so far prove global polynomial nonnegativity, but can be extended to establish local nonnegativity on a basic semialgebraic set defined by $m$ polynomial inequalities,
\begin{equation}\label{e:semialgebraic-set}
\mathbb{K} := \{ x \in \mathbb{R}^n : g_1(x) \geq 0, \ldots, g_m(x) \geq 0 \}.
\end{equation} 

Set $g_0(x)\equiv1$ for convenience. To verify that $f \in \mathbb{R}[x]_{n,d}$ (not necessarily of even degree) is nonnegative on $\mathbb{K}$, it suffices to find an integer $\omega$ (known as the \textit{relaxation order}) and exponent sets $\mathbb{B}_0,\ldots,\mathbb{B}_m \subseteq \mathbb{N}^n_{\omega}$ such that
\begin{equation}\label{e:wsos}
f(x) = \sum_{i=0}^{m} g_i(x) (x^{\mathbb{B}_i})^\tr \, Q_i \, x^{\mathbb{B}_i}, \qquad
Q_i \in \mathbb{S}^{\abs{\mathbb{B}_i}}_+.
\end{equation}
As before, these conditions define the feasible set of an SDP.
Generally, one chooses $\omega$ such that
\begin{equation*}
2 \omega \geq \max\{ \deg(f), \deg(g_1),\ldots, \deg(g_m)\}
\end{equation*}
and then takes
\begin{equation*}
\mathbb{B}_i = \mathbb{N}^n_{\omega_i}, \qquad \omega_i := \omega - \lceil \tfrac12 \deg(g_i) \rceil,
\end{equation*}
where $\lceil a \rceil$ is the smallest integer greater than or equal to $a$. This ensures that each term in the sum in~\cref{e:wsos} is a polynomial of degree at most $2\omega$. One can allow for $2\omega > \deg(f)$ because cancellations may occur when summing all terms.

Since each polynomial $\sigma_i(x) := (x^{\mathbb{B}_i})^\tr \, Q_i \, x^{\mathbb{B}_i}$ in \cref{e:wsos} is SOS, this condition gives a weighted SOS decomposition of $f$, meaning a representation of $f$ as a weighted sum of SOS polynomials where the weights are $g_0=1$ and the polynomials $g_1,\ldots,g_m$ appearing in the semialgebraic definition of $\mathbb{K}$. 
%Note that, for $m=0$, $\mathbb{K}=\mathbb{R}^n$ and one recovers the Gram matrix decomposition~\cref{e:sos-gram}.
Remarkably, this sufficient condition for local nonnegativity is also necessary if $f$ is strictly positive on $\mathbb{K}$ and this is a compact set satisfying the so-called Archimedean condition.

\begin{assumption}[Archimedean condition]\label{ass:archimedean}
	There exists an integer $\nu\geq 0$, SOS polynomials $\sigma_0,\ldots,\sigma_m \in \Sigma_{n,2\nu}$, and a constant $r \in \mathbb{R}$ such that $r^2 - \|x\|_2^2 = \sum_{i=0}^m g_i(x) \sigma_i(x)$.
\end{assumption}

\begin{theorem}[\citealp{putinar1993positive}]\label{th:putinar}
Suppose that $f \in \mathbb{R}[x]_{n,d}$ is strictly positive on a basic semialgebraic set $\mathbb{K}$ defined as in~\cref{e:semialgebraic-set} that satisfies the \nameref{ass:archimedean}. Then, there exists a relaxation order $\omega$ such that $f$ admits the weighted SOS decomposition~\cref{e:wsos}.
\end{theorem}

If $f$ and the polynomials $g_1,\ldots,g_m$ are term sparse, one can proceed as in \cref{s:sparse-sos-general-strategy} and attempt to reduce the computational complexity of the generic weighted SOS decomposition~\cref{e:wsos} by requiring the matrices $Q_0,\ldots,Q_m$ to be sparse and to admit a clique-based positive semidefinite matrix decomposition. The only new aspect is that one must consider how the sparse polynomial $(x^{\mathbb{B}_i})^\tr \, Q_i \, x^{\mathbb{B}_i}$ interacts with the corresponding $g_i$ in order to determine the overall structure of the sum on the right-hand side of~\cref{e:wsos}. This requires some care, especially if one hopes to recover sparse versions of \cref{th:putinar}.

To give an example of this general strategy, let us explain how to extend the correlative sparsity technique outlined in \cref{s:sos-csp-unconstrained}. In this case, one replaces the csp graph of $f$ constructed with the \textit{joint csp graph} of the polynomials $f,g_1,\ldots,g_m$, which has vertices $\{1,\ldots,n\}$ and an edge between vertices $i$ and $j$ if at least one of the following conditions hold:
\begin{enumerate}[{(a)}]
	\item\label{condition:csp-f} The variables $x_i$ and $x_j$ are multiplied together in $f$; %as in \cref{s:sos-csp-unconstrained},
	\item\label{condition:csp-gi} At least one of $g_1,\ldots,g_m$ depends on both $x_i$ and $x_j$, even if these variables are not multiplied together.
\end{enumerate}
The different treatment of $f$ and $g_1,\ldots,g_m$ reflects the asymmetric role these polynomials play in~\cref{e:wsos}.
Then, one imposes that each matrix $Q_i$ in~\cref{e:wsos} is the densest possible matrix such that the support of $g_i(x)(x^{\mathbb{B}_i})^\tr \, Q_i \, x^{\mathbb{B}_i}$ is consistent with the joint csp graph. Precisely, let $\mathcal{J}_1,\ldots,\mathcal{J}_t$ be the maximal cliques of the joint csp graph and, for each $i=0,\ldots,m$, let $\var(g_i) \subset \{1,\ldots,n\}$ be the set of indices of the variables on which $g_i$ depends, with the convention that $\var(g_0)=\var(1)=\emptyset$. By condition~\ref{condition:csp-gi} above, there is at least one clique $\mathcal{J}_k$ such that $\var(g_i) \subseteq \mathcal{J}_k$ and we denote the set of clique indices $k$ for which this holds by
\begin{equation}
\mathcal{N}_i := \left\{ k \in \{1,\ldots,t\}:\; \var(g_i) \subseteq \mathcal{J}_k \right\}.
\end{equation}
Observe in particular that $\mathcal{N}_0 = \{1,\ldots,t\}$ since $\var(g_0) = \emptyset \subset \mathcal{J}_k$ for all $k = 1, \ldots, t$.
%For notational convenience, we have $\var(g_0) = \emptyset \subseteq \mathcal{J}_k$ for all $k = 1, \ldots, t$ and $\mathcal{N}_0 = \{\mathcal{J}_1, \ldots, \mathcal{J}_t\}$.  
%
The sparsity graph $\mathcal{G}_i(\mathbb{B}_i,\mathcal{E}_i)$ of $Q_i$ is defined to have the edge set
\begin{equation}
\mathcal{E}_i := 
\!\!\bigcup_{k \in \mathcal{N}_i}%(g_i)} 
\!\!
\left\{ (\beta,\gamma) \in \mathbb{B}_i \times\mathbb{B}_i:\, \nnz(\beta+\gamma) \subseteq \mathcal{J}_k  \right\}.
\end{equation}
One can check that $\mathcal{G}_i(\mathbb{B}_i,\mathcal{E}_i)$ is chordal if so is the joint csp graph of $f,g_1,\ldots,g_m$. Moreover, it has maximal cliques $\mathcal{C}_{i,1},\ldots,\mathcal{C}_{i,\abs{\mathcal{N}_i}}$ with $\mathcal{C}_{i,k} := \{\beta \in \mathbb{B}_i:\; \nnz(\beta) \subseteq \mathcal{J}_k \}$. Consequently, the clique-based positive semidefinite decomposition of $Q_i$ reads
\begin{align}
Q_i = \sum_{k=1}^{\abs{\mathcal{N}_i}} E_{\mathcal{C}_{i,k}}^\tr S_k E_{\mathcal{C}_{i,k}},
&&
S_k \in \mathbb{S}^{\abs{\mathcal{C}_{i,k}}}_+.
\end{align}

If the cliques $\mathcal{C}_{i,k}$ are small, imposing this clique-based decomposition for each $i=0,\ldots,m$  in~\cref{e:wsos} allows one to search for a weighted SOS decomposition of $f$ by solving an SDP with low computational complexity. Moreover, arguing as in \Cref{remark:sparse-sos-decomposition}, one concludes that this process yields the representation
\begin{equation}\label{e:sparse-wsos-csp}
f(x) = \sum_{i=0}^m \sum_{k \in \mathcal{N}_i}  g_i(x)  \sigma_{i,k}(x_{\mathcal{J}_k}),
\end{equation}
where each SOS polynomial $\sigma_{i,k}$ depends only on variables indexed by a single clique of the joint csp graph. Crucially, the following sparse version of \cref{th:putinar} guarantees that such a sparse weighted SOS decomposition exists if the joint csp graph is chordal, the semialgebraic definition of the set $\mathbb{K}$ in~\cref{e:semialgebraic-set} includes inequalities of the form $r_k^2 - \|x_{\mathcal{J}_k}\|_2^2 \geq 0$ for all $k=1,\ldots,t$, and $f>0$ on $\mathbb{K}$.

\begin{theorem}[\citealp{lasserre2006convergent,grimm2007note}]\label{th:csp-putinar-scalar}
	Let $f$ be a polynomial that is strictly positive on a basic semialgebraic set $\mathbb{K} = \{x \in \mathbb{R}^n:g_1(x)\geq 0,\ldots,g_m(x)\geq 0\}$, whose definition includes the inequalities $r_k^2 - \|x_{\mathcal{J}_k}\|^2 \geq 0$ for some constants $r_1,\ldots,r_t$ and all $k=1,\ldots,t$. If the joint csp graph of $f,g_1,\ldots,g_m$ is chordal, $f$ has a sparse weighted SOS decomposition in the form~\cref{e:sparse-wsos-csp}.
\end{theorem}
\begin{remark}
    The assumption that the semialgebraic definition of $\mathbb{K}$ includes the inequalities $r_k^2 - \|x_{\mathcal{J}_k}\|_2^2 \geq 0$ can be weakened by requiring that the $\abs{\mathcal{J}_k}$-dimensional set
    \begin{equation*}
        \mathbb{K}_k := \{
        \hat{x} \in \mathbb{R}^{\abs{\mathcal{J}_k}}:\; g_i(\hat{x}) \geq 0 \;\; \forall i \;\text{s.t.}\; \var(g_i) \subseteq \mathcal{J}_k
        \}
    \end{equation*}
    satisfies the \nameref{ass:archimedean} for each $k=1,\ldots,t$. Moreover, the assumption is mild when $\mathbb{K}$ is compact because, in principle, the inequalities $r_k^2 - \|x_{\mathcal{J}_k}\|^2 \geq 0$ can be added with values of $r_k$ large enough not to change the set $\mathbb{K}$. Proving that $\mathbb{K}$ remains unchanged for candidate $r_k$, however, may not be easy in practice. \markendexample
\end{remark}

The TSSOS, chordal-TSSOS and CS-TSSOS hierarchies can also be extended to produce weighted SOS decomposition on basic semialgebraic sets \citep{Wang2020tssos,Wang2020chordal-tssos,Wang2020cs-tssos}. Interested readers are referred to these works for the details. Here, we simply observe that, just like their global counterparts described in \cref{s:tssos,s:cs-tssos}, these extended hierarchies stabilize after a finite number of steps. Upon stabilization, moreover, the extended TSSOS and CS-TSSOS hierarchies recover the block-diagonal structure of the matrices $Q_0,\ldots,Q_m$ implied by joint sign symmetries of the polynomials $f,g_1,\ldots,g_m$ 
(see \citealp[Theorem~6.5 and Corollary 6.8]{Wang2020tssos};
\citealp[Proposition 3.10]{Wang2020cs-tssos}). This observation can be combined with a symmetry-exploiting version of \cref{th:putinar} \citep[Theorem~3.5]{Riener2013} and with \cref{th:csp-putinar-scalar} to conclude that the TSSOS and CS-TSSOS hierarchies are \textit{guaranteed} to work for term-sparse polynomials that are strictly positive on compact sets whose semialgebraic definition satisfies suitable versions of the Archimedean condition.
%(cf.~\citealp[Theorem~6.11]{Wang2020tssos} and \citealp[Theorem~4.2]{Wang2020cs-tssos}).

%%%%%%%%%%%%%%%%%%%%%%%%%%%%%%%%%
\subsection{Decomposition of sparse polynomial matrices}\label{s:pop-sparse-matrix}
\noindent
Having studied sparsity-exploiting techniques to reduce the complexity of searching for SOS representations for term-sparse polynomials, we now switch gear and review how chordal sparsity can be exploited when looking for SOS representations of sparse polynomial matrices. \Cref{ss:global-matrix-sos} presents results by \cite{zheng2020sum} that partially extend the classical chordal decomposition theorem (\cref{T:ChordalDecompositionTheorem}) to SOS polynomial matrices with chordal sparsity. Decomposition results giving SOS certificates of matrix positivity on semialgebraic sets are briefly outlined in \cref{ss:local-matrix-sos}. All of these results are useful for static output controller design~\citep{henrion2006convergent}, robust stability region analysis~\citep{henrion2011inner}, and stability analysis of time-delay systems~\citep{peet2009positive}. %enable [mention some control-related applications, with references]~\cite{henrion2006convergent,henrion2011inner,peet2009positive}}. 
Note that all results presented in this section consider the \textit{structural} sparsity of polynomial matrices, \textit{not} their term sparsity. In principle, one could exploit both structural and term sparsity by combining the results reviewed below with those of \cref{s:pop-sparse-global}.

%%%%%%%%%%%%%%%%%%%%%%%%%%%%%%
\subsubsection{Global decomposition}\label{ss:global-matrix-sos}
\noindent
Consider a symmetric $n$-variate polynomial matrix $P \in \mathbb{R}[x]_{n,2d}^{r \times r}$ of degree $2d$ whose (structural) sparsity pattern is described by a chordal graph $\mathcal{G}(\{1,\ldots,r\},\mathcal{E})$, i.e.,
\begin{equation}
(i,j) \notin \mathcal{E} \quad \implies \quad P_{ij}(x) \equiv 0 \quad \forall x \in \mathbb{R}^n.
\end{equation}
Since checking whether $P$ is positive semidefinite globally via the SOS certificates described in \cref{s:matrix-sos} is expensive when $r$ is large, we seek to exploit the sparsity of $P$ and replace one large matrix SOS constraint with multiple smaller ones.

Let $\mathcal{C}_1,\ldots,\mathcal{C}_t$ be the maximal cliques of $\mathcal{G}$.  If $P$ is positive semidefinite globally, then applying \cref{T:ChordalDecompositionTheorem} for each $x \in \mathbb{R}^n$ reveals that there exists $x$-dependent positive semidefinite matrices $S_k:\mathbb{R}^n \to \mathbb{S}_+^{\abs{\mathcal{C}_k}}$ such that
\begin{equation}\label{e:pop-matrix-chordal-decomposition-basic}
P(x) = \sum_{i=1}^t E_{\mathcal{C}_k}^\tr \, S_k(x) \, E_{\mathcal{C}_k}.
\end{equation}
However, this decomposition is not immediately useful in practice because the matrices $S_k$ need not be polynomial, so they cannot be searched for using SOS methods. As an example, consider
\begin{equation*}
P(x) = \begin{pmatrix}
2+x^2 & x+x^2 & 0 \\x+x^2 & 1+2x^2 & x-x^2 \\ 0 & x-x^2 & 2 + x^2
\end{pmatrix},
\end{equation*}
whose sparsity graph is a simple three-node chain graph with two maximal cliques, $\mathcal{C}_1 = \{1,2\}$ and $\mathcal{C}_2 = \{2,3\}$. \cite{zheng2020sum} proved that this matrix is positive definite globally, but does not admit a chordal decomposition~\cref{e:pop-matrix-chordal-decomposition-basic} with polynomial $S_1$ and $S_2$. Using this example, and recalling that all positive semidefinite univariate polynomial matrices are also SOS, one can prove the following general statement.

\begin{proposition}[\citealp{zheng2020sum}]\label{th:sos-matrix-decomp-nonexistence}
	Let $\mathcal{G}$ be a connected and not complete chordal graph with $r\geq 3$ vertices and maximal cliques $\mathcal{C}_1,\ldots,\mathcal{C}_t$. For any positive integers $n$ and $d$, there exists a positive definite SOS matrix $P \in \Sigma_{n,2d}^{r}$ with sparsity graph $\mathcal{G}$ that does not admit a decomposition~\cref{e:pop-matrix-chordal-decomposition-basic} with polynomial matrices $S_1,\ldots,S_t$.
\end{proposition}

On the other hand, the direct proof of \cref{T:ChordalDecompositionTheorem} given by \cite{kakimura2010direct} can be combined with a diagonalization procedure for polynomial matrices due to \cite{schmudgen2009noncommutative} to show that~\cref{e:pop-matrix-chordal-decomposition-basic} holds with SOS matrices $S_1,\ldots,S_k$ for \textit{all} positive semidefinite polynomial matrices, up to multiplication by an SOS polynomial.

\begin{theorem}[\citealp{zheng2020sum}]\label{th:matrix-chordal-sos-decomp-global}
	Let $P \in \mathbb{R}[x]_{n,2d}^{r \times r}$ be positive semidefinite and let $\mathcal{C}_1,\ldots,\mathcal{C}_t$ be the maximal cliques of its sparsity graph. There exist $\nu \in \mathbb{N}$, an SOS polynomial $\sigma \in \Sigma_{n,2\nu}$, and SOS polynomial matrices $S_k \in \Sigma_{n,2d+2\nu}^{\abs{\mathcal{C}_k}}$ for $k=1,\ldots,t$, such that
	\begin{equation}\label{e:weighted-sos-matrix-decomp-global}
	\sigma(x) P(x) = \sum_{k=1}^t E_{\mathcal{C}_k}^\tr \, S_k(x) \, E_{\mathcal{C}_k}.
	\end{equation}
\end{theorem}
\noindent
When the maximal cliques of the sparsity graph of $P$ are small, this result enables one to construct an SOS certificate of global positive semidefinitess using small matrix SOS constraints, which have a much lower computational complexity than simply requiring $\sigma P$ to be SOS. 

%Since the SOS weight $\sigma$ must also be tuned, however, this can only be done when the matrix $P$ is fixed. When it depends on a vector of parameters $\lambda \in \mathbb{R}^\ell$ to be optimized such that $P$ is positive semidefinite, 
Implementation of the chordal SOS decomposition in \cref{th:matrix-chordal-sos-decomp-global} using SDPs requires the matrix $P$ to be fixed, because the SOS weight $\sigma$ must be determined alongside the SOS matrices $S_1,\ldots,S_t$. Often, however, $P$ depends on a vector of parameters $\lambda \in \mathbb{R}^\ell$ that must be optimized whilst ensuring that $P$ is positive semidefinite. In these cases, condition~\cref{e:weighted-sos-matrix-decomp-global} is not jointly convex in $\lambda$ and $\sigma$, so the latter must be fixed \textit{a priori}. This is generally restrictive because, when $\sigma$ is fixed arbitrarily, \cref{th:sos-matrix-decomp-nonexistence} implies that the decomposition~\cref{e:weighted-sos-matrix-decomp-global} may not exist. However, one can prove a sparse-matrix version of Reznick's Positivstellensatz~\citep{Reznick1995} 
to conclude that the weight $\sigma(x) = \|x\|_2^{2\nu}$ is guaranteed to work at least when $P$ is a homogeneous positive definite matrix. % with even support.

\begin{theorem}[\citealp{zheng2020sum}]\label{th:sos-matrix-decomp-even}
	Let $P \in \mathbb{R}[x]_{n,2d}^{r \times r}$ be homogeneous of degree $2d$ and positive definite on $\mathbb{R}^n \setminus \{0\}$. Let $\mathcal{C}_1,\ldots,\mathcal{C}_t$ be the maximal cliques of the sparsity graph of $P$. There exist $\nu \in \mathbb{N}$ and SOS polynomial matrices $S_k \in \Sigma_{n,2d+2\nu}^{\abs{\mathcal{C}_k}}$ for $k=1,\ldots,t$, such that
	$$\|x\|_2^{2\nu} P(x) = \sum_{k=1}^t E_{\mathcal{C}_k}^\tr \, S_k(x) \, E_{\mathcal{C}_k}.$$
\end{theorem}

Decomposition results such as this, where the SOS weight $\sigma$ is fixed, are of considerable interest because they enable the construction of convergent hierarchies of sparsity-exploiting SOS relaxations for optimization problems with global polynomial matrix inequalities (see \citealp{henrion2006convergent,henrion2011inner} and \citealp{peet2009positive} for particular examples). To illustrate the idea, let us consider the generic convex minimization problem
\begin{align}
b^* := \min_{\lambda \in \mathbb{R}^\ell} \quad & b(\lambda)
\nonumber\\[-1ex]
\text{s.t.} \quad & \underbrace{P_0(x) + \sum_{i=1}^\ell \lambda_i P_i(x) }_{=:P(x,\lambda)}\succeq 0 \; \forall x  \in \mathbb{R}^n,
\label{e:opt-problem-maxiz-ineq-global}
\end{align}
where $b: \mathbb{R}^\ell \to \mathbb{R}$ is a convex cost function and $P_0,\ldots,P_\ell \in \mathbb{R}[x]_{n,2d}^{r \times r}$ are symmetric polynomial matrices whose sparsity graph is chordal and has maximal cliques $\mathcal{C}_1,\ldots,\mathcal{C}_t$. Given any integer $\nu\geq 0$, a feasible vector $\lambda$ and an upper bound on the optimal cost $b^*$ may be found by solving the SOS relaxation
\begin{align}
b_\nu^* := \min_{\lambda \in \mathbb{R}^\ell} \; & b(\lambda)
\nonumber\\[-1ex]
\text{s.t.} \; 
&\|x\|_2^{2\nu} P(x;\lambda)  = \sum_{k=1}^t E_{\mathcal{C}_k}^\tr S_k(x) E_{\mathcal{C}_k},
\nonumber\\
&S_k(x) \in \Sigma_{n,2d+2\nu}^{\abs{\mathcal{C}_k}} \text{ for } k=1,\ldots,t,
\end{align}
which can be reformulated as a standard-form SDP.
If the polynomial matrices $P_0,\ldots,P_\ell$ are homogeneous of even degree, %and have even support, 
and there exists $\lambda_0 \in \mathbb{R}^\ell$ such that $P(x,\lambda_0)$ is positive definite, then one can use \cref{th:sos-matrix-decomp-even} to prove that $b^*_\nu \to b^*$ from above as $\nu \to \infty$; see \cite{zheng2020sum} for more details and numerical examples.  Under further technical assumptions (see \citealp{zheng2020sum} for details), asymptotic convergence when $P_0,\ldots,P_\ell$ are not homogeneous is preserved by replacing the SOS multiplier $\|x\|_2^{2\nu}$ with $(1 + \|x\|_2^2)^{\nu}$.

%%%%%%%%%%%%%%%%%%%%%%%%%%%%%%
\subsubsection{Decomposition on a semialgebraic set}\label{ss:local-matrix-sos}
\noindent
We now turn our attention to sparse polynomial matrix inequalities on a semialgebraic set $\mathbb{K}$ defined as in~\cref{e:semialgebraic-set} by $m$ polynomial inequalities $g_i(x) \geq 0$, $i=1,\ldots,m$. A sufficient condition for a symmetric polynomial matrix $P \in \mathbb{R}[x]^{r \times r}_{n,d}$ to be positive semidefinite on $\mathbb{K}$ is that there exist an integer $\nu \in \mathbb{N}$ and SOS matrices $S_0,\ldots,S_m \in \Sigma_{n,2\nu}^m$ such that
\begin{equation}\label{e:matrix-weighted-sos-local}
P(x) = S_0(x) + \sum_{i=1}^m g_i(x) S_i(x).
\end{equation}
A matrix version of \hyperref[th:putinar]{Putinar's Positivstellensatz} proved by \cite{scherer2006matrix} states that this condition~\cref{e:matrix-weighted-sos-local} is also necessary when $P$ is positive definite on $\mathbb{K}$ and this set satisfies the \nameref{ass:archimedean}.

The weighted matrix SOS decomposition~\cref{e:matrix-weighted-sos-local} can be searched for with semidefinite programming, but this is prohibitively expensive when $P$ is large. If it has chordal structural sparsity, however, 
%arguments by \cite{scherer2006matrix} can be combined with the chordal decomposition \cref{T:ChordalDecompositionTheorem} to
one can 
show that the SOS matrices $S_i$ admit a clique-based decomposition. This yields the following sparse matrix version of \hyperref[th:putinar]{Putinar's Positivstellensatz}.

\begin{theorem}[\citealp{zheng2020sum}]\label{th:sos-matrix-decomp-local}
	Let $\mathbb{K}$ be a semialgebraic set defined as in~\cref{e:semialgebraic-set} that satisfies the \nameref{ass:archimedean}. Suppose that the symmetric polynomial matrix $P \in \mathbb{R}[x]_{n,d}^{r \times r}$ is positive definite on $\mathbb{K}$ and that its sparsity graph has maximal cliques $\mathcal{C}_1,\ldots,\mathcal{C}_t$. There exist an integer $\nu \in \mathbb{N}$ and SOS matrices $S_{i,k} \in \Sigma_{n,2\nu}^{\abs{\mathcal{C}_k}}$ for $i=0,\ldots,m$ and $k=1,\ldots,t$ such that
	\begin{equation}\label{e:sparse-matrix-putinar}
	P(x) = \sum_{k=1}^t E_{\mathcal{C}_k}^\tr \bigg( S_{0,k}(x) + \sum_{i=1}^m  g_i(x)  S_{i,k}(x) \bigg) E_{\mathcal{C}_k}.
	\end{equation}
\end{theorem}

%This result can be proved either by combining the arguments by \cite{scherer2006matrix} with the classical chordal decomposition in \cref{T:ChordalDecompositionTheorem}, or by replicating the proof of the latter given by \cite{kakimura2010direct} using the Weierstrass polynomial approximation theorems.

This result can be used to construct sparsity-exploiting SOS relaxations of optimization problems with polynomial matrix inequalities on compact semialgebraic sets that satisfy the Archimedean condition. For example, consider an optimization problem analogous to~\cref{e:opt-problem-maxiz-ineq-global}, where the polynomial matrix inequality is enforced on $\mathbb{K}$ rather than on the full space $\mathbb{R}^n$, and denote its optimal value by $b^*$. If there exists $\lambda_0 \in \mathbb{R}^\ell$ such that the inequality is strict on $\mathbb{K}$ and this set satisfies the Archimedean condition, then the optimal value of the SOS problem
\begin{align*}
\min_{\lambda\in\mathbb{R}^\ell} \quad &b(\lambda)\\
\text{s.t.} \quad& P(x,\lambda) \text{ satisfies~\cref{e:sparse-matrix-putinar}}\\
&S_{i,k} \in \Sigma_{n,2\nu}^{\abs{\mathcal{C}_k}} \text{ for } i=0,\ldots,m \text{ and }k=1,\ldots,t
\end{align*}
converges to $b^*$ from above as $\nu \to \infty$. Interested readers are referred to \cite{zheng2020sum} for more details and computational examples.

%%%%%%%%%%%%%%%%%%%%%%%%%%%%%%%%%
\subsection{Other approaches}\label{s:other-pop-relaxations}
\noindent
The scalability of SOS approaches to polynomial inequalities and polynomial optimization problems can be improved using techniques beyond those described in this section. One example is to replace semidefinite conditions on a large Gram matrix with stronger conditions based on factor-width-$k$ decompositions, which are discussed in \cref{section:factor-wdith-k} below. For the particular case of $k=2$, one obtains  \emph{scaled diagonally dominant SOS} (SDSOS) certificates of nonnegativity \citep{ahmadi2019dsos}. Another approach is to use \textit{bounded-degree SOS conditions} \citep{Lasserre2017bsos}, in which (loosely speaking) one restrict the degree of the monomial basis $x^\mathbb{B}$ used in the Gram matrix representation and handles monomials of higher degree using positivity certificates that can be reformulated as linear programs. Term sparsity can be exploited in these frameworks, too: the relation between correlative sparsity and SDSOS conditions is discussed by \cite{zheng2018sparse}, while \cite{weisser2018sparse} develop sparsity-exploiting bounded-degree SOS hierarchies.

Finally, when working with polynomials that are invariant under groups of symmetry transformations, a large Gram matrix can be replaced with one that has a block-diagonal structure using symmetry reduction techniques \citep{gatermann2004symmetry,lofberg2009pre,Riener2013}. The block-diagonalization based on sign-symmetries, recovered by the TSSOS and CS-TSSOS hierarchies discussed in \cref{s:tssos,s:cs-tssos}, is only one particular example; more sophisticated strategies require using a ``symmetry-adapted'' basis for the space of polynomials in lieu of the monomial basis $x^\mathbb{B}$.

%%%%%%%%%%%%%%%%%%%%%%%%%%%%%%%%%
\subsection{Open-source software implementations}\label{s:pop-software}
\noindent
Many of the sparsity-exploiting techniques for polynomial optimization described in this section are implemented in open-source software. The Newton polytope reduction technique is implemented in almost all parsers for SOS optimization, including \software{SOSTOOLS} \citep{prajna2002introducing}, \software{YALMIP} \citep{lofberg2004yalmip}, \software{GloptiPoly} \citep{henrion2009gloptipoly} and \software{SumOfSquares.jl} \citep{legat2017sos,weisser2019polynomial}. 
Correlative sparsity techniques are implemented in the MATLAB toolboxes \software{SparsePOP} \citep{waki2008algorithm}  and \software{aeroimperial-yalmip} \citep{aeroimperial-yalmip}. %\footnote{\url{https://github.com/aeroimperial-optimization/aeroimperial-yalmip}}. 
The recent Julia package \software{TSSOS} \citep{Magron2021tssos} implements the TSSOS, chordal-TSOS and CS-TSSOS hierarchies. Term sparsity and symmetries in polynomial optimization can also be exploited through \software{SumOfSquares.jl} \citep{legat2017sos,weisser2019polynomial}.
%%%%%%%%%%%%%%%%%%%%%%%%%%%%%%%%%%%%%%%%%%%%%%%%%%%%
\section{Factor-width decomposition}%: connection and difference}
\label{section:factor-width-two}
\noindent
We have seen that the matrix decomposition approach can lead to significant efficiency improvements in the solution of sparse SDPs (cf. \cref{section:sparse-SDPs}) and sparse polynomial optimization problems (cf. \cref{section:polynomial_optimization}). 
We now turn our attention to the problem of testing positive-semidefinitess of matrices that are not necessarily sparse, for which similar matrix decomposition ideas can also be leveraged using approximation methods. This class of methods is known as \emph{factor-width decomposition} \citep{boman2005factor}. We will highlight its connections and differences with the chordal decomposition reviewed above.  

After reviewing some background in~\cref{section:factor-width-background}, %and defining notations in \cref{section:factor-width-notation},  
we discuss how a hierarchy of inner and outer approximations for positive semidefinite matrices can be constructed based on \emph{factor-width-k matrices} in \cref{section:factor-wdith-k}. We then discuss in \cref{section:blockFW2} how this can be extended further, leading to the notion of \emph{block} factor-width-two matrices \citep{zheng2019block}, which aims to strike a balance between numerical computation and approximation quality. Applications in semidefinite and SOS optimization are discussed in \cref{ss:sdp-applications,ss:sos-applications}. 

\subsection{Background} \label{section:factor-width-background}
\noindent
As emphasized in the previous sections, solving large-scale semidefinite programs is at the centre of many problems in control engineering and beyond, and the development of fast and reliable solvers has attracted significant attention recently, mainly focusing on sparsity exploiting and low-rank solution exploiting methods \citep{de2010exploiting,majumdar2020recent}. Some of these methods attempt to solve the problem exactly using, e.g., chordal decomposition (cf. \cref{section:sparse-SDPs,section:polynomial_optimization}) when sparsity is present, but others are trying to provide approximate solutions when these problems are large and dense. This section focuses on the latter case, i.e., the case of dense and large SDPs, and the general idea is still based on a certain matrix decomposition, similar to~\cref{section:sparse-SDPs,section:polynomial_optimization}.  

One basic approach  is to approximate the positive semidefinite cone $\mathbb{S}^n_+$ with the cone of factor-width-$k$ matrices \citep{boman2005factor}, which allows for a certain matrix decomposition discussed in~\Cref{section:factor-wdith-k} below. We will denote the cone of factor-width-$k$ matrices by $\mathcal{FW}^n_k$, where $n$ is the matrix dimension. The case $k = 2$ is of special interest: this is also the case of symmetric \emph{scaled diagonally dominant} matrices, and enforcing $\mathcal{FW}^n_2$ is equivalent to a number of second-order cone constraints, which implies that linear functions can be optimized over $\mathcal{FW}^n_2$ by solving a second-order cone program (SOCP). Compared to SDPs, SOCPs are much more scalable but this approximation is very conservative: the restricted problem may %have a different optimal point or 
even become infeasible. At the same time, attempting an approximation over $\mathcal{FW}^n_3$ will result into an $\mathcal{O}(n^3)$ number of positive semidefinite constraints, which may not strike a good balance between approximation and computational efficiency.  For this reason, most work has focused on the case of factor-width-two matrices and on some closely related extensions \citep{wang2021polyhedral,ahmadi2015sum,ahmadi2017optimization}.  

This notion of factor-width-two matrices was recently extended to the \emph{block} case by \cite{zheng2019block}, who showed that the approximation quality is significantly improved compared to $\mathcal{FW}^n_2$ and remains computationally feasible unlike the approximation using $\mathcal{FW}^n_3$. At the same time, block factor-width-two matrices can form a new hierarchy of approximations using a ``coarsening'' of the decomposition results (cf.~\Cref{definition:partition_order}). %Moreover, reducing the number of partitions can be shown to improve the quality of the approximation. 
An alternative approach that results in an improved approximation is based on the use of \emph{decomposed structured subsets} \citep{miller2019decomposed}.

\subsection{Factor-width-\texorpdfstring{$k$}{k} decompositions}%and connection with chordal decomposition} 
\label{section:factor-wdith-k}
\noindent
We now introduce the concept of \emph{factor-width-$k$ matrices}, originally defined in \cite{boman2005factor}.
\begin{definition}
    The \emph{factor width} of a matrix $X \in \mathbb{S}^n_+$ is the smallest integer $k$ such that there exists a matrix $V$ where $A = VV^\tr$ and each column of $V$ has at most $k$ nonzeros.
\end{definition}
%As mentioned before 
The factor width of $X$ is also the smallest integer $k$ for which $X$ is the sum of positive semidefinite matrices that are non-zero at most on a $k \times k$ principal submatrix:
\begin{equation}\label{eq:fwk}
    Z = \sum_{i=1}^s E_{\mathcal{C}_i}^\tr Z_i E_{\mathcal{C}_i}
\end{equation}
for some matrices $Z_i \in \mathbb{S}^k_+$, where $\mathcal{C}_i$ is a set of $k$ distinct integers from 1 to $n$ and $s = \binom{n}{k}$. We use $\mathcal{FW}_k^n$ to denote the set of $n \times  n$ matrices with factor-width at most $k$. The dual of $\mathcal{FW}^n_k$ with respect to the normal trace inner product is
$$
    (\mathcal{FW}^n_k)^* = \left\{X \in \mathbb{S}^n \mid E_{\mathcal{C}_i}XE_{\mathcal{C}_i}^\tr \in \mathbb{S}^k_+, \forall i = 1, \ldots, s\right\}.
$$
The following hierarchy of inner/outer approximations of $\mathbb{S}^n_+$ follows directly from these definitions:
%\begin{eqnarray*} 
\begin{equation}
\label{eq:fwkinner}
        \begin{aligned}
        \mathcal{FW}_1^n &\subseteq \mathcal{FW}_2^n \subseteq \ldots \subseteq \mathcal{FW}_n^n = \\ &  \mathbb{S}^n_+ = (\mathcal{FW}_n^n)^* \subseteq \ldots \subseteq  (\mathcal{FW}_2^n)^* \subseteq (\mathcal{FW}_1^n)^*.
    \end{aligned}
\end{equation}
%\end{eqnarray*}

The set $\mathcal{FW}_2^n$ is of particular interest because it is equivalent to the set of symmetric scaled diagonally dominant matrices \citep{boman2005factor}. Furthermore, linear optimization over $\mathcal{FW}_2^n$ can be converted into an SOCP, for which efficient algorithms exist. The better scalability of SOCPs compared to SDPs makes inner approximations of positive semidefinite cones based on $\mathcal{FW}_2^n$ very attractive, and form the basis of the SDSOS framework for polynomial optimization proposed by \cite{ahmadi2019dsos}.

{\begin{remark}[Factor-width decomposition vs chordal decomposition]
The decomposition~\cref{eq:fwk} is formally the same as the chordal decomposition in \cref{T:ChordalDecompositionTheorem}, and the two differ only in the choice of ``cliques'' $\mathcal{C}_1,\ldots,\mathcal{C}_s$. For chordal decomposition, they are the maximal cliques of (a chordal extension of) the sparsity graph of $Z$. For factor-width-$k$ decomposition, instead, they are all $\binom{n}{k}$ sets of $k$ distinct indices from $\{1,\ldots,n\}$. These two different choices, however, have considerably different implications: while chordal decomposition is necessary and sufficient for a sparse matrix to be positive semidefinite, factor-width-$k$ decomposition is only sufficient unless $k=n$. The quality of the approximation of positive semidefinite cones by $(\mathcal{FW}^n_k)^*$ was recently investigated by \cite{song2021approximations} and \cite{blekherman2020sparse}.
\markendexample
\end{remark}
}

\subsection{Block factor-width-two decomposition} \label{section:blockFW2}
\noindent
The representation~\cref{eq:fwk} reveals that checking whether a matrix $Z$ belongs to $\mathcal{FW}^n_k$ for any values of $n$ and $k$ is equivalent to an SDP. When $k<n$, this SDP has smaller semidefinite cones than $\mathbb{S}^n_+$, but may be more expensive than checking whether $Z \in \mathbb{S}^n_+$ directly because of the combinatorial number of cones, $\binom{n}{k}$. Setting $k = 2$ does lead to efficiency gains, but the gap between  $\mathcal{FW}^n_2$ and $\mathbb{S}^n_+$ might be unacceptably large in some applications. For this reason, \emph{block factor-width-two matrices} are of interest.

Recall from~\Cref{ss:block-matrices-notation} the notion of a block-partition of a matrix $Z \in \mathbb{S}^n$ subordinate to a partition $\alpha$ of $n$. Recall also the definition of the index matrix $E_{\mathcal{C}_k,\alpha}$ in~\cref{Eq:IndexMatrixBlock}. We here further define 
\begin{subequations}
\begin{align} 
    E_{i,\alpha} &:= \begin{bmatrix} 0 &\ldots & I_{\alpha_i} & \ldots & 0 \end{bmatrix} \in \mathbb{R}^{\alpha_i \times n},\\
    E_{ij,\alpha} &:= \begin{bmatrix} (E_{i,\alpha})^\tr & (E_{j,\alpha})^\tr \end{bmatrix}^\tr \in \mathbb{R}^{(\alpha_i + \alpha_j) \times n}.
    %i \neq j.
\label{eq:blockbasis}
\end{align}
\end{subequations}
The set of block factor-width-two matrices, denoted by $\mathcal{FW}_{\alpha,2}^n$, is defined as follows \citep{zheng2019block}. 

\begin{definition}\label{def:block-FW}
    For any partition $\alpha = \{\alpha_1, \ldots, \alpha_p\}$ of $n$, a symmetric matrix $Z \in \mathbb{S}^n$ belongs to the class $\mathcal{FW}_{\alpha,2}^n$ of block factor-width-two matrices if and only if %we have
    \begin{equation} \label{eq:BlkFW}
        Z = \sum_{i=1}^{p-1}\sum_{j=i+1}^p (E_{ij,\alpha})^\tr X_{ij} E_{ij,\alpha}
    \end{equation}
    for some $X_{ij} \in \mathbb{S}^{\alpha_i+\alpha_j}_+$, where $E_{ij,\alpha}$ is defined in~\cref{eq:blockbasis}. %in~\cref{eq:blockbasis}.
\end{definition}

It is clear that~\cref{eq:BlkFW} is a direct block extension of~\cref{eq:fwk} when $k = 2$. Also, it is not hard to check that $\mathcal{FW}_{\alpha,2}^n$ is a cone. Its dual (with respect to the trace inner product) is characterized by the following proposition.
\begin{proposition}[\cite{zheng2019block}] \label{prop:set-properties}
For any partition $\alpha = \{\alpha_1, \ldots, \alpha_p\}$ of $n$, the dual of $\mathcal{FW}_{\alpha,2}^n$  is
\begin{multline*}
    (\mathcal{FW}_{\alpha,2}^n)^* = \{X \in \mathbb{S}^n \mid \, E_{ij,\alpha} X (E_{ij,\alpha})^\tr \succeq 0, \\ 1\leq i < j \leq p\}.
\end{multline*}
Furthermore, both $\mathcal{FW}_{\alpha,2}^n$ and $(\mathcal{FW}_{\alpha,2}^n)^*$ are proper cones, \emph{i.e.}, they are convex, closed,  solid, and  pointed cones.
\end{proposition}

It should be clear from \Cref{def:block-FW} and \Cref{prop:set-properties} that semidefinite programming can be used to verify whether a matrix belongs to $\mathcal{FW}_{\alpha,2}^n$ or to $(\mathcal{FW}_{\alpha,2}^n)^*$. While a gap between these cones and the positive semidefinite cone $\mathbb{S}^n_+$ remains, the next theorem states that the size of the gap can be reduced by coarsening the partition $\alpha$ (cf.~\Cref{definition:partition_order}), generally at the expense of increasing the computational complexity of the semidefinite representations of $\mathcal{FW}_{\alpha,2}^n$ and $(\mathcal{FW}_{\alpha,2}^n)^*$. This tradeoff between approximation gap and complexity is the main advantage of using block factor-width-two cones.
\begin{theorem}[\cite{zheng2019block}] \label{theo:inclusion}
    Let $\gamma  \sqsubset \beta \sqsubset \alpha$ be partitions of $n$ with $\alpha = \{\alpha_1,\alpha_2\}$, and let $\mathbf{1} = \{1, \ldots, 1\}$ denote the uniform unit partition. Then,
    \begin{multline*}
        \mathcal{FW}^n_2 = \mathcal{FW}_{\mathbf{1},2}^n \subseteq  \mathcal{FW}_{\gamma,2}^n \subseteq  \mathcal{FW}_{\beta,2}^n \subseteq  \mathcal{FW}_{\alpha,2}^n \equiv \mathbb{S}^n_{+}
        \\
        \equiv (\mathcal{FW}_{\alpha,2}^n)^* \subseteq (\mathcal{FW}_{\beta,2}^n)^* \subseteq (\mathcal{FW}_{\gamma,2}^n)^*, \subseteq  (\mathcal{FW}_{\mathbf{1},2}^n)^*.
    \end{multline*}
\end{theorem}

This result does not quantify how well $\mathcal{FW}_{\alpha,2}^n$ and $(\mathcal{FW}_{\alpha,2}^n)^*$ approximate the positive semidefinite cone. Such information is clearly not only of theoretical interest, but also of practical importance, especially for dense positive semidefinite cone that cannot be studied using chordal decomposition. Some progress in this direction was recently made by \cite{zheng2019block}, who leveraged results by \cite{blekherman2020sparse} to show that the normalized distance between either $\mathcal{FW}^n_{\alpha,2}$ or $(\mathcal{FW}^n_{\alpha,2})^*$ and $\mathbb{S}^n_+$ is at most $\frac{p-2}{p}$, where $p$
is the number of blocks in the partition $\alpha$.

\begin{figure}
	\centering
    \includegraphics[scale=0.9]{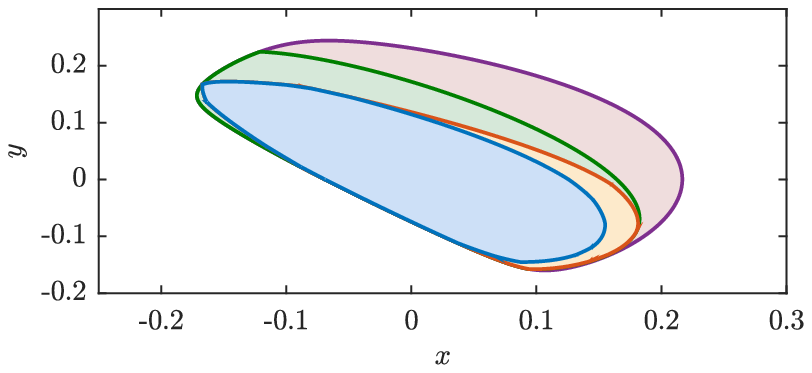}\\[1ex]
    \includegraphics[scale=0.9]{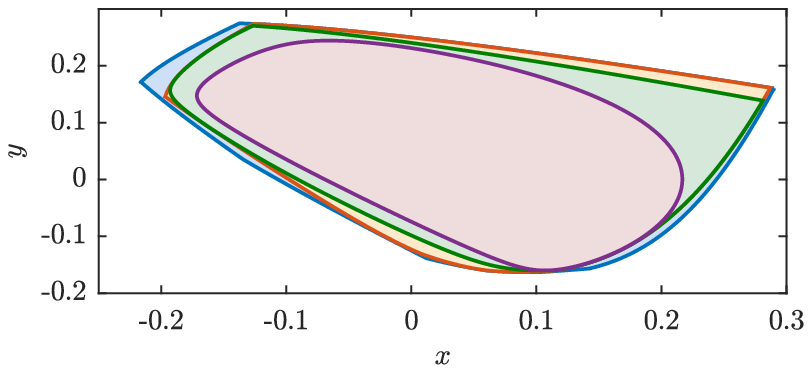}\\[0.25ex]
    \begin{tikzpicture}
    \draw[fill=fwFillBlue, draw=fwLineBlue, thick] (-2.2,0) rectangle ++(0.75,0.35);
    \draw[fill=fwFillRed, draw=fwLineRed, thick] (1.8,0) rectangle ++(0.75,0.35);
    \draw[fill=fwFillGreen, draw=fwLineGreen, thick] (-2.2,-0.5) rectangle ++(0.75,0.35);
    \draw[fill=fwFillMagenta, draw=fwLineMagenta, thick] (1.8,-0.5) rectangle ++(0.75,0.35);
    \node[anchor=west] at (-1.4,0.15) {\footnotesize$\mathcal{FW}^5_{\mathbf{1},2}$, $(\mathcal{FW}^5_{\mathbf{1},2})^*$};
    \node[anchor=west] at (2.6,0.15) {\footnotesize$\mathcal{FW}^5_{\alpha,2}$, $(\mathcal{FW}^5_{\alpha,2})^*$};
    \node[anchor=west] at (-1.4,-0.35) {\footnotesize$\mathcal{FW}^5_{\beta,2}$, $(\mathcal{FW}^5_{\beta,2})^*$};
    \node[anchor=west] at (2.6,-0.35) {\footnotesize$\mathcal{FW}^5_{\gamma,2}$, $(\mathcal{FW}^5_{\gamma,2})^*$};
    \draw[fill=white, draw=white] (-3.0,0) rectangle ++(0.5,0.35);
    \end{tikzpicture}
	\caption{Regions of the $(x,y)$ plane for which the $5 \times 5$ matrix in \cref{ex:fw-example} belongs to the block factor-width-two cones $\mathcal{FW}^5_{\mathbf{1},2} \subset \mathcal{FW}^5_{\alpha,2}\subset\mathcal{FW}^5_{\beta,2}\subset\mathcal{FW}^5_{\gamma,2} \equiv \mathbb{S}^5_+$ (top panel), and the dual cones $\mathbb{S}^5_+ \equiv (\mathcal{FW}^5_{\gamma,2})^* \subset (\mathcal{FW}^5_{\beta,2})^* \subset (\mathcal{FW}^5_{\alpha,2})^*$ (bottom panel). The partitions are $\mathbf{1}=\{1,1,1,1,1\}$, $\alpha = \{2,1,1,1\}$, $\beta =\{2, 1, 2\}$ and $\gamma = \{2, 3\}$. The inclusions of the plotted regions reflect the inclusions of the cones and the order relation $\mathbf{1}\sqsubset \alpha \sqsubset \beta \sqsubset \gamma$.}
	\label{FIG:3}
\end{figure}

Compared to~\cref{eq:fwkinner}, one main advantage of the hierarchy of inner/outer approximations using block factor-width-two cones in \Cref{theo:inclusion} is that the number of basis matrices in the representation~\cref{eq:BlkFW} remains $\mathcal{O}(p^2)$, instead of a combinatorial number $\binom{n}{k}$. Moreover, the value of $p$ decreases when coarsening the partition. Therefore, the cone $\mathcal{FW}^n_{\alpha,2}$ is often computationally more tractable than the cone $\mathcal{FW}^n_k$ with $k \geq 3$. 
\begin{example}\label{ex:fw-example}
Consider the $5\times 5$ matrix
\begin{equation*}
    \small
    \begin{bmatrix}
    1+6x+4y & 3x+y & 2x+y & x+4y & 3x+3y\\
    3x+y & 1+6y & 5x+3y & y & 2x+2y\\
    2x+y & 5x+3y & 1+2x+2y & x+2y & 5x+6y\\
    x+4y & y & x+2y & 1+2x & 3x+3y\\
    3x+3y & 2x+2y & 5x+6y & 3x+3y & 1+6x+2y
    \end{bmatrix}
\end{equation*}
and the progressively coarser partitions $\mathbf{1}=\{1,1,1,1,1\}$, $\alpha = \{2,1,1,1\}$, $\beta =\{2, 1, 2\}$ and $\gamma = \{2, 3\}$.
The regions of the $(x,y)$ plane for which the matrix is in the cones $\mathcal{FW}^5_{\mathbf{1},2} \subset \mathcal{FW}^5_{\alpha,2}\subset\mathcal{FW}^5_{\beta,2}\subset\mathcal{FW}^5_{\gamma,2} \equiv \mathbb{S}^5_+$ are shown in the top panel in \cref{FIG:3}. The bottom panel of the same figure, instead, shows the regions of the plane for which the matrix is in the dual cones
$\mathbb{S}^5_+ \equiv (\mathcal{FW}^5_{\gamma,2})^* \subset (\mathcal{FW}^5_{\beta,2})^* \subset (\mathcal{FW}^5_{\alpha,2})^*$. It is evident from these figures that all of the inclusions are strict. However, the block factor-width-two cones approximate well the positive semidefinite one along some directions.
\markendexample
\end{example}

\subsection{Applications to semidefinite programming}\label{ss:sdp-applications}
\noindent
Recall from  \cref{theo:inclusion} that the cones $\mathcal{FW}^n_{\alpha,2}$ and $(\mathcal{FW}^n_{\alpha,2})^*$ approximate the positive semidefinite cone $\mathbb{S}^n_+$ from the inside and from the outside, respectively, and that the approximation improves as the partition $\alpha$ is coarsened. This allows one to compute convergent sequences of upper and lower bounds on the optimal value of an SDP in the primal standard form~\cref{eq:SDP_primal}, which we denote by $J^*$ for simplicity, using optimization problems of increasing computational complexity that are always simpler to solve than~\cref{eq:SDP_primal} itself. Precisely, since $\mathcal{FW}^n_{\alpha,2} \subseteq \mathbb{S}^n_+$ for any partition $\alpha$ of $n$, the optimal value of the \textit{block factor-width cone program}
\begin{equation} \label{eq:SDPfw}
\begin{aligned}
    U_\alpha := \min_{X} \quad & \langle C, X \rangle \\
    \text{subject to} \quad & \langle A_i, X \rangle  = b_i, \; i = 1, \ldots, m, \\
    &  X \in  \mathcal{FW}^n_{\alpha,2}
\end{aligned}
\end{equation}
bounds the optimal value of the SDP~\cref{eq:SDP_primal} from above. A complementary lower bound is given by
\begin{equation} \label{eq:SDPfw-lower}
\begin{aligned}
    L_\alpha :=\min_{X} \quad & \langle C, X \rangle \\
    \text{subject to} \quad & \langle A_i, X \rangle  = b_i, \; i = 1, \ldots, m, \\
    &  X \in  (\mathcal{FW}^n_{\alpha,2})^*
\end{aligned}
\end{equation}
because $\mathbb{S}^n_+ \subseteq (\mathcal{FW}^n_{\alpha,2})^*$. By \cref{theo:inclusion}, replacing $\alpha$ with a coarser partition can only improve these upper and lower bounds, and we have the following corollary.

\begin{corollary}
        Let $J^*$ denote the optimal value of the SDP~\cref{eq:SDP_primal} and let $\alpha_1 \sqsubseteq \alpha_2 \sqsubseteq \ldots \sqsubseteq \alpha_k = \{\alpha_{k1},\alpha_{k2}\}$ be a sequence of partitions of $n$. Then, $L_{\alpha_1}\leq \cdots \leq L_{\alpha_k}=J^*= U_{\alpha_k}\leq \cdots \leq U_{\alpha_1}$.
\end{corollary}

When $\alpha = \mathbf{1} =\{1, \ldots, 1\}$ is the finest possible partition, problems~\cref{eq:SDPfw,eq:SDPfw-lower} can be reformulated as SOCPs. This case was studied extensively by \cite{ahmadi2019dsos}, and numerical experiments show that the optimal values $L_{\mathbf{1}}$ and $U_{\mathbf{1}}$ can often be very poor bounds for $J^*$. %As shown in~\Cref{theo:inclusion}, we can create a coarser partition $\alpha \sqsubseteq \beta$ that improves the solution accuracy $J^{\alpha} \geq J^{\beta} \geq J^*$.
To obtain better results using coarser partitions, one can leverage the definition of $ \mathcal{FW}^n_{\alpha,2}$ and rewrite the upper bound problem~\cref{eq:SDPfw} as
% \begin{equation} \label{eq:SDPfw_s1}
% \begin{aligned}
%     \min_{X_{jl}} \quad & \sum_{j=1}^{p-1} \sum_{l = j+1}^{p}\left\langle E_{jl,\alpha}C(E_{jl,\alpha})^\tr, X_{jl} \right\rangle \\
%     \text{s.t.} \quad & \sum_{j=1}^{p-1} \sum_{l = j+1}^{p} \left\langle E_{jl,\alpha}A_i(E_{jl,\alpha})^\tr, X_{jl} \right\rangle  = b_i, i = 1, \ldots, m, \\
%     &  X_{jl} \in  \mathbb{S}^{\alpha_j + \alpha_l}_{+},  1 \leq j < l \leq p.
% \end{aligned}
% \end{equation}
\begin{equation} \label{eq:SDPfw_s1}
\begin{aligned}
    \min_{X_{jl}} \quad & \sum_{j=1}^{p-1} \sum_{l = j+1}^{p}\left\langle C_{jl,\alpha}, X_{jl} \right\rangle \\
    \text{s.t.} \quad 
    & \sum_{j=1}^{p-1} \sum_{l = j+1}^{p} \left\langle A_{ijl,\alpha}, X_{jl} \right\rangle  = b_i, \;i = 1, \ldots, m, \\
    &  X_{jl} \in  \mathbb{S}^{\alpha_j + \alpha_l}_{+},  1 \leq j < l \leq p,
\end{aligned}
\end{equation}
where
%\begin{align*}
 $    C_{jl,\alpha} := E_{jl,\alpha} C (E_{jl,\alpha})^\tr,  A_{ijl,\alpha} := E_{jl,\alpha}A_i(E_{jl,\alpha})^\tr.$
%\end{align*}
This is a standard-form SDP and can be solved with general-purpose solvers. Observe that the number of equality constraints in this SDP is the same as for the original problem~\cref{eq:SDP_primal}, but the dimension of semidefinite cones has been reduced. Since general-purpose SDP solvers can handle multiple small semidefinite cones much more efficiently than a single large one, problem~\cref{eq:SDPfw_s1} can often be solved much faster than~\cref{eq:SDP_primal}. For instance, the numerical experiments in \cite{zheng2019block} show that {useful} upper bounds $U_\alpha$ on the optimal value of SDP relaxations of polynomial optimization problems can be found with a reduction of up to 80\% in CPU time. 

%%%%%%%%%%%%%%%%%%%%%%%%%%%%%%%%%%%%%%%%%%%%%
\subsection{Applications to SOS optimization}\label{ss:sos-applications}
\noindent
Block factor-width-two decompositions can also be applied to reduce the computational cost of SOS optimization. 
As discussed in \cref{section:polynomial_optimization}, an $n$-variate polynomial $p \in \mathbb{R}[x]_{n,2d}$ of even degree $2d$ is SOS if and only if there exists an exponent set $\mathbb{B}\subseteq \mathbb{N}^n_d$ and a positive semidefinite matrix $Q$ such that \citep{parrilo2000structured}
\begin{equation} \label{eq:SOSpsd}
    p(x) = (x^\mathbb{B})^\tr \, Q \, x^\mathbb{B}.
\end{equation}
The fundamental computational challenge in optimization over the cone $\Sigma_{n,2d}$ of $n$-variate SOS polynomials of degree at most $2d$ is that the parameterization~\cref{eq:SOSpsd} requires in general an $N \times N$ positive semidefinite matrix with $N=\binom{n+d}{d}$. This may be prohibitive even for moderate values of $n$ and $d$. 

For polynomials characterized by term sparsity, the computational complexity can be reduced dramatically using the approaches reviewed in \cref{section:polynomial_optimization}, which are based on chordal decomposition. To handle polynomials that are not term sparse, \cite{ahmadi2019dsos} introduced the notion of \emph{scaled diagonally dominant sum-of-squares (\SDSOS)}. These are special SOS poynomials whose Gram matrix $Q$ in~\cref{eq:SOSpsd} belongs to the factor-width-two cone $\mathcal{FW}^{\abs{\mathbb{B}}}_2$. As in the case of semidefinite programming, defining block-\SDSOS\ polynomials by replacing $\mathcal{FW}^{\abs{\mathbb{B}}}_2$ with its superset $\mathcal{FW}^{\abs{\mathbb{B}}}_{\alpha,2}$ for any partition $\alpha$ of $\abs{\mathbb{B}}$ offers an improved inner approximation of $\Sigma_{n,2d}$.
\begin{definition}
    Given a partition $\alpha = \{\alpha_1, \ldots, \alpha_g\}$ of $\abs{\mathbb{B}}$, a polynomial $p \in \mathbb{R}[x]_{n,2d}$ is said to be $\alpha$-\SDSOS\ if and only if there exists coefficient vectors $f_{ij,t} \in \mathbb{R}^{\alpha_i + \alpha_j}$ and exponent sets $\mathbb{B}_{ij} \subseteq \mathbb{N}^n_d$ such that
    \begin{equation} \label{eq:sdsos}
        p(x) = \sum_{1 \leq i < j \leq g}\left(\sum_{t=1}^{\alpha_i + \alpha_j} \left(f_{ij,t}^\tr x^{\mathbb{B}_{ij}}\right)^2\right).
    \end{equation}
\end{definition}

The set of all $\alpha$-\SDSOS\ polynomials in $n$ independent variables and degree no larger than $2d$ will be denoted by $\alpha$-$\SDSOS_{n,2d}$. It is not difficult to check that it is a cone. Moreover, since definition~\cref{eq:sdsos} is considerably more structured that the definition~\cref{e:sos-definition} of general SOS polynomials, the inclusion $\alpha$-$\SDSOS_{n,2d} \subseteq \Sigma_{n,2d}$ is immediate.

For the uniform unit partition $\alpha = \{1,\ldots, 1\}$ of $\binom{n+d}{d}$, the cone  $\alpha$-$\SDSOS_{n,2d}$ reduces to the normal \SDSOS\ cone studied by \cite{ahmadi2019dsos}. At the other hand of the spectrum, for any partition in the form $\alpha = \{\alpha_1, \alpha_2\}$ one has  $\alpha$-$\SDSOS_{n,2d} = \Sigma_{n,2d}$. This second statement is a direct consequence of the following result, which reveals a connection between the polynomial cone $\alpha$-$\SDSOS_{n,2d}$ and the block factor-width-two cone $\mathcal{FW}^{\abs{\mathbb{B}}}_{\alpha, 2}$.

\begin{theorem}[\citealp{zheng2019block}] \label{theo:sdsos}
    A polynomial $p \in \mathbb{R}[x]_{n,2d}$ belongs to the cone $\alpha$-$\SDSOS_{n,2d}$ if and only if it admits a Gram matrix representation~\cref{eq:SOSpsd} with $\mathbb{B}\subseteq\mathbb{N}^n_d$ and $Q \in \mathcal{FW}^{\abs{\mathbb{B}}}_{\alpha, 2}$.
\end{theorem}
Similar to~\Cref{theo:inclusion}, we can build a hierarchy of inner approximations for the SOS cone $\Sigma_{n,2d}$. 
%This is stated in the following corollary.
\begin{corollary}\label{coro:sdsosinclusion}
    Let $\mathbf{1} = \{1, \ldots, 1\}$, $\alpha = \{\alpha_1, \ldots, \alpha_g\}$, $\beta = \{\beta_1, \ldots, \beta_h\}$ and $\gamma = \{\gamma_1,\gamma_2\}$ be partitions of $\abs{\mathbb{B}}$ such that $\alpha  \sqsubseteq \beta$. Then,
    \begin{multline}
        \SDSOS_{n,2d} = \mathbf{1}\text{-}\SDSOS_{n,2d} \subseteq  \alpha\text{-}\SDSOS_{n,2d} \\ 
        \subseteq  \beta\text{-}\SDSOS_{n,2d} \subseteq  \gamma\text{-}\SDSOS_{n,2d} = \Sigma_{n, 2d}.
    \end{multline}
\end{corollary}
Consider now an optimization problem of the form 
\begin{equation}
\label{E:generalSOS}
    \begin{aligned}
        w^* := \min_{u}\quad & w^\tr u \\[-1ex]
        \text{s.t.} \quad  & p(x) := p_0(x) + \sum_{i=1}^t u_ip_i(x) \geq 0, \forall x \in \mathbb{R}^n,%\in %\alpha\text{-}\SDSOS_{n,2d},
    \end{aligned}
\end{equation}
where {$p_0,\,\ldots,\,p_t\in \mathbb{R}[x]_{n,2d}$} are given polynomials, $w \in \mathbb{R}^t$ is a given cost vector, and  $u \in \mathbb{R}^t$ is the decision variable. Let $\alpha$ be any partition of $\binom{n+d}{d}$. To compute an upper bound on the optimal cost $w^*$, one can strengthen the nonnegavity constraint on $p$ with the SOS constraints $p \in \Sigma_{n,2d}$, the SDSOS constraint $p \in \SDSOS_{n,2d}$, or the block-SDSOS constraint $p \in \alpha\text{-}\SDSOS_{n,2d}$. The first approach replaces~\cref{E:generalSOS} with an SDP, the second one leads to an SOCP, and the third yields a block-factor-width cone program that can be reformulated as a standard-form SDP. According to \cref{coro:sdsosinclusion}, the SOS constraint provides the best upper bound on $w^*$, but is the most computationally expensive. At the other extreme is the SDSOS constraint, which offers the fastest computations but may be too restrictive---in fact, the corresponding SOCP may even be infeasible. The  block-SDSOS constraint $p \in \alpha\text{-}\SDSOS_{n,2d}$, instead, can balance the computational speed and upper bound quality thanks to the freedom one has in choosing the partition $\alpha$. This expectation is confirmed by the numerical experiments of \cite{zheng2019block}, but the problem of choosing an optimal partition for given computational resources remains an open problem. 

%%%%%%%%%%%%%%%%%%%%%%%%%%%%%%%%%%%%%%%%%%%%%%%%%%%%

\begin{table*}[t]
     \setlength{\abovecaptionskip}{-2mm}
     \setlength{\belowcaptionskip}{0mm}
     \renewcommand\arraystretch{1.0}
     \caption{Applications of exploiting chordal sparsity in control, machine learning, relaxation of QCQP (Quadratically-constrained quadratic program), fluid dynamics, and beyond.} 
     \begin{center}
     \small
  \begin{tabular}{l l l}%{c C{15mm} C{20mm} C{24mm} C{18mm} c }
  \toprule
   \multicolumn{1}{l}{Area}   &  Topic  &  References \\ %&  Key feature  \\
    \hline
    Control & Linear system analysis &  \multicolumn{1}{p{9.4cm}}{\cite{mason2014chordal,andersen2014robust,pakazad2017distributed,ZKSP2018scalable,deroo2015distributed}}  \\
            & Decentralized control & \multicolumn{1}{p{9.2cm}}{\cite{ZMP2018Scalable,ZKSP2020Distributed,heinke2020distributed,deroo2014distributed}}  \\
            & Nonlinear system analysis & \multicolumn{1}{p{9.2cm}}{\cite{Tacchi2019,Schlosser2020,zheng2018sparse}; \citet[Chapter 5]{mason2015chordal}} \\
            & Model predictive control & \cite{hansson2018exploiting,ahmadi2019efficient} \\[1ex]
    Machine learning & Verification of neural networks &  \multicolumn{1}{p{9.2cm}}{\cite{BPLZ2021neural,zhang2020tightness,dvijotham2020efficient,NewP21}}  \\
    & Lipschitz constant estimation &\cite{latorre2020lipschitz,chen2020semialgebraic} \\ 
                     & Training of support vector machine & \cite{andersen2010support} \\
                     & Geometric perception \& coarsening &\cite{yang2020one,chen2020chordal,liu2019spectral} \\
                     & Covariance selection & \cite{dahl2008covariance,zhang2018large} \\ 
                     & Subspace clustering &\cite{miller2019chordal} \\ [1ex]
    \multirow{2}{2.5cm}{Relaxation of QCQP and POPs} & Sensor network locations & \multicolumn{1}{p{9cm}}{\cite{kim2009exploiting,nie2009sum,jing2019angle}} \\
    & Max-Cut problem &  \cite{ZFPGW2020chordal,andersen2010implementation,garstka2019cosmo}\\
    & Optimal power flow (OPF) & \multicolumn{1}{p{9cm}}{\cite{andersen2014reduced,dall2013distributed,jabr2011exploiting,molzahn2013implementation,molzahn2014sparsity,jiang2017minimum}} \\
    & State estimation in power systems & \cite{zhu2014power,zhang2017conic,weng2013distributed} \\[1ex]
    Others & 
               Fluid dynamics &\cite{Fantuzzi2018beanard-marangoni,Arslan2021}\\
& Partial differential equations & \multicolumn{1}{p{9cm}}{\cite{Mevissen2008solving-pdes,Mevissen2009cavity-flow,Mevissen2011smooth-pde,mevissen2010sparse}} \\
           & Robust quadratic optimization & \cite{andersen2010linear} \\
           & Binary signal recovery & \cite{fosson2019recovery} \\
           & Solving polynomial systems &  \multicolumn{1}{p{9cm}}{\cite{mou2021chordal,li2021choosing,cifuentes2016exploiting,cifuentes2017chordal,Tacchi2019sparse-sets}}\\
           & Other problems  &\multicolumn{1}{p{9cm}}{\cite{yang2020exploiting,jeyakumar2016semidefinite,khoshfetrat2017distributed,baltean2019scoring,madani2017finding}} \\
    \bottomrule
  \end{tabular}
  \end{center}
\label{tab:applications}
\end{table*}

\section{Applications}
%in control, machine learning, fluid dynamics, and beyond}
\label{section:applications}
\noindent
The matrix decomposition techniques reviewed in the previous sections can be used to reduce the computational complexity of a wide variety of analysis and control problems that can be formulated as SDPs or SOS programs.  
%In this section, we present some applications from a wide variety of areas, which exploit chordal sparsity to improve computational efficiency and/or design distributed solutions. 
As anticipated in~\Cref{subsection:aggregate_sparsity_pattern}, complex large-scale dynamical systems at the heart of modern technology often possess a natural graph-like structure, due for example to sparse interactions between subsystems in a network \citep{riverso2014plug,zheng2020smoothing,ZMP2018Scalable,andersen2014reduced,dall2013distributed}. %\red{(Add references?)}. 
The key to enabling efficient numerical treatment of control problems for such systems is to devise SDP or SOS relaxations that preserve this graph structure as much as possible. Precisely, one aims to obtain SDPs with aggregate sparsity (cf. \cref{section:sparse-SDPs}) or polynomial optimization problems with term sparsity (cf. \cref{section:polynomial_optimization}). If this can be done, %chordal decomposition in \cref{T:ChordalCompletionTheorem,T:ChordalDecompositionTheorem,T:GeneralCompletionTheorem,T:GeneralDecompositionTheorem} and 
then the sparsity exploiting techniques discussed in \cref{section:sparse-SDPs,section:polynomial_optimization} can bring considerable computational gains and enable the study of very large systems.

%\cref{tab:applications} lists some representative applications.~
This section describes how chordal sparsity can be exploited for a small selection of problems in control and machine learning. \Cref{subsection:control_applications} focuses on stability analysis for linear and nonlinear systems, and on decentralized control of networked linear systems. In \cref{subsection:QCQP}, we review sparsity-promoting relaxations of nonconvex quadratically constrained quadratic programs (QCQPs) and apply them to the well-known Max-Cut problem from graph theory, as well as to a network sensor location problem. Finally, \Cref{subsection:verification} shows how chordal sparsity allows for efficient verification of neural networks in machine learning. %, and fluid dynamics. 
We stress that these are only a few of the application domains in which chordal decomposition has enabled considerable progress in recent years; other fields include, for instance, fluid mechanics, model predictive control, and optimal power flow. \Cref{tab:applications} provides a (non-exhaustive) list of references.

\subsection{Stability analysis and decentralized control} \label{subsection:control_applications}
\noindent
Stability analysis and control synthesis problems for dynamical systems governed by ordinary differential equations can often be reformulated as SDPs or SOS programs using Lyapunov functions \citep{boyd1994linear,zhou1996robust,papachristodoulou2005tutorial, parrilo2000structured, lasserre2010moments}. If the interactions between individual components of the system have a sparse graph structure, considering Lyapunov functions with a separable or nearly-separable structure can lead to sparse SDPs and SOS programs, which can be solved efficiently using the techniques in \cref{section:sparse-SDPs,section:polynomial_optimization}. Here, we give three simple examples of this fact.

\subsubsection{Stability of linear networked systems}  \label{subsection:stability_analysis}
\noindent
Consider a continuous-time linear autonomous system
\begin{equation} \label{eq:LTI_autonomous}
    \dot{x}(t) = Ax(t),
\end{equation}
where $x (t)\in \mathbb{R}^n$ is the system state at time $t$ and $A \in \mathbb{R}^{n\times n}$ is the system matrix. It is well known \citep{boyd1994linear,zhou1996robust} that the equilibrium state $x(t)=0$ is asymptotically stable if and only if all eigenvalues of $A$ have negative real part. Classical Lyapunov stability theory guarantees that this is true if and only if there exists a positive definite matrix $P$ such that the (positive definite) Lyapunov function $V(x) = x^\tr P x$ decays monotonically along all system trajectories $x(t)$. Equivalently, $P$ must satisfy the strict LMIs
\begin{equation} \label{eq:LyapunovLMI}
    P \succ 0, \qquad A^\tr P + PA \prec 0.
\end{equation}

Now, suppose that~\cref{eq:LTI_autonomous} is a compact representation of a network of $l$ linear subsystems with states $x_1 \in \mathbb{R}^{n_1}$, $\ldots$, $x_l \in \mathbb{R}^{n_l}$, whose interactions can be represented by a static undirected graph $\mathcal{G}_d(\{1,\ldots,l\}, \mathcal{E}_d)$ with $(i,j) \in \mathcal{E}_d$ if and only if systems $i$ and $j$ are directly coupled. In particular, the dynamics of each subsystem are given explicitly by
\begin{equation}
    \label{eq:network-linear-system}
    \dot{x}_i = A_{ii}x_i + \sum_{j \in \mathcal{N}_i} A_{ij} x_j, \qquad i = 1, \ldots, l,
\end{equation}
where $\mathcal{N}_i:=\{j:\, (j,i) \in \mathcal{E}_d\}$ denotes the neighbors of system $i$. Systems of this type are encountered, for example, when modelling power grids \citep{riverso2014plug} and traffic systems \citep{zheng2020smoothing,wang2020leading}.

If the matrix $P$ in~\cref{eq:LyapunovLMI} is assumed to be block-diagonal with $l$ blocks of size $n_1,\ldots,n_l$, meaning that we consider a quadratic Lyapunov function in the separable form \citep{boyd1989structured, ZMP2018Scalable,ZKSP2020Distributed,geromel1994decentralized}
\begin{equation} \label{eq:Separable_Lyapunov}
    V(x) = \sum_{i=1}^l x_i^\tr P_i x_i, %= x^\tr \mathrm{diag}(P_1, \ldots, P_l) x,
\end{equation}
then it is not hard to see that the block-sparsity graph of the matrix $A^\tr P + PA$ in~\cref{eq:LyapunovLMI} is the same as the system graph $\mathcal{G}_d$. When this graph is chordal with small maximal cliques, or admits a chordal extension with the same property, a feasible block-diagonal matrix $P$ satisfying~\cref{eq:LyapunovLMI} can be constructed for significantly larger networks than that can be handled without sparsity exploitation. Equivalently, for a given network size, CPU time requirements can be reduced dramatically.

\begin{figure}
    \centering
    \subfigure[]{
    \begin{tikzpicture}
        \node[draw, circle] (1) at (0,0) {\footnotesize 1};
        \node[draw, circle] (2) at (0:1.8cm) {\footnotesize 2};
        \node[draw, circle] (3) at (-{90/4}:1.8cm) {\footnotesize 3};
        \node[draw, circle] (4) at (-{2*90/4}:1.8cm) {\footnotesize 4};
        \node (a) at (-{2.75*90/4}:1.8cm) {\tiny $\bullet$};
        \node (b) at (-{3.0*90/4}:1.8cm) {\tiny $\bullet$};
        \node (c) at (-{3.25*90/4}:1.8cm) {\tiny $\bullet$};
        \node[draw, circle] (l) at (-90:1.8cm) {\footnotesize $l$};
        \draw (1)--(2);
        \draw (1)--(3);
        \draw (1)--(4);
        \draw (1)--(l);
    \end{tikzpicture}
    \label{f:LyapunovLinearStability_a}
    }
    \begingroup % keep the change local
    \setlength\arraycolsep{1pt}
    \def\arraystretch{0.6}
    \subfigure[]{
    \begin{tikzpicture}
        \node at (0,00) {
        \large
        $\left[\begin{array}{ccc c c}\\[-1ex]
        \blacksquare & \blacksquare & \blacksquare & \cdots & \blacksquare \\
        \blacksquare & \blacksquare & & & \\
        \blacksquare &  & \blacksquare & & \\[-1ex]
        \vdots &  & & \ddots& \\
        \blacksquare & & & &\blacksquare \\
        \end{array}\right]
    $};
    \end{tikzpicture}
    \label{f:LyapunovLinearStability_b}
    }
    \endgroup
    \caption{\textit{(a)} Graph $\mathcal{G}_d$ for the network of 10-dimensional linear systems used to generate the results reported in \cref{Table:LyapunovLinearStability}. \textit{(b)} Block sparsity pattern of the matrix $PA + A^\tr P$ when $P$ is block-diagonal; each block has size $10\times 10$, and there are $l$ diagonal blocks.}
    \label{fig:my_label}
\end{figure}
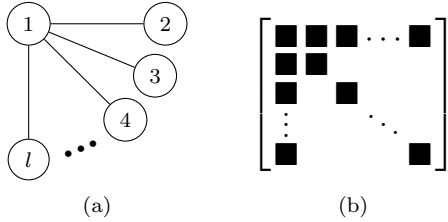

As an example, consider a network with a master node and $l-1$ independent subsystems connected to it, sketched in \cref{f:LyapunovLinearStability_a}. For simplicity, suppose that the subsystems have size $n_1=\cdots=n_l=10$. With a block-diagonal $P$, the second LMI in~\cref{eq:LyapunovLMI} has the chordal ``arrow-type'' block sparsity shown in \cref{f:LyapunovLinearStability_b}. \Cref{Table:LyapunovLinearStability} reports the CPU time required to construct a feasible $P$ with \software{MOSEK} as a function of the number $l$ of subsystems when the sparsity of this LMI is and is not exploited.\footnote{Computations were performed using the MATLAB toolboxes \software{YALMIP} and \software{SparseCoLO} on a laptop with 16GB RAM and an Intel i7 processor. The nonzero system matrices $A_{ij}$ were generated randomly whilst ensuring the existence of a feasible block-diagonal~$P$.} It is evident that exploiting chordal sparsity using the methods described in \cref{section:sparse-SDPs} leads to a significant reduction in CPU time. Similar results are obtained for systems with more realistic network graphs if its maximal cliques are small; see \cite{mason2014chordal}, \cite{deroo2015distributed} and \cite{ZKSP2018scalable,ZMP2018Scalable}.

\begin{table}
    \centering
    \caption{CPU time (in seconds) required by the SDP solver \software{MOSEK} to construct a block-diagonal Lyapunov matrix $P$ satisfying the LMIs in~\cref{eq:LyapunovLMI} for a network of $l$ 10-dimensional systems with connectivity graph $\mathcal{G}_d$ shown in \cref{f:LyapunovLinearStability_a}.}
    \label{Table:LyapunovLinearStability}
    \begin{tabular}{c c c}
        \toprule
        $l$ & No sparsity exploitation & Sparsity exploitation \\ \hline
        10  &      0.55  &      0.26 \\
        50  &     14.92  &      0.90 \\
        100  &     86.09  &      1.21 \\
        125  &    113.06  &      1.17 \\
        150  &    185.42  &      1.96 \\
        175  &    334.13  &      2.69 \\
        200  &    498.49  &      3.55 \\
         \bottomrule
    \end{tabular}
\end{table}

\begin{remark}[Separable Lyapunov functions]
    Searching for a Lyapunov function $V(x)$ with the separable structure~\cref{eq:Separable_Lyapunov} is convenient to ensure that the sparsity of the system matrix $A$ is inherited by the LMI $A^\tr P + PA \prec 0$. The existence of such a separable Lyapunov function can be guaranteed for special classes of stable linear systems \citep{carlson1992block,sootla2017block,sootla2019existence}, but not in general. When a separable Lyapunov function $V(x)$ fails to exist, the structure of the network graph $\mathcal{G}_d$ may be still be leveraged to promote sparsity in~\cref{eq:LyapunovLMI}; for instance, the case of banded graphs, cycles and trees was studied by \cite{mason2014chordal}. Determining a suitable structure for $V(x)$ (equivalently, for the matrix $P$) \textit{a priori} for general graph structures, however, remains a challenging problem. \markendexample
\end{remark}

%%%%%%%%%%%%%%%%%%%%%%%%%%%%%%%%%%%%%%%%%%%%%%%%%%%%%%
\subsubsection{Stability of sparse polynomial systems}
\noindent
Structured Lyapunov functions can bring computational advantages also when studying the asymptotic stability of sparse nonlinear systems with polynomial dynamics. As an example, consider a nonlinear system with the structure \cite[Section VI.D]{zheng2018sparse}
\begin{equation} \label{eq:nonlinear_systems}
\begin{aligned}
    \dot{x}_1 &= f_1(x_{1},x_2), \\
    \dot{x}_i &= f_i(x_{i-1},x_i,x_{i+1}), \qquad i=2,\ldots,l-1, \\
    \dot{x}_l &= f_l(x_{l-1},x_l),
\end{aligned}
\end{equation}
where each vector field $f_i$ depends polynomially on its arguments and $x_i \in \mathbb{R}^{n_i}$. Let $x=(x_1,\ldots,x_l)$ be the collection of all system states and write $f=(f_1,\ldots,f_l)$. Suppose the system has an equilibrium at the origin. This equilibrium is locally asymptotically stable if there exist a region $\mathcal{D} \subset \mathbb{R}^{n_1} \times \cdots \mathbb{R}^{n_l}$ containing the origin, a constant $\epsilon > 0$, and a Lyapunov function $V:\mathbb{R}^{n_1} \times \cdots \times \mathbb{R}^{n_l} \to \mathbb{R}$ such that
\begin{subequations}
\begin{align}
    V(0)&=0,\\
    \label{e:poly-stability-b}
    V(x)&\geq \|x\|_2^2 &&\forall x \in \mathcal{D},\\
    \label{e:poly-stability-c}
    -f(x) \cdot \nabla V(x) &\geq \epsilon \|x\|_2^2 &&\forall x \in \mathcal{D}.
\end{align}
\end{subequations}
Upon fixing $\mathcal{D} = \{x: r_i^2 - \|x_i\|^2 \geq 0 \; \forall i = 1,\ldots,l \}$, which has a fully separable structure, and requiring $V$ to be a polynomial, the last two inequalities become polynomial inequalities on a basic semialgebraic set. One can therefore search for $V$ using SOS optimization. Moreover, the structure of $V$ can be chosen to ensure that these polynomial inequalities are correlatively sparse (cf. \cref{s:sos-csp-unconstrained,s:pop-sparse-constrained}), enabling efficient implementation.

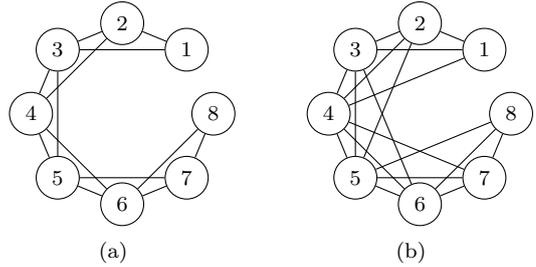
\begin{figure}
    \centering
    \subfigure[]{
    \label{fig:nonlinear-stability-graphs-a}
    \begin{tikzpicture}
        \foreach \a in {1,2,...,8}{
        \node[draw, circle] (\a) at ({\a*360/8}:1.2cm) {\footnotesize \a};
        }
        \draw (1)--(2);
        \draw (2)--(3);
        \draw (3)--(4);
        \draw (4)--(5);
        \draw (5)--(6);
        \draw (6)--(7);
        \draw (7)--(8);
        \draw (1)--(3);
        \draw (2)--(4);
        \draw (3)--(5);
        \draw (4)--(6);
        \draw (5)--(7);
        \draw (6)--(8);
        \end{tikzpicture}
    }
    \hspace{1em}
    \subfigure[]{
    \label{fig:nonlinear-stability-graphs-b}
    \begin{tikzpicture}
        \foreach \a in {1,2,...,8}{
        \node[draw, circle] (\a) at ({\a*360/8}:1.2cm) {\footnotesize \a};
        }
        \draw (1)--(2);
        \draw (2)--(3);
        \draw (3)--(4);
        \draw (4)--(5);
        \draw (5)--(6);
        \draw (6)--(7);
        \draw (7)--(8);
        \draw (1)--(3);
        \draw (2)--(4);
        \draw (3)--(5);
        \draw (4)--(6);
        \draw (5)--(7);
        \draw (6)--(8);
        \draw (1)--(4);
        \draw (2)--(5);
        \draw (3)--(6);
        \draw (4)--(7);
        \draw (5)--(8);
        \end{tikzpicture}
    }
    \caption{Correlative sparsity graphs for the polynomial inequality in~\cref{e:poly-stability-c} when the Lyapunov function $V$ has \textit{(a)} the separable form~\cref{e:poly-nnl-stability-a}, and \textit{(b)} the partially separable form~\cref{e:poly-nnl-stability-b}.}
\end{figure}

For example, if one takes
\begin{equation}\label{e:poly-nnl-stability-a}
    V(x) = \sum_{i=1}^l V_i(x_i)
\end{equation}
to have a fully separable structure as in the case of linear systems considered previously, then the correlative sparsity graph of inequalities~\cref{e:poly-stability-b} is a graph with no edges, while that of~\cref{e:poly-stability-c} is the same chain graph characterizing the cascaded interactions between the state vectors $x_1,\ldots,x_l$, shown in \cref{fig:nonlinear-stability-graphs-a} for $l=8$. If this choice for $V$ is insufficient, one can try the structured choice
\begin{equation}\label{e:poly-nnl-stability-b}
    V(x) = \sum_{i=2}^{l-1}V_i(x_{i-1}, x_i, x_{i+1}).
\end{equation}
In this case, the correlative sparsity graph of~\cref{e:poly-stability-b} is the chain graph mentioned above, while that of~\cref{e:poly-stability-c} is a chordal graph with maximal cliques $\{i,i+1,i+2,i+3\}$ for $i=1,\ldots,l-3$, which is shown in \cref{fig:nonlinear-stability-graphs-b} for $l=8$. One can of course build an entire hierarchy of structured Lyapunov functions with increasing degree of couplings between subsystem variables, at the expense of increasing the number of edges in the correlative sparsity graph of the polynomial inequalities~\cref{e:poly-stability-b,e:poly-stability-c}. Numerical experiments by~\cite{zheng2018sparse} for the structured Lyapunov function in~\cref{e:poly-nnl-stability-b}, which we report in~\cref{Table:LyapunovTime}, show that this approach can significantly reduce the computation time and resources required to prove stability of nonlinear systems compared to standard SOS techniques. 

Similar ideas can be used to partition nonlinear systems into subsystems \citep{anderson2011decomposition} and can be adapted to problems beyond stability analysis, such as the estimation of region of attractions, positively invariant sets, and global attractors \citep{Tacchi2019, Schlosser2020}. 

\begin{table}[t]
	\centering
	\caption{{CPU time, in seconds, required by \software{MOSEK} to construct a structured quadratic Lyapunov function~\cref{e:poly-nnl-stability-b} for a locally asymptotically stable, degree-3 polynomial system of the form~\cref{eq:nonlinear_systems}. Entries marked {\sc oom} indicate memory errors.}}
	\label{Table:LyapunovTime}
	\begin{tabular*}{\linewidth}{@{\extracolsep{\fill}}r| rrrrrrrrr}
		\toprule[1pt]
		$l$ & 10 & 15 & 20 & 30 & 40  & 50 \\
		\midrule% \\
		Standard SOS  &1.4 & 21.3 & 262.1 & {\sc oom} & {\sc oom} & {\sc oom}\\
		Sparse SOS &  0.6 & 0.7  &  0.8 & 1.0 & 1.2 & 1.4  \\
		\bottomrule[1pt]
	\end{tabular*}
\end{table}

\subsubsection{Decentralized control of linear networked systems} 
\noindent
Consider a network of linear system with control inputs and disturbances,
$$
    \dot{x}_i = A_{ii}x_i + \sum_{j \in \mathcal{N}_i} A_{ij} x_j + B_i u_i + M_id_i, \quad i = 1, \ldots, l,
$$
where $x_i \in \mathbb{R}^{n_i}$, $u_i \in \mathbb{R}^{m_i}$ and $d_i \in \mathbb{R}^{q_i}$ denote the local state, input, and disturbance of subsystem $i$, respectively, and $\mathcal{N}_i$ is the index set of all systems connected to system $i$. % (cf.~\cref{eq:network-linear-system}). 
Setting $x = (x_1,\ldots,x_l)$, $u = (u_1,\ldots,u_l)$ and $d = (d_1,\ldots,d_l)$, the system can be written compactly as
\begin{equation*}
    \dot{x}(t) = Ax(t) + Bu(t) + Md(t),
\end{equation*}
where $A$ has block sparsity induced by the system graph (cf. \cref{subsection:stability_analysis}), while $B = \text{diag}(B_1,\ldots,B_l)$ and $M = \text{diag}(M_1,\ldots,M_l)$ are block-diagonal. 

The \emph{optimal decentralized control} problem \citep{geromel1994decentralized} seeks to design static state feedback laws,
\begin{equation} \label{E:DeController}
    u_i(t) = -K_{ii} x_i(t), \qquad \forall i = 1, \ldots, l,
\end{equation}
that minimize the $\mathcal{H}_2$ norm of the transfer function from disturbance $d$ to the output
$$
    z = \begin{bmatrix} Q^{\frac{1}{2}} \\0 \end{bmatrix} x + \begin{bmatrix} 0\\ R^{\frac{1}{2}} \end{bmatrix} u,
$$
where $Q := \text{diag}(Q_1,\ldots, Q_l)$ and $R := \text{diag}(R_1,\ldots, R_l)$ are given block-diagonal matrices. % and diagonal block $Q_i, R_i$ correspond to the subsystem $i$. 
The decentralized constraint~\cref{E:DeController} makes the control problem challenging to solve \citep{geromel1994decentralized, furieri2019separable}. One simple strategy is to enforce %assume 
that the closed-loop system admits a separable Lyapunov function in the form~\cref{eq:Separable_Lyapunov}. This allows translating the decentralized constraint on the controller to other auxiliary design variables
%, known as the notion of \emph{sparsity invariance}~
\citep{furieri2020sparsity,furieri2019separable}. In particular, a suboptimal decentralized controller can be computed using the formula $K_{ii} = Z_iX_i^{-1}$ for each $i = 1, \ldots, l$ (\citealp{geromel1994decentralized}; \citealp[Section II.B]{ZKSP2020Distributed}), where the matrices $Z_1,\ldots,Z_l$ and $X_1,\ldots,X_l$ solve the SDP
\begin{subequations}\label{E:LMIHnorm3}
  \begin{align}
    \min_{X_i,Y_i,Z_i} \;\; & \sum_{i=1}^l \langle Q_i,X_i \rangle+ \langle R_i,Y_i \rangle \nonumber \\
    \text{s.t.} \;\; & (AX\!-\!BZ)+(AX\!-\!BZ)^{\tr} \!+ \! MM^{\tr} \!\preceq \!0, \label{E:LMIHnorm3_eq1}\\
    & \begin{bmatrix} Y_i & Z_i \\ Z_i^{\tr} & X_i \end{bmatrix} \succeq 0, \; X_i \succ 0\quad \forall i = 1, \ldots, l \label{E:LMIHnorm3_eq2}
  \end{align}
\end{subequations}
and $X = {\rm diag}(X_1,\ldots,X_l)$ and $Z = {\rm diag}(Z_1,\ldots,Z_l)$ are block-diagonal concatenations of the matrix variables.

%This problem is an SDP and, in principle, can be solved using general-purpose solvers. Moreover, 
The cost function of this SDP and the constraints in~\cref{E:LMIHnorm3_eq2} are fully separable, as they depend only on variables corresponding to a single subsystem. The coupling constraint~\cref{E:LMIHnorm3_eq1}, instead, has a block sparsity pattern induced by the system graph by virtue of the block-diagonal structure of $B$, $M$, $X$ and $Z$. As in \cref{subsection:stability_analysis}, therefore, the chordal decomposition techniques of~\cref{section:sparse-SDPs} allow for a fast numerical solution when the underlying system graph is sparse, which enables control synthesis for large-scale but sparse networks. In addition, customized distributed design methods that combine chordal decomposition with ADMM can solve~\cref{E:LMIHnorm3} in a privacy-safe way, without requiring subsystems to share information about their local dynamics \citep{ZKSP2020Distributed}.

\subsection{Relaxation of nonconvex QCQPs} \label{subsection:QCQP}
\noindent
A (nonconvex) quadratically constrained quadratic program (QCQP) is an optimization problem in the form 
\begin{equation} \label{eq:QCQP}
    \begin{aligned}
    \min_{x} &\quad x^\tr P_0 x + 2q^\tr_0 x + r_0 \\
    \text{subject to}&\quad x^\tr P_i x + 2q^\tr_i x + r_i \leq 0, \; i = 1,\ldots, m, 
    \end{aligned}
\end{equation}
where $x \in \mathbb{R}^n $ is the optimization variable, and $P_i \in \mathbb{S}^{n}, q_i \in \mathbb{R}, r_i\in \mathbb{R}, i = 0, 1, \ldots, m$ are given problem data. QCQPs have very powerful modeling capabilities; for instance, many hard combinatorial and discrete optimization problems can written in the form~\cref{eq:QCQP} \citep{nesterov2000semidefinite}. This also means that QCQPs are hard to solve in general, so many different relaxation strategies have been proposed to find approximate bounds and feasible values for the optimization variable  $x$ \citep{nesterov2000semidefinite,park2017general}. One approach that provides good bounds, both empirically and theoretically \citep{nesterov2000semidefinite}, is to introduce the positive semidefinite matrix $X = xx^\tr$ and rewrite~\cref{eq:QCQP} as
\begin{equation*} %\label{eq:QCQP_matrix}
    \begin{aligned}
    \min_{x, X} &\quad \langle P_0, X \rangle + 2q^\tr_0 x + r_0 \\
    \text{subject to}&\quad \langle P_i, X \rangle + 2q^\tr_i x + r_i \leq 0, \; i = 1,\ldots, m, \\
    &\quad X = xx^\tr.
    \end{aligned}
\end{equation*}

Upon relaxing the intractable constraint $X = xx^\tr$ into the inequality $X \succeq  xx^\tr $ and applying Schur's complement to rewrite the latter as an LMI, we arrive at the semidefinite relaxation
\begin{equation*} %\label{eq:QCQP_sdp}
    \begin{aligned}
    \min_{x, X} &\quad \langle P_0, X \rangle + 2q^\tr_0 x + r_0 \\
    \text{subject to}&\quad \langle P_i, X \rangle + 2q^\tr_i x + r_i \leq 0, \; i = 1,\ldots, m, \\
    &\quad \begin{bmatrix} 1 & x^\tr \\ x & X \end{bmatrix} \succeq 0,
    \end{aligned}
\end{equation*}
which is equivalent to the following primal-form SDP with nonnegative variables%inequality-constrained SDP
\begin{equation} \label{eq:QCQP_sdp}
    \begin{aligned}
    \min_{Z \in \mathbb{S}^{n+1},w} &\quad \left\langle \left(\begin{smallmatrix}r_0& q_0^\tr\\ q_0 & P_0\end{smallmatrix}\right), Z \right\rangle\\
    \text{s.t.}&\quad \left\langle \left(\begin{smallmatrix}r_i& q_i^\tr\\ q_i & P_i\end{smallmatrix}\right), Z \right\rangle + w_i = 0, \; i = 1,\ldots, m, \\
    &\quad Z_{11} = 1,\\
    &\quad Z\succeq 0, w \geq 0.
    \end{aligned}
\end{equation}
It is not difficult to see that the optimal value of problem~\cref{eq:QCQP_sdp} bounds that of the QCQP~\cref{eq:QCQP} from  below and that, if an optimal solution $Z_\star$ has rank one, then the relaxation is exact and $Z_\star = \left( \begin{smallmatrix} 1 & x_\star^\tr \\ x_\star & x_\star x_\star^\tr \end{smallmatrix}\right)$ where $x_\star$ solves~\cref{eq:QCQP}.

If the data matrices $P_0, \ldots, P_m$ are sparse, then the aggregate sparsity pattern $\mathcal{E}$ of the SDP~\cref{eq:QCQP_sdp} is also sparse, and the positive semidefinite constraint on $Z$ can be replaced with the conic constraint $Z \in \mathbb{S}^{n+1}_+(\mathcal{E},?)$. The chordal decomposition techniques described in \cref{section:sparse-SDPs} can therefore be applied to solve~\cref{eq:QCQP_sdp} efficiently. 
The following subsections briefly discuss two types of problem for which sparsity can be exploited effectively: Max-Cut problems \citep{goemans1995improved} and sensor network location problems \citep{kim2009exploiting,nie2009sum,so2007theory,jing2019angle}. 

\subsubsection{Max-Cut problem}
\noindent
The maximum cut (Max-Cut) problem is a classic problem in graph theory \citep{goemans1995improved}. %\red{[Does it have applications to control that can be pointed out?]}
Consider an undirected graph $\mathcal{G}(\mathcal{V},\mathcal{E})$ with $n$ vertices such that each edge $(i,j) \in \mathcal{E}$ is assigned a nonzero weight $W_{ij}$, and set $W_{ij} =0$ if $(i,j) \notin \mathcal{E}$. 
The Max-Cut problem aims to %find a cut of the graph with the largest possible weight, i.e., a 
partition the graph's vertices into two complementary sets $\mathcal{V}_1$ and $\mathcal{V}_2$ such that the total weight of all edges linking $\mathcal{V}_1$ and $\mathcal{V}_2$ is maximized. Given a binary variable $x \in \{-1, +1\}^n$ assigning nodes to one of the two partitions, one seeks to maximize
$$
    \frac{1}{2} \sum_{i,j: x_i x_j = -1} W_{ij} = \frac{1}{4} \sum_{i,j}W_{ij}(1 - x_ix_j).
$$
This is equivalent to solving
\begin{equation} \label{eq:maxcut}
    \begin{aligned}
        \min_{x} &\quad x^\tr W x \\
        \text{subject to} & \quad x_i^2 = 1, \quad i = 1, \ldots n,
    \end{aligned}
\end{equation}
where $W$ is the given matrix of weights.

This problem is a particular QCQP, and can easily be rewritten in the generic form~\cref{eq:QCQP} using data matrices $P_0, P_1, \ldots, P_n$ whose aggregate sparsity graph coincides with the original graph $\mathcal{G}$. If $\mathcal{G}$ is sparse with small maximal cliques,therefore, SDP relaxations of~\cref{eq:maxcut} can be solved efficiently using the sparsity-exploiting techniques in~\cref{section:sparse-SDPs}. Indeed, numerical experiments by \cite{andersen2010implementation} and \cite{ZFPGW2020chordal} demonstrated that the sparsity-exploiting solvers \software{SMCP} and \software{CDCS} can solve benchmark Max-Cut problems from the SDPLIB problem library \citep{borchers1999sdplib} order of magnitude faster than standard conic solvers.

\subsubsection{Sensor network location}
\noindent
The sensor network location problem, also known as \emph{Graph Realization} \citep{so2007theory}, has important applications such as inventory management and environment monitoring. At a basic level, the problem is to find unknown \emph{sensor points} $x_1, \ldots, x_n \in \mathbb{R}^d $ ($d=2$ or $3$) satisfying some specified distance constraints, as well as distance constraints with respect to $m$ known \emph{anchor points} $a_1, \ldots, a_m \in \mathbb{R}^d$. Precisely, given pairing sets
\begin{align*}
    \mathcal{E}_x \subseteq \{1,\ldots,n\} \times \{1,\ldots,n\},\\
    \mathcal{E}_a \subseteq \{1,\ldots,m\}\times \{1,\ldots,n\},
\end{align*}
we seek to find sensor locations $x_1, \ldots, x_n \in \mathbb{R}^d$ such that 
\begin{equation} \label{eq:sensor_constraints}
    \begin{aligned}
    \|x_i - x_j\|^2 = {d}^2_{ij}, \quad (i,j) \in \mathcal{E}_x, \\
    \|a_i - x_j\|^2 = f^2_{ij}, \quad (i,j) \in \mathcal{E}_a,
    \end{aligned}
\end{equation}
where the numbers $d_{ij}$ and $f_{ij}$ are specified distances.

One way to relax the sensor location problem into an SDP is to consider~\cref{eq:sensor_constraints} as a set of quadratic constraints for $x_1, \ldots, x_n$, and apply the generic SDP relaxation strategy to the QCQP \citep{kim2009exploiting} 
\begin{equation} \label{eq:sensor_network_problem}
\begin{aligned}
    \min_{x_1, \ldots, x_n \in \mathbb{R}^d} & \quad 0 \\
    \text{subject to} & \quad \cref{eq:sensor_constraints}.
\end{aligned}
\end{equation}
%Upon denoting $x = \begin{bmatrix} x_1^\tr, \ldots, x_n^\tr  \end{bmatrix}^\tr \in \mathbb{R}^{nd}$, and $X = xx^\tr$, all constraints in~\cref{eq:sensor_network_problem} becomes linear in the elements of $X$, and a semidefinite relaxation can be derived using a relaxed constraint $X \succeq xx^\tr$ \citep{kim2009exploiting}. 
It is clear that the data matrices and vectors of this QCQP are very sparse, and that the aggregate sparsity pattern of the corresponding SDP relaxation is determined only by the edge sets $\mathcal{E}_a$ and $\mathcal{E}_x$. Then, the techniques in~\cref{section:sparse-SDPs} can be applied to solve the relaxed problem quickly; we refer the interested reader to~\cite{kim2009exploiting} for more detailed discussions and experiment results. Similar ideas can be used to analyze sensor location problems where the distance measurements $d_{ij}$ and $f_{ij}$ are affected by noise \citep{kim2009exploiting}. 

\begin{remark}
There are other ways to formulate an SDP relaxation for~\cref{eq:sensor_network_problem}. One \citep{so2007theory} is to introduce a matrix variable $Y = XX^\tr$ with $X = \begin{bmatrix}x_1, x_2, \ldots, x_n \end{bmatrix} \in \mathbb{R}^{d \times n}$, rewrite all the constraints in~\cref{eq:sensor_network_problem} as linear equalities in $X$ and $Y$, relax the nonconvex relation between these variables into the inequality $Y \succeq XX^\tr$ and apply Schur's complement to obtain an SDP. %; see~\cite{so2007theory}, where the exactness of this relaxation is also discussed. 
A sparsity-exploiting version of this approach is described by \citet[Section 3.3]{kim2009exploiting}. 
Another option \citep{nie2009sum} is to formulate the search for the sensor locations as an unconstrained polynomial optimization problem,
$$
    \min_{x_1, \ldots, x_n} \sum_{(i,j) \in \mathcal{E}_a} (\|a_i - x_j\|^2 - f^2_{ij})^2 + \sum_{(i,j) \in \mathcal{E}_x}(\|x_i - x_j\|^2 - {d}^2_{ij})^2.
$$
The polynomial objective is term-sparse when the coupling set $\mathcal{E}_x$ contains only a small subset of all pairs $(i,j)$ (in fact, correlatively sparse; see \cref{s:pop-sparse-global} for definitions of these concepts). Therefore, the sparse SOS techniques outlined in \cref{s:pop-sparse-global} can be applied to solve the problem efficiently. The interested reader is referred to~\cite{nie2009sum} for experiment results. \markendexample
\end{remark}

\subsection{Machine learning: Verification of neural networks} \label{subsection:verification}

\noindent
Neural networks are one of the fundamental building blocks of modern machine-learning methods. For safety-critical applications, it is essential to ensure that they are provably robust to input perturbations. Given a neural network $f(x_0):\mathbb{R}^d \rightarrow \mathbb{R}^m$, a nominal input $\bar x \in \mathbb R^d$, a linear function $\phi:\mathbb{R}^m \to \mathbb{R}$ on the network's output, and a perturbation radius $\epsilon \in \mathbb{R}$, the network verification problem \citep{raghunathan2018semidefinite,salman2019convex,tjandraatmadja2020convex} asks to either verify that
\begin{equation}\label{eq:NNverification}
	\phi(f(x_0)) > 0 \quad \forall x_0:\, \| x_0 -	\bar{x} \|_\infty \leq \epsilon,
\end{equation} 
or to identify at least one counterexample to this relation. 

Consider an $L$-layer feedforward neural network where
\begin{align*}
f(x_0) &= W_Lx_L+b_L,\\
x_{i+1} &= \RELU(W_{i}x_i + b_i), \quad i = 0, \ldots, L-1,
\end{align*}
where $W_i \in \mathbb{R}^{n_{i+1} \times n_i}$ and $b_i \in\mathbb{R}^{n_{i+1}}$ are the network weights and biases, respectively, and the so-called Rectified Linear Unit (ReLU) activation function $\RELU:\mathbb{R}^k \to \mathbb{R}^k$ is the element-wise positive part of its argument, $\RELU(z) = [\max(z_i,0)]_{i=1}^k$. 
Condition~\cref{eq:NNverification} can be decided by solving the optimization problem
\begin{subequations} \label{eq:NNoptimization}
\begin{align}
	\gamma^\star := \min_{x_0, \ldots, x_L} &\quad c^\tr x_L + c_0
	\nonumber \\
	\text{subject to}& \quad x_{i+1} = \RELU(W_{i}x_i + b_i), \, i \in
	[L], \label{eq:NNoptimization_ReLU} \\
	&\quad  \|x_0 - \bar{x}\|_\infty \leq \epsilon
	\label{eq:NNoptimization_input},
\end{align}
\end{subequations}
where $[L]:=\{0, 1, \ldots, L-1\}$ and $c$, $c_0$ are problem data related to the linear function $\phi(\cdot)$. If $\gamma^\star >0$, then~\cref{eq:NNverification} holds, otherwise counterexamples can be found.  

Since the action of the ReLU function can be described by quadratic constraints,
\begin{align*}
y = \RELU(z) \quad \iff\quad
y \geq z, \;
y \geq 0, \;
y(y-z) = 0,
\end{align*}
problem~\cref{eq:NNoptimization} can be reformulated into a QCQP with variable $x = \begin{bmatrix} x_0^\tr, x_1^\tr, \ldots, x_L^\tr \end{bmatrix}^\tr$ \citep{raghunathan2018semidefinite}, and subsequently relaxed into an SDP as described in \cref{subsection:QCQP} above. 
If the optimal value of this SDP is positive, the network is verified; otherwise, nothing can be said.

\begin{figure}
	\centering
	
	\setlength{\abovecaptionskip}{0mm}
	\setlength{\belowcaptionskip}{0mm}
	\centering
	\subfigure[]{\includegraphics[width=0.9\linewidth]{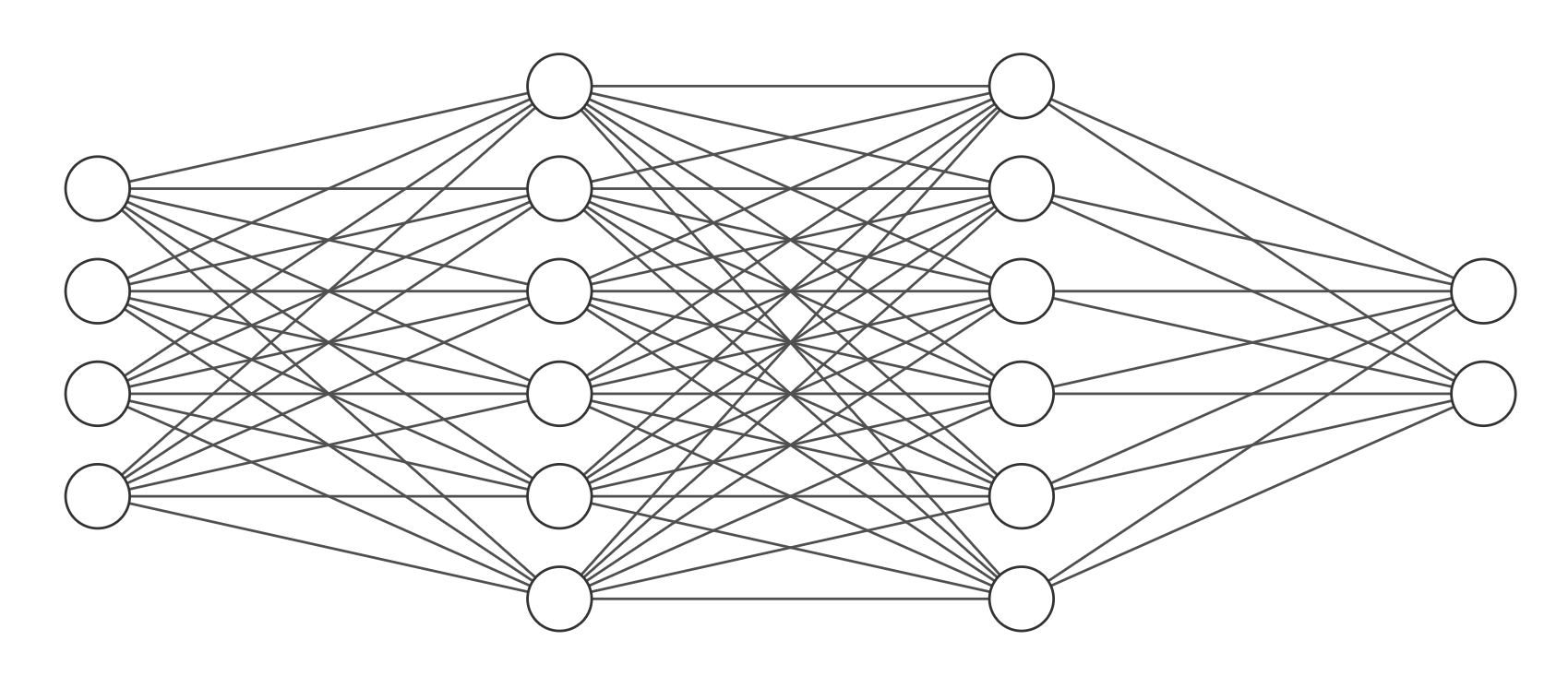} \label{figure:chordal_network}}
	\subfigure[]{  
		\begin{tikzpicture}[scale=1]
		{\small
		\node[draw, circle, inner sep=0pt, minimum size=0.8cm] (1) at (0,0) {\footnotesize $i-1$};
		\node[draw, circle, inner sep=0pt, minimum size=0.8cm] (2) at (2,0) {\footnotesize $i$};
		\node[draw, circle, inner sep=0pt, minimum size=0.8cm] (3) at (4,0) {\footnotesize $i+1$};
		\node[draw, circle, inner sep=0pt, minimum size=0.8cm] (4) at (6,0) {\footnotesize $i+2$};
		\draw[black, thick]  (1)--(2);
		\draw[black, thick]  (2)--(3);
		\draw[black, thick]  (3)--(4);
		}
		\end{tikzpicture}
		\label{figure:chordal_fourcliques}
	}
	\caption{Abstraction of \textit{(a)} a 4-layer neural network into \textit{(b)} a chordal chain graph with four vertices and maximal cliques $\{i-1,i\}$,  $\{i,i+1\}$ and $\{i+1,i+2\}$.}
	\label{figures:chordal_nn}
\end{figure}
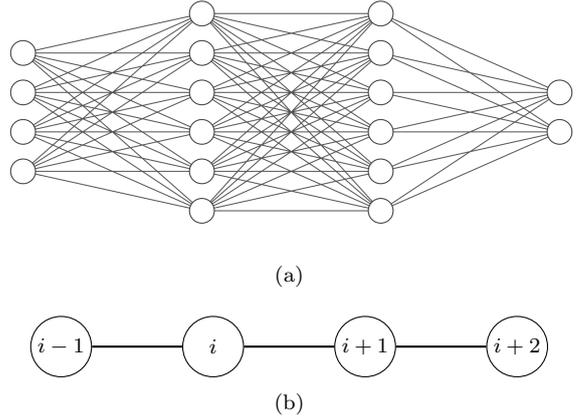

\begin{figure}[t]
  \centering
  \includegraphics[width=0.85\linewidth]{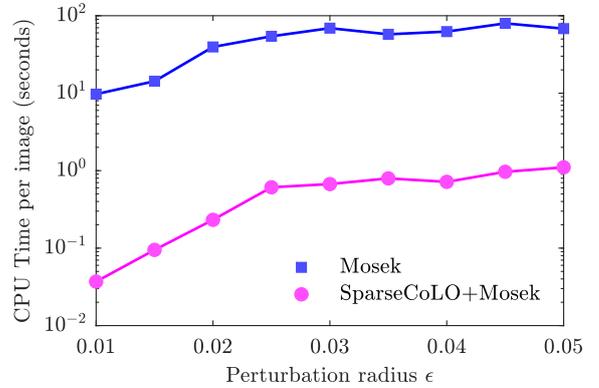}\\
  \caption{CPU time (seconds) required to solve SDP relaxations of the neural network verification problem~\cref{eq:NNoptimization} for image verification, with and without sparsity exploitation. The SDP solver was \software{MOSEK}, and sparsity was exploited using \software{SparseCoLO}. The neural network, with $L=2$ layers and $n_i = 64$ neurons per layer, was trained for image classification on the MNIST dataset.
  }
  \label{fig:neural_network_sparsecolo}
\end{figure} 

Since the constraints~\cref{eq:NNoptimization_ReLU} have a very natural cascading structure, the interaction among variables $x_0,\ldots, x_L$ can be modeled by a line graph with maximal cliques $\mathcal{C}_i=\{i,i+1\}$ for $i = 0, \ldots, L-1$ (see~\cref{figures:chordal_nn} for illustration with $L=4$). The SDP relaxation of~\cref{eq:NNoptimization} inherits this cascading structure, in addition to any sparsity coming from the structure of the weight matrices $W_i$. The chordal decomposition techniques described in~\cref{section:sparse-SDPs} can therefore be applied to solve it efficiently. This idea has been recently validated by~\cite{BPLZ2021neural}, who considered robustness verification in the context of image classifiers. {For instance, the results reproduced in \cref{fig:neural_network_sparsecolo} for a neural network with $L=2$ layers and $n_i = 64$ neurons per layer show that exploiting sparsity reduced by two orders of magnitude the CPU time required to verify the robustness of an image classifier on the MNIST dataset.}
Similar results were obtained by~\cite{NewP21}, and interested readers are invited to consult~\cref{tab:applications} for references to more machine learning applications where sparsity exploitation can dramatically reduce computational complexity. 

%%%%%%%%%%%%%%%%%%%%%%%%%%%%%%%%%%%%%%%%%%%%%%%%%%%%
\section{Conclusion and outlook}
\label{section:conclusion}
\noindent
In this paper, we reviewed theory and applications of decomposition methods for large-scale semidefinite and polynomial optimization. Specifically, we presented classical chordal decomposition results for sparse positive semidefinite matrices (cf.~\cref{T:ChordalDecompositionTheorem,T:ChordalCompletionTheorem,T:GeneralDecompositionTheorem,T:GeneralCompletionTheorem}) and we discussed how they can be exploited to implement efficient first- and second-order algorithms for SDPs (\cref{section:sparse-SDPs}). We showed also how matrix decomposition (primarily, but not necessarily, chordal) can be leveraged to exploit \textit{term sparsity} and \textit{structural sparsity} in large-scale polynomial optimization (\cref{section:polynomial_optimization}). In particular, we demonstrated that many sparsity-exploiting techniques for polynomial inequalties---including the well-known correlatively sparse SOS representations and the recent TSSOS, CS-TSSOS and chordal-TSSOS hierarchies---are based on the general matrix decomposition strategy outlined in~\Cref{s:sparse-sos-general-strategy}. We also discussed how the classical chordal decomposition theorem (\cref{T:ChordalDecompositionTheorem}) can be generalized in different ways to obtain SOS chordal decomposition theorems for sparse polynomial matrices (cf.~\cref{th:matrix-chordal-sos-decomp-global,,th:sos-matrix-decomp-even,,th:sos-matrix-decomp-local} and further results by \cite{zheng2020sum}). In~\Cref{section:factor-width-two}, we reviewed factor-width decompositions for SDPs with dense semidefinite constraints, to which chordal decomposition cannot be applied. Finally, in~\Cref{section:applications} we demonstrated how all of these techniques can be used to reduce the computational complexity of SDPs and polynomial optimization problems encountered in some control and machine learning applications. References to these and other applications are summarized in \cref{tab:applications}.

Despite the considerable progress made in recent years, numerical methods for semidefinite and polynomial optimization are still far from being mature. The most pressing open challenge, in our opinion, lies in bridging the gap between the size of SDPs that can currently be solved with tractable computational resources, and the size of the SDPs that arise from complex control applications. Indeed, the state-of-the-art decomposition techniques reviewed in this article are often still not enough to enable the use of semidefinite programming to analyze and control large-scale nonlinear systems. The same is true for control problems with systems of smaller size, but which require real-time computations. 

Achieving significant progress is likely to require theoretical extensions of the decomposition approaches we have discussed, as well as the development of efficient software that can effectively exploit modern multi-core and distributed-memory computer architectures. We conclude this article by outlining some possible research directions that may bear fruit in the near future.

\subsection*{Combining matrix decomposition with other structures}
\noindent
SDPs encountered in applications often have structural properties beyond sparsity, which can also be leveraged to reduce computational complexity; examples are symmetries, the existence of low-rank solutions, and low-rank data matrices~\citep{de2010exploiting,gatermann2004symmetry,majumdar2020recent}. %There are also extensive results in exploiting symmetry and low-rank properties to improve the scalability of solving SDPs~\citep{de2010exploiting,gatermann2004symmetry,majumdar2020recent}. 
It is natural to try and combine the exploitation of such additional structure with matrix decomposition, but, to the best of our knowledge, a unified and theoretically robust framework to do so is yet to be developed. Particular questions to be answered in this context include whether there exist symmetry reduction techniques that preserve (or even promote) sparsity in SDPs, and whether low-rank positive semidefinite completions \cite[Theorem 1.5]{dancis1992positive} can be exploited in SDPs with aggregate sparsity and low-rank optimal solutions (see \citealp{jiang2017minimum} and \citealp{miller2019chordal} for some results in this direction).

In addition, although we have presented chordal and factor-width decompositions separately, they can be combined if either one, applied in isolation, does not reduce the complexity of a large-scale SDP enough. 
A relatively straightforward approach \citep{miller2019decomposed} is to first apply the standard chordal decomposition, and then enforce positive semidefinite constraints associated to large maximal cliques using factor-width approximations. This idea can be taken forward in various directions; for instance, one could use block-chordal and block-factor-width decompositions, or extend ideas by \cite{gartska2020merging} to formulate adaptive strategies wherein cliques are either combined or factor-width decomposed, depending on their relative sizes and on the available computational resources. Both ideas remain largely unexplored, and further work is required to determine if they can be brought to bear on real-life control problems.

\subsection*{Tailored hierarchies for sparse polynomial optimization}
\noindent
Almost all existing methods for exploiting term sparsity in polynomial optimization rely on the general matrix decomposition approach presented in \Cref{s:sparse-sos-general-strategy}, where the Gram matrix associated with SOS certificates of nonnegativity is decomposed according to the maximal cliques of a sparsity graph to be prescribed \textit{a priori}. While the correlatively sparse, TSSOS, 
and related hierarchies described in \cref{s:pop-sparse-global} give useful general strategies to select this sparsity graph, there is ample scope for tailoring the graph structure in particular control applications. It is not unreasonable to expect that problem-specific choices, motivated for example by physical intuition on the dynamical system one is trying to analyse or control, may bring significant further gains. However, it remains to be seen whether this expectation can be met in practice. Better integration between the development of optimization tools and application-related modeling, discussed further below, seems key to achieving progress in this direction.

\subsection*{Decomposition and completion of polynomial matrices}
\noindent
The exploitation of sparsity for polynomial matrix inequalities can be improved in various directions, reducing computational complexity beyond what can be achieved using only the SOS chordal decomposition results summarized in \cref{s:pop-sparse-matrix}. For instance, those results can be combined in a natural way with techniques to leverage term-sparsity in scalar polynomial inequalities. Indeed, when a polynomial matrix inequality $P(x) \succeq 0$ is ``scalarized'' into a nonnegativity condition for the polynomial $p(x,y)=y^\tr P(x) y$, the structural sparsity of $P$ translates into correlative sparsity of $p$ with respect to $y$. The matrix decomposition results of \Cref{s:pop-sparse-matrix} have equivalent statement at the scalar level \cite[Section 4]{zheng2020sum} that can be used to refine or extend term-sparse SOS decomposition hierarchies for polynomials. The latter, in turn, can be used to efficiently handle (scalarized) polynomial matrix inequalities.

It would also be interesting to establish SOS completion results for sparse polynomial matrices, in the spirit of \Cref{T:ChordalCompletionTheorem}. Preliminary results in this direction exist \citep{zheng2018decomposition}, but are far from complete. Extension of the results in this reference will contribute to building a comprehensive theory for SOS chordal decomposition and completion of polynomial matrices, which can be used to build tractable SDP approximations of large-scale optimization problems with sparse polynomial matrix inequalities.

\subsection*{To chordality and beyond}
\noindent
Exploiting sparsity in semidefinite and polynomial optimization without modifying the problem usually requires chordality (cf. \cref{T:ChordalDecompositionTheorem,,T:ChordalCompletionTheorem,,T:GeneralDecompositionTheorem,,T:GeneralCompletionTheorem} for SDPs, and \cref{th:csp-putinar-scalar,th:sos-matrix-decomp-even,th:sos-matrix-decomp-local} for polynomial optimization). Enforcing chordality with traditional chordal extension strategies, even if approximately minimal, may lead to graphs with unacceptably large maximal cliques. The largest maximal clique size plays a major role in determining the computational complexity of a decomposed SDP (or SDP relaxation of a polynomial optimization problem). Therefore, systematic techniques to produce chordal extensions that approximately minimize the largest maximal cliques size would be very valuable.

If good chordal extensions prove hard to find, a compelling alternative is to sacrifice chordality and use nonchordal graphs with small cliques that can be determined analytically. This was done, for instance, by \cite{nie2009sparse} and \cite{kovcvara2020decomposition}. 
While clique decompositions of matrix inequalities based on nonchordal graphs are conservative in general, it may still be possible to identify classes of matrices for which the equivalence between the original and decomposed inequalities can be guaranteed.
For example, sparse (scaled)-diagonally dominant matrices always admit a clique decomposition, even when their sparsity graph is not chordal \cite[Proposition 1]{miller2019decomposed}. The same is true for certain positive semidefinite matrices whose sparsity pattern can be extended to be of a ``block-arrow'' type \citep{kovcvara2020decomposition}. Necessary and sufficient cycle conditions for positive semidefinite completion problem with nonchordal sparsity graphs were investigated by \cite{barrett1996real}. Extensions of these results, even if limited to particular application domains, are likely to enable considerable progress in the solution of large-scale SDPs with nonchordal sparsity.

\subsection*{Efficient software for modern computers}
\noindent
Reliable and user-friendly implementations of the cutting-edge decomposition techniques for SDPs and polynomial optimization problems reviewed in this paper are, in our opinion, just as important as further theoretical advances. Most of the available open-source packages mentioned in \Cref{subsection:implementations_SDP,s:pop-software} have not yet reached the level of maturity required to solve robustly a wide range of SDPs or polynomial optimization problems arising from real-life applications. Moreover, many of the commonly-used optimization modeling environments on which these packages rely are by now over a decade old, and often cannot handle extremely large problems of industrial relevance efficiently. 

The lack of very-high-performance software currently limits the scale of problems that can be solved without ad-hoc implementations. Since such implementations require considerable expertise in large-scale optimization, the deployment of SDP-based frameworks for system analysis and control to real-world problems is currently hindered. We expect that improvements in software reliability, efficiency, user-friendliness, and the ability to leverage modern multi-processor and/or distributed computing platforms will considerably increase the practical impact of decomposition methods for SDPs, bringing great benefit to the community of application-oriented researchers.

\subsection*{Blending application-driven modeling with optimization}
\noindent
The decomposition techniques reviewed in~\Cref{section:sparse-SDPs,,section:polynomial_optimization,,section:factor-width-two} apply to generic standard-form SDPs and polynomial optimization problems, irrespective of the context in which they arise. In control-related application, however, SDPs and polynomial optimization problems often come from modeling or relaxation frameworks for the study of dynamical systems, the details of which strongly affect the structure of the eventual optimization problem. Bridging the existing gaps between application-driven modeling and the development of large-scale optimization algorithm promises to enable significant progress in the study of linear and nonlinear systems.  On the one hand, it may be possible to implement tailored SDP solvers that target special structures arising in particular applications. On the other hand, given a particular control or analysis task, one should attempt to formulate modelling approaches that lead to optimization problems with a ``computationally friendly" structure. For example, when studying fluid flows using semidefinite programming (see, e.g., \citealp{Fantuzzi2018beanard-marangoni} and \citealp{Arslan2021}), a smart discretization of the flow field leads to SDPs with chordal aggregate sparsity that can be solved in minutes even though their linear matrix inequalities have more than $10\,000$ rows/columns. Similarly, using structured Lyapunov (or Lyapunov-like) functions as explained in \Cref{subsection:control_applications} can lead to structured SDPs, enabling the analysis of increasingly large systems in fields such as robotics, smart energy grid, and autonomous transportation.

Of course, the design of analysis and control frameworks that combine system-level modeling with algorithmic considerations will present a number of challenges. %, which may prove hard to overcome but also present numerous opportunities for interesting research. More importantly, however,
Resolving these challenges, however, promises to remove long-standing barriers to the study of complex systems, especially nonlinear ones. Success seems likely to require a collaborative effort between researchers working in different areas and an increasing awareness of outstanding problems in particular application domains, as well as of state-of-the-art tools for large-scale optimization. We hope that the present review of decomposition methods for semidefinite and polynomial optimization takes a step in the right direction and can inspire new discoveries in the near future.

%%%%%%%%%%%%%%%%%%%%%%%%%%%%%%%%%%%%%%%%%%%%%%%%%%%%
\section*{Acknowledgements}
\noindent
Y.Z was supported in part by Clarendon Scholarship. G.F. gratefully acknowledges funding from an Imperial College Research Fellowship. A.P. was supported in part by the Engineering and Physical Sciences Research Council (EPSRC) under project EP/M002454/1.

%%%%%%%%%%%%%%%%%%%%%%%%%%%%%%%%%%%%%%%%%%%%%%%%%%%%
\appendix
% \appendices
%%%%%%%%%%%%%%%%%%%%%%%%%%%%%%%%%%%%%%%%%%%%%%%%%%%%%%%%%%%%%%%%%%%%
%\section{Perfect elimination ordering and the zero fill-in property}
\section{Cholesky factorization with no fill-in}
\label{appendix:zero_fillin}
\begin{figure}[t]
  \centering
  \newcommand{\fighspace}{\hspace{0.5cm}}
  \subfigure[]
  { \label{Fig:ChordalDecomposition_a}
    \raisebox{2mm}{\includegraphics[height=0.25\columnwidth]{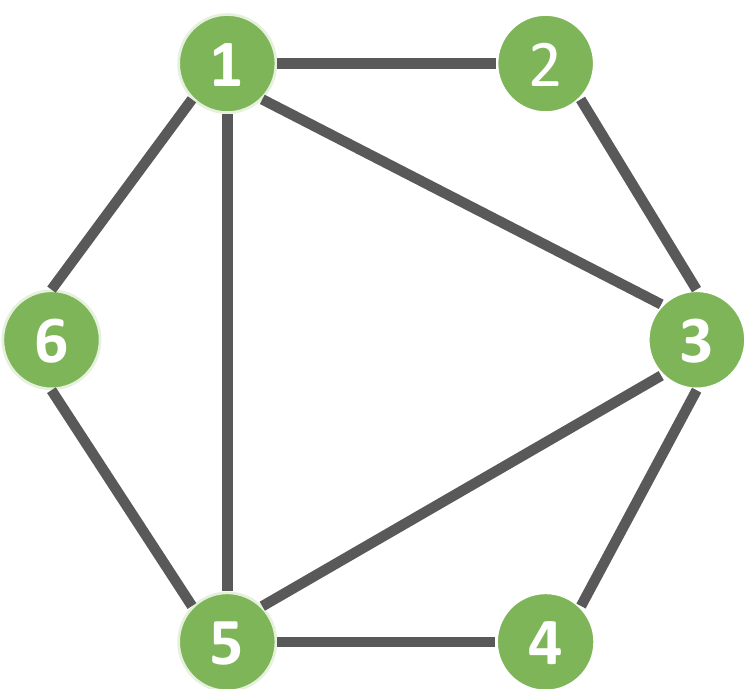}}
  }
  \fighspace
  \subfigure[]
  { \label{Fig:ChordalDecomposition_b}
    \includegraphics[height=0.3\columnwidth]{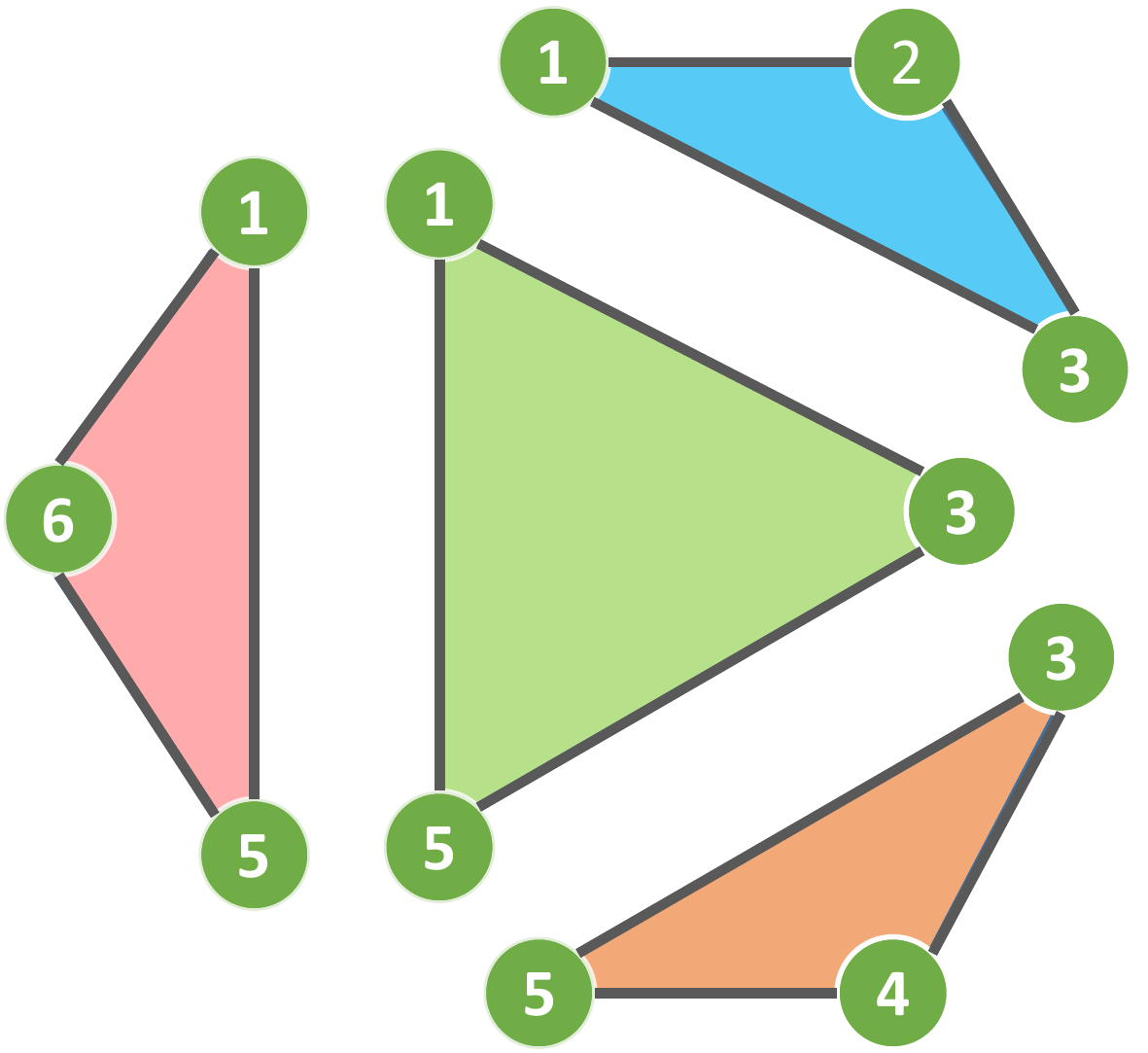}
  }
  \fighspace
  \subfigure[]
  { \label{Fig:ChordalDecomposition_c}
  \begin{tikzpicture}
    %\tikzset{>=latex}
    \footnotesize
    % place nodes
    \node[draw] at (0, 0)   (a) {$\mathcal{C}_1 = \{1,2,3\}$};
    %\node[draw] at (3, 0)   (b) {B};
    \node[draw] at (0, -1)  (b)     {$\mathcal{C}_4 = \{1,3,5\}$};
    \node[draw] at (-1.5, -2)  (c)     {$\mathcal{C}_3 = \{1,5,6\}$};
    \node[draw] at (1.5, -2)  (d)     {$\mathcal{C}_2 = \{3,4,5\}$};

    % draw edges
    \draw[] (a) -- (b);
    \draw[] (b) -- (c);
    \draw[] (b) -- (d);
\end{tikzpicture}
  }
  \caption{Chordal graph decomposition: (a) a chordal graph with six nodes; (b) maximal cliques; (c) a clique tree that satisfies the clique intersection property. %A chordal graph can be equivalently represented by a set of complete subgraphs.
  }
  \label{Fig:ChordalDecomposition}
\end{figure}
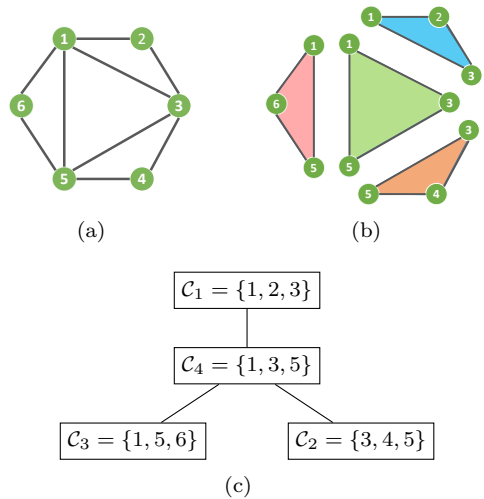

\noindent
The no fill-in property of the Cholesky factorization for positive definite matrices with chordal sparsity is one of the most important results for sparsity exploitation in matrix calculations; for instance, it enables a simple proof of~\cref{T:ChordalDecompositionTheorem} and efficient computations involving barrier functions for sparse matrix cones (cf.~\Cref{subsection:nonsymmetric_algorithms}). To formally introduce this no fill-in property, we first define the notions of simplicial vertices and perfect elimination ordering for graphs.
\begin{definition}
A vertex $v$  in a graph $\mathcal{G}(\mathcal{V},\mathcal{E})$ is called \emph{simplicial} if all its neighbors are connected to each other.
\end{definition}
\begin{definition}
An ordering $\sigma = \{v_1, \ldots, v_n\}$ of the vertices in a graph $\mathcal{G}$ is a \emph{perfect elimination ordering} if each $v_i, i = 1, \ldots, n,$ is a simplicial vertex in the subgraph induced by nodes $\{v_{i}, v_{i+1}, \ldots , v_n\}$.
\end{definition}
\noindent
For example, vertices $2, 4, 6$ are simplicial for the graph in~\cref{Fig:ChordalDecomposition_a}, and the ordering $\sigma = \{2,4,6,1,3,5\}$ is a perfect elimination ordering. 
A graph $\mathcal{G}$ is chordal if and only if it has at least one perfect elimination ordering~\citep[Theorem 4.1]{vandenberghe2015chordal}. The maximal cardinality search (\Cref{algorithm:maximal_cardinality_search}) either returns one of the perfect elimination orderings or certifies that none exists in $\mathcal{O}(\mathcal{|\mathcal{V}| + |\mathcal{E}}|)$ time \citep{tarjan1984simple}. %The algorithm assigns each node in the graph with a number that defines its place in the elimination ordering $\alpha$. It initializes all the nodes with zero weight $w \in \mathbb{R}^n$, starts from an arbitrary node, and generate the elimination ordering $\alpha$ in reverse. The weight vector $w$ is updated by increasing the weights of the nodes adjacent to the node that was previously labeled at each iteration. For a chordal graph,~\Cref{algorithm:maximal_cardinality_search} is guaranteed to return a perfect elimination ordering~\citep{tarjan1984simple}. If the resulting ordering is not a perfect elimination ordering, the graph is not chordal. Thus, this algorithm can also efficiently check whether a graph is chordal.

\begin{algorithm}[t]
  \caption{Maximal cardinality search}
  \label{algorithm:maximal_cardinality_search}
  \begin{algorithmic}
    \Require A graph $\mathcal{G}(\mathcal{V},\mathcal{E})$
    \Ensure  An elimination ordering $\alpha$ of $\mathcal{G}$
    \For{all vertices $v$ in $\mathcal{G}$}
       \State $w(v) = 0$.
    \EndFor
    
    \For{$i=n$ to $1$}
        \State pick an unnumbered vertex $v$ with maximum weight in $w$;
        \State set $\alpha(v) = i$;
        \For{all unnumbered vertex $u$ adjacent to $v$}
            \State $w(u) \leftarrow w(u) + 1$;
        \EndFor
    \EndFor
  \end{algorithmic}
\end{algorithm}

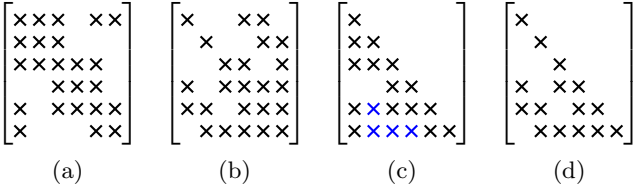
\begin{figure}[t]%[!b]
    \centering
    \setlength{\abovecaptionskip}{0.5em}
    \setlength{\belowcaptionskip}{0em}
    \newcommand{\fighspace}{\hspace{0cm}}
    %\subfigure[]{
    \begingroup % keep the change local
    \setlength\arraycolsep{1pt}
    \def\arraystretch{0.7}
    \begin{tikzpicture}
    \node at (0,0) {$\left[\begin{matrix}\\[-1.2ex]
    	\mycross{black} & \mycross{black} & \mycross{black} &  & \mycross{black} & \mycross{black} \\
    	 \mycross{black} & \mycross{black} & \mycross{black} & & &    \\ 
    	 \mycross{black} & \mycross{black} &\mycross{black} & \mycross{black} & \mycross{black} &  \\ 
    	 &  &\mycross{black} &\mycross{black} & \mycross{black} &      \\ 
    	\mycross{black} &  &\mycross{black} &\mycross{black} &\mycross{black} & \mycross{black}  \\
    	 \mycross{black}& & & &\mycross{black} &\mycross{black} 
    	\end{matrix}\right]$};
    \node at (0,-1.3) {\small (a)};
    	
    	\node at (2.2,0) {$\left[\begin{matrix}\\[-1.2ex]
    	\mycross{black} &  &  & \mycross{black} & \mycross{black} &  \\
    	  & \mycross{black} &  & &\mycross{black} & \mycross{black}   \\ 
    	  &  &\mycross{black} & \mycross{black} &  & \mycross{black} \\ 
    	\mycross{black} &  &\mycross{black} &\mycross{black} & \mycross{black} & \mycross{black}     \\ 
    	\mycross{black} & \mycross{black} & &\mycross{black} &\mycross{black} & \mycross{black}  \\
    	 &\mycross{black} &\mycross{black}&\mycross{black} &\mycross{black} &\mycross{black} 
    	\end{matrix}\right]$};
    \node at (2.2,-1.3) {\small (b)};
    
    \node at (4.4,0) {$\left[\begin{matrix}\\[-1.2ex]
    	\mycross{black} &  &  &  & &  \\
    	\mycross{black}  & \mycross{black} &  & &  &    \\ 
    	 \mycross{black} & \mycross{black} &\mycross{black} &  &  &  \\ 
    	 &  &  \mycross{black} & \mycross{black} & &    \\ 
    	\mycross{black} & \mycross{blue} & \mycross{black} &\mycross{black} &\mycross{black} &  \\
    	\mycross{black} &\mycross{blue} &\mycross{blue}&\mycross{blue} &\mycross{black} &\mycross{black} 
    	\end{matrix}\right]$};
    \node at (4.4,-1.3) {\small (c)};
    	\node at (6.6,0) {$\left[\begin{matrix}\\[-1.2ex]
    	\mycross{black} &  &  &  &  &  \\
    	  & \mycross{black} &  & & &    \\ 
    	  &  &\mycross{black} &  &  & \\ 
    	\mycross{black} &  &\mycross{black} &\mycross{black} &  &      \\ 
    	\mycross{black} & \mycross{black} & &\mycross{black} &\mycross{black} &   \\
    	 &\mycross{black} &\mycross{black}&\mycross{black} &\mycross{black} &\mycross{black} 
    	\end{matrix}\right]$};
    \node at (6.6,-1.3) {\small (d)};
    \end{tikzpicture}
    \endgroup
    \caption{(a) A symbolic $6 \times 6$ sparse positive definite matrix $Z$ with chordal sparsity graph shown in \cref{Fig:ChordalDecomposition_a}. (b) Sparsity pattern of $P_{\sigma}ZP_{\sigma}^\tr$ for the perfect elimination ordering $\sigma=\{2,4,6,1,3,5\}$. (c) Cholesky factor of $Z$; the entries marked by $\mycross{blue}$ denote nonzero fill-ins. (d) Cholesky factor of $P_{\sigma}ZP_{\sigma}^\tr$.
    %before/after performing : (a) pattern of $Z$; (b) pattern of ; (c) pattern of the Cholesky factor of $Z$; (d) pattern of the Cholesky factor of $P_{\sigma}ZP_{\sigma}^\tr$.
    }
    \label{Fig:ExampleCholesky}
\end{figure}

Now, given a positive definite matrix $Z \in \mathbb{S}^n_{+}(\mathcal{E},0)$ with a chordal sparsity pattern $\mathcal{E}$, we have a sparse Cholesky factorization with zero fill-in \citep{rose1970triangulated}, \citep[Theorem 9.1]{vandenberghe2015chordal}
\begin{equation} \label{eq:app_sparseCholesky}
    P_{\sigma}Z P_{\sigma}^\tr = L L^\tr, \qquad P_{\sigma}^\tr (L + L^\tr) P_{\sigma} \in \mathbb{S}^n(\mathcal{E},0),
\end{equation}
where $P_{\sigma}$ is a permutation matrix corresponding to the perfect elimination ordering $\sigma$ and $L$ is a lower-triangular matrix. This can be proven using an elimination process according to the perfect elimination ordering $\sigma$; see~\citep[Chapter 9.1]{vandenberghe2015chordal} and~\cite{kakimura2010direct} for details.  \cref{Fig:ExampleCholesky} illustrates the process of sparse Cholesky factorization for a %generic 
$6 \times 6$ positive definite matrix with chordal sparsity graph shown in \cref{Fig:ChordalDecomposition_a}.
% the positive definite matrix
% \begin{equation} \label{eq:app_counter_ex2}
%     Z = \begin{bmatrix}4  &  2  &  2 &    0&    1&    1\\
%      2  &  4    &2&    0&    0&    0\\
%      2   & 2   & 3&    2   & 2   & 0\\
%      0   & 0&    2&    3&    2 &   0\\
%      1    &0   & 2&    2   & 3  &  1\\
%      1    &0&    0&    0&    1&    3\end{bmatrix},
% \end{equation}
% whose chordal sparsity graph is shown in \cref{Fig:ChordalDecomposition_a}.

% \begin{figure}[t]%[!b]
%     \centering
%     \setlength{\abovecaptionskip}{0em}
%     \setlength{\belowcaptionskip}{0em}
%     \newcommand{\fighspace}{\hspace{1cm}}
%     \subfigure[]
%     { %\label{Fig:Ch1ExampleChordala}
%       \includegraphics[scale=0.8]{figures/figuressp_s1.eps}
%     } \fighspace
%     \subfigure[]
%     { %\label{Fig:Ch1ExampleChordalb}
%       \includegraphics[scale=0.8]{figures/figuressp_s2.eps}
%     } \fighspace
%     \subfigure[]
%     { %\label{Fig:Ch1ExampleChordalc}
%       \includegraphics[scale=0.8]{figures/figuressp_s3.eps}
%     }
%     \fighspace
%     \subfigure[]
%     { %\label{Fig:Ch1ExampleChordalc}
%       \includegraphics[scale=0.8]{figures/figuressp_s4.eps}
%     }
%     \caption{Sparsity pattern of a positive definite matrix $A \in \mathbb{S}^n_{+}(\mathcal{E},0)$ with a chordal pattern and its Cholesky factor before/after performing a perfect elimination ordering permutation: (a) pattern of $A$; (b) pattern of $P_{\sigma}AP_{\sigma}^\tr$; (c) pattern of the Cholesky factor of $A$; (d) pattern of the Cholesky factor of $P_{\sigma}AP_{\sigma}^\tr$. %\red{[GF: are we sure about (d)? The matrix $P_{\sigma}AP_{\sigma}^\tr$ should be symmetric.]}
%     %\red{[GF: would be better to have a smaller matrix here? It's difficult to read the figures]}
%     }
%     \label{Fig:ExampleCholesky}
% \end{figure}
%%%%%%%%%%%%%%%%%%%%%%%%%%%%%%%%%%%%%%%%%%%%%%%%%%%%
\section{A proof of \texorpdfstring{\cref{T:ChordalDecompositionTheorem}}{Theorem \ref{T:ChordalDecompositionTheorem}}} \label{appendix:proof}
\noindent
The sparse Cholesky factorization~\cref{eq:app_sparseCholesky} with zero fill-in allows for a simple proof of \cref{T:ChordalDecompositionTheorem}. For simplicity, but without loss of generality, assume that the matrix $Z$ has already been permuted in such a way that $\sigma = \{1,2,\ldots, n\}$ is a perfect elimination ordering, so $P_\sigma = I$ in~\cref{eq:app_sparseCholesky}. We denote the columns of $L$ by $l_1, l_2, \ldots, l_n$, and write 
$$
    Z = L L^\tr = \sum_{i=1}^n l_i l_i^\tr. 
$$
Since $L + L^\tr$ has the same sparsity pattern $\mathcal{E}$, the non-zero elements of each column vector $l_i$ must be indexed by a maximal clique $\mathcal{C}_{h_i}$ for some $h_i \in \{1, \ldots, t\}$. Thus, the non-zero elements of $l_i$ can be extracted through multiplication by the matrix $E_{\mathcal{C}_i}$, and we have
$$
    l_i = E^\tr_{\mathcal{C}_{h_i}} E_{\mathcal{C}_{h_i}} l_i \quad \Rightarrow \quad l_il_i^\tr = E^\tr_{\mathcal{C}_{h_i}} \underbrace{\left(E_{\mathcal{C}_{h_i}} l_il_i^\tr E_{\mathcal{C}_{h_i}}^\tr \right)}_{Q_{i}}E_{\mathcal{C}_{h_i}}.
$$
Now, let $J_{k} = \{i: h_i = k\}$ be the set of column indices $i$ such that column $i$ is indexed by clique $\mathcal{C}_k$. These index sets are disjoint and $\cup_k J_{k} = \{1,\ldots,n\}$, so we obtain
\begin{equation*}
\begin{aligned}
    Z = LL^\tr 
    &= \sum_{i =1}^n E_{\mathcal{C}_{h_i}}^\tr Q_{i} E_{\mathcal{C}_{h_i}}\\
    &= \sum_{k=1}^t \sum_{i \in J_{k}} E_{\mathcal{C}_k}^\tr Q_{i} E_{\mathcal{C}_k}\\
    &= \sum_{k=1}^t E_{\mathcal{C}_k}^\tr \bigg( \sum_{i \in J_{k}} Q_{i} \bigg) E_{\mathcal{C}_k}.
\end{aligned}
\end{equation*}
This is exactly~\cref{eq:ChordalMatrixDecomposition} in \cref{T:ChordalDecompositionTheorem} with matrices $Z_{k} = \sum_{i \in J_{k}} Q_{i}$ that is in $\mathbb{S}^{|\mathcal{C}_k|}_+$.

%%%%%%%%%%%%%%%%%%%%%%%%%%%%%%%%%%%%%%%%%%%%%%%%%%%%
\section{Some properties of maximal cliques}
\label{app:maximal_cliques}
\noindent
A connected chordal graph $\mathcal{G}(\mathcal{V},\mathcal{E})$ with $n$ vertices has at most $n-1$ maximal cliques that can be identified in linear time---more precisely, with a complexity of $\mathcal{O}(|\mathcal{V}| + |\mathcal{E}|)$ \citep{tarjan1984simple,berry2004maximum}. \Cref{algorithm:maximal_clique} is a simple strategy with a complexity $\mathcal{O}(|\mathcal{V}| + |\mathcal{E}|)$  to find all maximal cliques based on a perfect elimination ordering.
% For example, consider the chordal graph in~\cref{Fig:ChordalDecomposition_a}. A perfect elimination ordering is $\sigma = \{2,4,6,1,3,5\}$. If we run \Cref{algorithm:maximal_clique} on this graph with this elimination ordering, the node sets are 
% \begin{equation*}
%     \begin{aligned}
%         \mathcal{C}_1 &= \{2,1,3\} \qquad \leftarrow \text{maximal clique}\\
%         \mathcal{C}_2 &= \{4,3,5\}  \qquad \leftarrow \text{maximal clique} \\
%         \mathcal{C}_3 &= \{6,5,1\}  \qquad \leftarrow \text{maximal clique} \\
%         \mathcal{C}_4 &= \{1,3,5\}  \qquad \leftarrow \text{maximal clique} \\
%         \mathcal{C}_5 &= \{3,5\} \quad \qquad \leftarrow \text{subset of $\mathcal{C}_4$, not maximal clique} \\
%         \mathcal{C}_6 &= \{5\} \qquad \qquad \leftarrow \text{subset of $\mathcal{C}_4$, not maximal clique} \\
%     \end{aligned}
% \end{equation*}
For example, the chordal graph in~\cref{Fig:ChordalDecomposition_a} has the perfect elimination ordering $\sigma = \{2,4,6,1,3,5\}$, and \Cref{algorithm:maximal_clique} constructs the sets
\begin{equation*}
    \begin{aligned}
        \mathcal{C}_1 &= \{2,1,3\}, &% \qquad \leftarrow \text{maximal clique}\\
        \mathcal{C}_2 &= \{4,3,5\}, &%  \qquad \leftarrow \text{maximal clique} \\
        \mathcal{C}_3 &= \{6,5,1\},\\%  \qquad \leftarrow \text{maximal clique} \\
        \mathcal{C}_4 &= \{1,3,5\}, &%  \qquad \leftarrow \text{maximal clique} \\
        \mathcal{C}_5 &= \{3,5\}, &% \quad \qquad \leftarrow \text{subset of $\mathcal{C}_4$, not maximal clique} \\
        \mathcal{C}_6 &= \{5\}. % \qquad \qquad \leftarrow \text{subset of $\mathcal{C}_4$, not maximal clique} \\
    \end{aligned}
\end{equation*}
The sets $\mathcal{C}_1,\ldots,\mathcal{C}_4$ are maximal cliques, while $\mathcal{C}_5, \mathcal{C}_6$ are not because they are subsets of $\mathcal{C}_4$.

The maximal cliques of a chordal graph can be arranged in a so-called \textit{clique tree}, that is, a graph $ \mathcal{T}(\Gamma, \Xi) $ with the maximal cliques $\Gamma = \{\mathcal{C}_1, \ldots, \mathcal{C}_t\}$ as its vertices and an edge set $ \Xi \subseteq \Gamma \times \Gamma $. %such that $(\mathcal{C}_i,\mathcal{C}_j) \in \Xi$ if and only if $\mathcal{C}_i \cap \mathcal{C}_i \neq \emptyset$. 
In particular, the clique tree can be chosen to satisfy the \emph{clique intersection property}, meaning that $ \mathcal{C}_i \cap \mathcal{C}_j \subseteq \mathcal{C}_k $ if clique $ \mathcal{C}_k $ lies on the path between cliques $ \mathcal{C}_i $ and $ \mathcal{C}_j $ in the tree and the intersection $\mathcal{C}_i \cap \mathcal{C}_j$ is nonempty \citep{blair1993introduction}. For example, the clique tree in \cref{Fig:ChordalDecomposition_c} satisfies the clique intersection property.  

%as well as finding a clique tree that satisfies the clique intersection property

\begin{algorithm}[t]
  \caption{Maximal clique search}
  \label{algorithm:maximal_clique}
  \begin{algorithmic}
    \Require A chordal graph $ \mathcal{G}(\mathcal{V}, \mathcal{E})$, and a perfect elimination ordering $ \alpha = \{v_1, \ldots,v_n\}$
    \Ensure  All its maximal cliques $\mathcal{C}_1,\mathcal{C}_2, \ldots, \mathcal{C}_t$
    \State Initialize $\mathcal{C}_0 = \emptyset$;
    \For{$i=1$ to $n$}
        \State $\mathcal{C}_i = \{ v_i\} \cup \{ u \text{ adjacent to }v_i \text{ and behind } v_i\text{ in } \alpha \}$;
        \If {$\mathcal{C}_i$ is not a subset of $\mathcal{C}_{0}$}
            \State $\mathcal{C}_i$ is a maximal clique;
            \State $ \mathcal{C}_0=\mathcal{C}_i$;
        \EndIf
    \EndFor
  \end{algorithmic}
\end{algorithm}

The maximal cliques of a chordal graph play a central role in the sparse matrix decomposition results stated in \cref{T:ChordalDecompositionTheorem,,T:ChordalCompletionTheorem,,T:GeneralDecompositionTheorem,,T:GeneralCompletionTheorem}. It is important to remember that these require one to use \textit{all} maximal cliques in the (chordal) sparsity graph of a matrix $X$, even when a subset of cliques already covers all nonzero entries of $X$. For example, consider the indefinite matrix
$$
    X = \begin{bmatrix} 
     2 &    2&    2&      0 &    1 &    1 \\
    2 &    2&    2&       0&      0 &  0     \\
    2 &   2&    2&    2&    2&       0 \\
       0    &   0&    2&    2&    2&      0\\
    1&       0 &    2&    2&    2&    1\\
    1&       0  &0 &     0&    1&    2
    \end{bmatrix},
$$
whose chordal sparsity graph is shown in \cref{Fig:ChordalDecomposition} and has the four maximal cliques identified above. Even though the maximal cliques $\mathcal{C}_1$, $\mathcal{C}_2$, and $\mathcal{C}_3$ already cover all nonzero entries of the matrix, the maximal clique $\mathcal{C}_4$ is necessary when applying \cref{T:ChordalCompletionTheorem} to check whether $X$ admits a positive semidefinite completion. Indeed, observing that
$$
    E_{\mathcal{C}_1}XE_{\mathcal{C}_1}^\tr = E_{\mathcal{C}_2}XE_{\mathcal{C}_2}^\tr = \begin{bmatrix} 2&2&2\\
                    2&2&2\\
                    2&2&2 \end{bmatrix} \in \mathbb{S}^3_+
$$
and 
$$
E_{\mathcal{C}_3}XE_{\mathcal{C}_3}^\tr=\begin{bmatrix} 2&1&1\\
                    1&2&1\\
                    1&1&2 \end{bmatrix} \in \mathbb{S}^3_+,
$$
is not sufficient to conclude $X \in \mathbb{S}^6_+(\mathcal{E},?)$ because the submatrix
$$
E_{\mathcal{C}_4}XE_{\mathcal{C}_4}^\tr=\begin{bmatrix} 2&2&1\\
                    2&2&2\\
                    1&2&2 \end{bmatrix}
$$
indexed by clique $\mathcal{C}_4$ has one negative eigenvalue. Similarly, the matrix 
\begin{equation*} %\label{eq:app_counter_ex2}
    Z = \begin{bmatrix}4  &  2  &  2 &    0&    1&    1\\
     2  &  4    &2&    0&    0&    0\\
     2   & 2   & 3&    2   & 2   & 0\\
     0   & 0&    2&    3&    2 &   0\\
     1    &0   & 2&    2   & 3  &  1\\
     1    &0&    0&    0&    1&    3\end{bmatrix}
\end{equation*}
%the matrix $Z$ in~\cref{eq:app_counter_ex2} 
is positive semidefinite and has the same sparsity graph as above, but it does not admit a decomposition $$Z = \sum_{k=1}^3 E_{\mathcal{C}_k}^\tr Z_k E_{\mathcal{C}_k}, \qquad Z_k \succeq 0$$ that uses only cliques $\mathcal{C}_1$, $\mathcal{C}_2$, and $\mathcal{C}_3$; the last maximal clique $\mathcal{C}_4$ is necessary for \cref{T:ChordalDecompositionTheorem} to apply. Indeed, any decomposition using only the first three maximal cliques requires
\begin{align*}
    Z_1 &= \begin{pmatrix}\alpha & 2 & 2\\ 2 & 4 & 2 \\ 2 & 2 & \beta\end{pmatrix}, \\
    Z_2 &= \begin{pmatrix}3-\beta & 2 & 2\\ 2 & 3 & 2 \\ 2 & 2 & \gamma\end{pmatrix},\\
    Z_3 &= \begin{pmatrix}4-\alpha & 1 & 1\\ 1 & 3-\gamma & 1 \\ 1 & 1 & 3\end{pmatrix},
\end{align*}
where $\alpha$, $\beta$ and $\gamma$ must be selected to make these three matrices positive semidefinite. For this, it is necessary that the diagonal elements and all $2\times2$ principal minors of $Z_1$, $Z_2$ and $Z_3$ be nonnegative; in particular,
\begin{subequations}
\begin{gather}
    \begin{align}
    \alpha,\beta,\gamma &\geq 0, &
    3-\beta &\geq 0, & 
    4-\alpha &\geq 0,\\
    4\alpha -4&\geq 0, &
    \alpha\beta -4 & \geq 0, &
    (3-\beta) \gamma -4 & \geq 0, %\label{eq:inequalties_A1}
    \end{align}
    \\
    (4-\alpha)(3-\gamma) -1  \geq 0. \label{eq:inequalties_A2}
\end{gather}
\end{subequations}
However, this set of inequalities is infeasible. Specifically, inequality~\cref{eq:inequalties_A2} can be rearranged to show that
\begin{equation*}
    \gamma \leq \frac{11 - 3\alpha}{4-\alpha}.
\end{equation*}
Moreover, we must have $3-4/\alpha \geq 3-\beta \geq 0$, so $\alpha \geq 4/3$, and therefore
\begin{equation*}
        \frac{4\alpha}{3\alpha-4} \leq \frac{4}{3-\beta} \leq \gamma \leq \frac{11 - 3\alpha}{4-\alpha}.
\end{equation*}
But this cannot be true because $4/3 \leq \alpha \leq 4$, so $\frac{4\alpha}{3\alpha-4} > \frac{11 - 3\alpha}{4-\alpha}$ strictly.

\balance
%%%%%%%%%%%%%%%%%%%%%%%%%%%%%%%%%%%%%%%%%%%%%%%%%%%%
\bibliography{references}
\end{document}